\documentclass[11pt]{amsart}
\usepackage{amsmath, amssymb}
 \usepackage[colorlinks=true,linkcolor=blue,urlcolor=blue,citecolor=blue, pagebackref]{hyperref}
\usepackage[alphabetic]{amsrefs}
\usepackage{amsfonts, epsfig,color,wrapfig,epstopdf}
\usepackage{ xypic,verbatim,amscd,color}
\usepackage{mathrsfs}
\usepackage{youngtab, ytableau}
\usepackage{graphicx}
\usepackage{tikz-cd}
\usepackage{colordvi}
\usepackage{relsize}
\usepackage[bbgreekl]{mathbbol}
\usepackage{amsfonts}
\DeclareSymbolFontAlphabet{\mathbb}{AMSb} 
\DeclareSymbolFontAlphabet{\mathbbl}{bbold}

\numberwithin{equation}{section}

\newcommand\op{\operatorname}
\newcommand\bull{{\bullet}}

\newcommand\frg{\mathfrak{g}}
\newcommand\frh{\mathfrak{h}}

\newcommand\mf{\mathcal{F}}

\newcommand\mk{\mathcal{K}}

\newcommand\Pic{\operatorname{Pic}}

\newcommand\me{\mathcal{E}}

\newcommand\mv{\mathcal{V}}
\newcommand\tensor{\otimes}
\newcommand\ml{\mathcal{L}}

\newcommand\codim{\operatorname{codim}}

\newcommand\mg{\mathcal{G}}

\newcommand\rk{\operatorname{rk}}

\newcommand\mn{\mathcal{N}}

\newcommand\SL{\operatorname{SL}}

\newcommand{\leto}[1]{\stackrel{#1}{\to}}

\newcommand\pic{\operatorname{Pic}}
\newcommand\Bun{\operatorname{Bun}}

\newcommand\frl{\mathfrak{l}}

\newcommand\Parbun{\op{Parbun}}

\newcommand{\C}{\Bbb{C}}
\newcommand{\Q}{\Bbb{Q}}
\renewcommand{\P}{\mathbb{P}}
\newcommand{\Z}{\mathbb{Z}}

\newtheorem{theorem}{Theorem}[section]

\newtheorem{remark}[theorem]{ Remark}

\newtheorem{corollary}[theorem]{Corollary}

\newtheorem{proposition}[theorem]{Proposition}

\newtheorem{lemma}[theorem]{Lemma}

\newtheorem{definition}[theorem]{Definition}
\newtheorem{definition/lemma}[theorem]{Definition/Lemma}
\newtheorem{defi}[theorem]{Definition}

\begin{document}
\title[The multiplicative eigenvalue problem]{Vertices in multiplicative eigenvalue problem for arbitrary groups}
\author{Prakash Belkale and Joshua Kiers}
\maketitle
\begin{abstract}
We determine, in an inductive framework, the vertices of the polytope $P(s,K)$ controlling the conjugacy classes of elements  which product to one in the maximal compact subgroup $K$ of a simple complex algebraic group $G$. This extends earlier work of the authors in related contexts. One feature of this work is
the use of Kontsevich compactifications of the moduli of $P$-bundles (replacing the use of quot schemes in type A) which are related to semi-infinite geometry. We also obtain a quantum generalization of Fulton's conjecture valid for all $G$.
\end{abstract}

\section{Introduction}

The well-known Hermitian eigenvalue problem asks for the possible eigenvalues of a sum of Hermitian matrices $A+B$ given the eigenvalues of $A$ and $B$. There is a natural unitary analogue of this problem, which asks for the possible eigenvalues of a product of unitary matrices $UV$ given the eigenvalues of $U$ and $V$. Both of these problems can be generalized to a larger number of factors. In its Lie theoretic formulation, the unitary problem concerns conjugacy classes in the maximal compact subgroup $K = SU(n)$ of $G = SL(n)$. The set of possible eigenvalues forms a compact, convex polytope $P(s,K)$ inside the $s$-fold product of the Lie algebra of the (real) maximal torus.

The problem admits a natural generalization to any simply-connected, simple Lie group $G$. The polytopes $P(s,K)$ are important in the study of genus-zero conformal blocks: specifically, the nonvanishing spaces of conformal blocks, under saturation, are controlled by the rational points inside $P(s,K)$ (see Theorem \ref{killing} below).

The main result of this work is an explicit method for naming the vertices of $P(s,K)$ for general $G$. This extends the work of the first author in the case $G = SL(n)$ \cite{BRigid}, and broadly parallels the authors' prior work in a number of related contexts (\cite{BRays,BKiers,KBranch}). Some of the vertices of $P(s,K)$ in the type A case were related in \cite{BRigid}, via strange duality, to the theory of rigid local systems \cite{Katz}. Strange duality phenomena exist for a larger collection of groups, (see \cite[Appendix D]{KumarConformalBlocks}) and it would be interesting to see whether these type A relations generalise. Other potential applications include the study of quantum saturation for arbitrary groups.

\subsection{A note on our methods}
The paper \cite{BRigid} used the quot scheme compatification of the moduli stack of $P$-bundles where $P$ is a maximal parabolic subgroup of $G=SL(n)$. It is  well-known that the ``far-alcove" wall in type A can be treated by diagram automorphisms, or geometrically by passing to shifted bundles. These operations played an essential role in \cite{BRigid}, and some of the results initially seemed special to type A, including the formula for the degrees of the basic divisors from enumerative geometry (especially Prop 6.0.1 in loc. cit.). With some effort we were able to formulate general $G$ versions of the results, and the resulting geometry was considerably different.

For general $G$ there are two compactifications of the moduli stack of $P$-bundles on a curve  that have been well studied: These are the Drinfeld and Kontsevich compactifications (see \cite{BrGa}, and \cite{FFKM,Campbell}). The former is not smooth and is desingularized by the latter \cite{Givental,FFKM,Campbell}. The geometry of the Drinfeld compactification is connected to the semi-infinite geometry of flag varieties (\cite{FFKM, Kato}, also Remark \ref{LBO} below). Since some of our methods require analysis of tangent spaces, we have chosen to work throughout with Kontsevich compactifications. We have been influenced by the combinatorics of Lusztig's Bruhat order  and semi-infinite lengths \cite{Lusztig,LamS} which arise in semi-infinite geometry. One of our proofs (the proof of Theorem \ref{T2}) involves the deformation theory of a nodal curve in a finite dimensional context, but is related to semi-infinite geometry.

We develop induction techniques from \cite{BKq} for making the induction maps explict. The general outline of the description of extremal rays follows the classical cases from \cite{BRays,BKiers,KBranch}, as well as tangent space techniques from \cite{BKR}.

\subsection{Facets of $P(s,K)$}

Let $B$ denote a fixed Borel subgroup of $G$ and $H$ a fixed maximal torus inside $B$. The Lie algebras of $G$, $B$, and $H$ will be denoted by $\mathfrak{g}$, $\mathfrak{b}$, $\mathfrak{h}$. Let $R\subset \mathfrak{h}^*$ be the set of roots of $\mathfrak{g}$ with respect to $\mathfrak{h}$, and $R^+$ the set of positive roots with respect to $\mathfrak{b}$. Let $W = N_G(H)/H$ be the Weyl group, where $N_G(H)$ is the normalizer of $H$ in $G$.
The set of simple roots (minimal positive roots) is denoted $\Delta = \{\alpha_1,\hdots,\alpha_r\}$, and the corresponding simple reflections by $\{s_1,\dots,s_r\}\subset W$. For each $1\le j\le r$, we define the fundamental coweight $x_j\in \mathfrak{h}$ by the requirement
$$
\alpha_i(x_j) = \delta_{i,j}, ~\forall 1\le i\le \ell.
$$

As $G$ is simple, there is a unique highest root $\theta\in R^+$. The \emph{fundamental alcove} $\mathscr{A}\subset \mathfrak{h}$ is
$$
\mathscr{A} = \{\mu\in \mathfrak{h} : \alpha_i(\mu)\ge 0 ~ \forall i, \theta(\mu)\le 1\}.
$$
The conjugacy classes of $K$ are bijectively parametrized by $\mathscr{A}$ via
\begin{align*}
C:\mathscr{A}&\to K/\op{Ad}K\\
\mu&\mapsto \text{the class of }\op{Exp}(2\pi i \mu).
\end{align*}
Fix an integer $s\ge 3$. The \emph{multiplicative eigenpolytope} $P(s,K)$ is by definition
$$
P(s,K) = \{(\mu_1,\hdots,\mu_s)\in \mathscr{A}^s : 1\in C(\mu_1)\cdots C(\mu_s)\}.
$$

To name the inequalities defining $P(s,K)$ inside $\mathscr{A}^s$, we first describe the (deformed) small quantum cohomology ring associated to the homogeneous spaces $G/P$. 

Let $P\supseteq B$ be a standard parabolic subgroup of $G$ and $W_P$ the Weyl group of $P$, which is naturally a subgroup $W_P\subseteq W$. Let $W^P$ denote the set of minimal-length coset representatives for $W/W_P$. To each $w\in W^P$ is associated a Schubert variety $X^P_w:=\overline{BwP/P}\subset G/P$ of dimension $\ell(w)$. Let $\sigma_w^P$ denote the class in $H^{2(\dim G/P-\ell(w))}(G/P,\Z)$ which is Poincar\'e-dual to the fundamental class of $X^P_w$ in $H_{2\ell(w)}(G/P,\Z)$.

Let $L=L_P$ be the Levi subgroup of $P$ containing $H$, and $\Delta_P$ the set of simple roots contained in $R_{\mathfrak{l}}$, where $\mathfrak{l}$ is the Lie algebra of $L$. Let $S_P=\Delta\setminus \Delta_P$.

\begin{defi}
Recall that $H_2(G/P)$ has a basis given by the fundamental classes $\mu(X^P_{s_i})$, where $\alpha_i\in S_P$. An element $d\in H_2(G/P)$ can be written as
$$d=(a_i)_{\alpha_i\in S_P}=\sum_{\alpha_i\in S_P} a_i\mu(X^P_{s_i})$$
We say $d\geq 0$ if all $a_i\geq 0.$
\end{defi}

Fix distinct points $p_1,\hdots,p_k\in \P^1$, general $(g_1,\hdots,g_k)\in G^k$ and $d\geq0 \in H_2(G/P)$. The Gromov-Witten invariant
$$
\langle \sigma_{u_1}^P, \hdots,\sigma_{u_k}^P \rangle_d,
$$
counts the number of maps $f:\P^1\to G/P$ of degree $d$ (see Def. \ref{stacky}) such that $f(p_i)\in g_iC_{u_i}^P$ (if infinitely many, the invariant is defined to be $0$).
\begin{remark}
We will simplify notation by denoting  $
\langle \sigma_{u_1}^P, \hdots,\sigma_{u_k}^P \rangle_d,
$ simply by $\langle {u_1}, \hdots,{u_k} \rangle_d$ when $P$ is clear from the context.
\end{remark}

Introduce variables $q_i$ for each $\alpha_i\in S_P$ and let $q^d=\prod_{\alpha_i\in S_P} q_i^{a_i}$
These invariants define a product in $H^*(G/P,\Z)\otimes \Z[q]$ by setting:
$$
\sigma_u^P*\sigma_v^P = \sum_{w\in W^P} \langle \sigma_u^P,\sigma_v^P,\sigma_{w_0 ww_0^P}^P\rangle_d ~q^d~ \sigma_w^P
$$
Moreover $*$ is associative (and clearly commutative) by \cite{FP,FW}, and graded, provided we assign $q_i$ the degree $2-2\rho^L(\alpha_i^\vee)$.

The first author and S. Kumar introduced a deformation of this product in \cite{BKq} by setting
\begin{align*}
\langle \sigma_{u_1}^P, \hdots, \sigma_{u_k}^P \rangle_d^{\circledast_0} &=
\left\{
\begin{array}{rr}
 & \text{ if $(u_1,\hdots,u_k,d)$ is}\\
 \langle \sigma_{u_1}^P, \hdots, \sigma_{u_k}^P \rangle_d,  & \text { ``quantum Levi-movable''}\\
    & \text{(see below)}\\\\
0, & \text{else};
\end{array}
\right. \\\\
&\text{and}\\\\
\sigma_u^P\circledast_0 \sigma_v^P &= \sum_{w\in W^P} \langle \sigma_u^P,\sigma_v^P,\sigma_{w_0 ww_0^P}^P\rangle_d^{\circledast_0} ~q^d~ \sigma_w^P.
\end{align*}

\begin{remark}
The condition of Levi-movability comes from geometry and we review it in Section \ref{LeviM2}. Note that Levi-movability is characterised by  $\langle \sigma_{u_1}^P, \hdots, \sigma_{u_k}^P \rangle_d\neq 0$  and a numerical condition (see \cite[Theorem 3.15]{BKq}): For all $\alpha_i\in S_P$,
$$(\chi_e-\sum_{j=1}^k \chi_{u_j})(x_i)+ \sum_{\alpha\in R^+-R_{\frl}^+}\alpha(x_i)\alpha(\tilde{d})=0$$
where $\tilde{d}=\sum_{\alpha_i\in S_P} a_i\alpha_i^{\vee}$ and $\chi_w=\rho-2\rho_L +w^{-1}\rho.$
\end{remark}

They further showed that these deformed invariants exactly parametrize the minimal inequalities (i.e., the facets) of the eigenpolytope $P(s,K)$. This result extended and refined the prior work of Biswas \cite{Bis}, Agnihotri-Woodward \cite{AW}, Belkale \cite{Blocal}, and Teleman-Woodward \cite{TW}.

\begin{theorem}[\cite{BKq}*{Theorem 1.3}]
Let $(\mu_1,\hdots,\mu_s)\in \mathscr{A}^s$. Then the following are equivalent:
\begin{enumerate}
\item $(\mu_1,\hdots,\mu_s)\in P(s,K).$
\item For any standard maximal parabolic subgroup $P$, any $u_1,\hdots, u_s\in W^P$, and any $d\ge 0$ such that
$$
\langle \sigma_{u_1}^P, \hdots, \sigma_{u_s}^P \rangle_d^{\circledast_0} = 1,
$$
the following inequality is satisfied:
$$
\sum_{k=1}^s \omega_P(u_k^{-1} \mu_k) \le d.
$$
\end{enumerate}
Moreover, none of these inequalities may be removed.
\end{theorem}

Given the data of a maximal parabolic $P$ and $(u_1,\hdots,u_s,d)$ such that $\langle \sigma_{u_1}^P, \hdots, \sigma_{u_s}^P \rangle_d^{\circledast_0} = 1$, we define the corresponding facet $\mathcal{F}(P,u_1,\hdots,u_s,d)$ ($\mathcal{F}$ for short) of $P(s,K)$:
$$
\mathcal{F}(P,u_1,\hdots,u_s,d) = \{(\mu_1,\hdots,\mu_s)\in P(s,K): \sum_{k=1}^s \omega_P(u_k^{-1} \mu_k) = d\}.
$$

Here are the main results of this paper. 
\begin{enumerate}
\item We provide a formula for systematically producing {\em some} vertices of $P(s,K)$ lying on $\mathcal{F}$, using the input data $(P,u_1,\hdots,u_s,d)$.
\item We also provide an explicit linear map $P(s,K_L)\to \mathcal{F}$ whose image contains the remaining vertices of $\mathcal{F}$ (those not covered by (1)); here $K_L$ is the maximal compact subgroup of the Levi subgroup $L$ associated with $P$.
\item In contrast to the type A case \cite[Lemma 12.4.1]{BRigid}, or the classical case of the (additive) Hermitian eigenvalue problem for  arbitrary groups \cite{BKiers},  not all vertices of $P(s,K)$ lie on the ``regular facets'' $\mf$ as above.
Examples are given in Section \ref{non-regular}. These extra vertices are necessarily vertices of $\mathscr{A}^s$, and can be determined by computation. Such extra extremal ray phenomena are already known to occur in the subgroup embedding problem \cite{KBranch}. Here we note that, as first observed by the second named author in a different context, one such ``extra'' vertex in the case of $D_4$ gives an example where quantum saturation fails (see Proposition \ref{non-sat}).
\end{enumerate}
These formulas are stated in Sections \ref{formulae} -- \ref{tII-intro}, in slightly different terms, after we describe $P(s,K)$ in terms of conformal blocks, since this is closer to the geometry of the rest of the paper.

\subsection{Conformal blocks}

Let $\op{Parbun}_G$ be the moduli-stack of quasi-parabolic principal $G$-bundles on $\P^1$ with marked points $\{p_1,\hdots,p_s\}$, i.e., parametrizing tuples $\tilde {\mathcal{E}} = (\mathcal{E},\bar g_1,\hdots, \bar g_s)$ such that $\mathcal{E}$ is a principal $G$-bundle on $\P^1$ and $ \bar g_i \in \mathcal{E}_{p_i}/B$.

The Picard group of $\op{Parbun}_G$ is identified as follows:
$$
\op{Pic}(\op{Parbun}_G)=
\prod_{i=1}^s X^*(T) \times \Z \mathcal{L},
$$
where $\mathcal{L}$ is the pull-back of the generator of $\op{Pic}(\op{Bun}_G)$ (see \cite{LS,Sorg}).  The line bundle associated to $(\lambda_1,\hdots,\lambda_s,\ell \mathcal{L})$ will be denoted $\mathcal{B}(\vec \lambda,\ell)$ for short.

Let $\kappa:\mathfrak{h}^*\to \mathfrak{h}$ denote the isomorphism induced by the Killing form, normalized by $(\theta,\theta)=2$ where $\theta$ is the highest root.

Let $\lambda_1,\hdots,\lambda_s$ be a collection of dominant weights, and $\ell\ge 0$ a nonnegative integer. Assume that the $\lambda_i$ ``are at level $\ell$'', i.e., $\lambda({\theta}^{\vee})\leq \ell$ for all $i$. This is equivalent to assuming that $\kappa(\lambda_i)/\ell \in  \mathscr{A}$ for all $i$.

The \emph{space of conformal blocks} is
$$
V(\vec \lambda,\ell) = H^0(\op{Parbun}_G,\mathcal{B}(\vec \lambda,\ell))^*.
$$
One of the main questions regarding these vector spaces is when they have positive dimension (as a function of the inputs $\vec \lambda,\ell$). The connection between this question and the eigenvalue problem is given by the following result (see \cite{BKq}*{Theorem 5.2}, where our $\ell$ is their $2 g^*m$).

\begin{theorem}\label{killing}
For dominant weights $\lambda_1,\hdots,\lambda_s$ and level $\ell> 0$, the following are equivalent:
\begin{enumerate}
\item There exists an integer $N\ge 1$ such that the conformal blocks vector space $V(N\vec \lambda,N\ell)$ has positive dimension.
\item $(\frac{\kappa(\lambda_1)}{\ell}, \frac{\kappa(\lambda_2)}{\ell}, \hdots, \frac{\kappa(\lambda_s)}{\ell})$ belongs to $P(s,K)$.
\end{enumerate}
\end{theorem}

Therefore vertices of $P(s,K)$ correspond bijectively to extremal rays of the cone
 $$
 \mathcal{C} =\mathcal{C}_G = \op{Pic}^+_\Q(\op{Parbun}_G) =  \{(\vec \lambda,\ell): V(N\vec \lambda, N\ell)\ne (0) \text{ for some $N\ge 1$}\}.
 $$

\subsection{Formulas for divisor classes}\label{formulae}

To produce extremal rays of $\mathcal{C}$, we will name certain $G$-invariant divisors in $\op{Parbun}_G$ whose associated line bundles $\mathcal{B}(\vec \lambda, \ell)$ generate extremal rays. These divisors all arise in a common fashion, and their corresponding parameters $(\vec \lambda,\ell)$ are obtained through enumerative (Gromov-Witten) invariants. We describe this process now.

Let $v_1,\dots,v_s\in W^P$ and $d'\in \Bbb{Z}=H_2(G/P)$. Assume
\begin{equation}\label{excess}
\sum_{i=1}^s \codim \sigma^P_{v_i} = \dim X +\int_{d'} c_1(T_X) +1.
\end{equation}

\begin{definition}\label{stacky2}
Let $\Omega^0(\vec{v},d')$ 
denote the set of tuples $(\mathcal{E},f,\overline{g_1},\hdots,\overline{g_s})$ where $\mathcal{E}$ is a principal $G$-bundle on $\mathbb{P}^1$ and $\overline{g_i}\in \mathcal{E}_{p_i}/B$ for each $i$ and  $f:\mathbb{P}^1\to \mathcal{E}/P$ satisfy
\begin{itemize}
\item for each $i$, $\overline{g_i}$ and $f(p_i)$ are in relative position $v_i$;
\item $f$ has degree $d'$.
\end{itemize}
The relative dimension of  $\Omega^0(\vec{v},d')$  over $\Parbun_G$ is $-1$.

We let $\overline{\widetilde{D}}$ $=\Omega(\vec{v},d')$ denote the closure of $\widetilde{D}$ in a suitable Kontsevich compactification  (see Section \ref{KCompact}) which is proper over
$\Parbun_G$, and let
$D$ denote the effective divisor class on $\Parbun_G$, obtained by pushing forward  $\overline{\widetilde{D}}$ to $\Parbun_G$.

\end{definition}

Set $\mathcal{O}(D) = \mathcal{B}(\vec\lambda,\ell)$. In order to write formulas for the $\lambda_i$, we first recall the notation for covering relations in the Bruhat order of $W$ from \cite{BGG}:
\begin{definition}
Let $v,w\in W$ and $\beta\in R^+$. If $w = s_\beta v $ and $\ell(w) = \ell(v)+1$, then we write $v\xrightarrow{\beta}w$. Note that if $v\xrightarrow{\beta}w$, then $w^{-1}\beta\in R^-$ (and $v^{-1}\beta \in R^+$).
\end{definition}

Then we have the following formulas for finding $\lambda_1, \hdots, \lambda_s$ and $\ell$. Write $\lambda_i = \sum c_i^b \omega_b$, where $\omega_b$ is the $b^\text{th}$ fundamental weight.
\begin{theorem}\label{formula11}
\begin{enumerate}
\item Set $v_{s+1} = \overline{s_\theta w_0}$ and $m' = \omega_P(\theta^\vee)$. Then
$$
\ell = \langle v_1,\hdots, v_s, v_{s+1}\rangle_{d'+m'}.
$$
\item $c_i^b=0$ if $\ell(s_{b}v_i)<\ell(v_i)$ or if $s_{b}v_i\not\in W^P$. Otherwise, $v_i\xrightarrow{\alpha_b}s_bv_i\in W^P$ and we define $w_k = v_k$ for $k\ne i$ and $w_i = s_bv_i$. Then
$$
c_i^b = \langle w_1,\hdots,w_s\rangle_{d'}.
$$
\end{enumerate}
\end{theorem}
Theorem \ref{formula11} is proved in Section \label{cycleclasses} (see Theorem \ref{firstformula}) in a generalised form applicable to arbitrary standard parabolics $P$.

\subsection{Type I rays}\label{tI-intro}

Fix the data of a face of $P(s,K)$: i.e., a maximal parabolic $P$, Weyl group elements $u_1,\hdots,u_s\in W^P$, and degree $d$ such that $\langle \sigma_{u_1}^P,\hdots,\sigma_{u_s}^P\rangle_d^{\circledast_0}=1$. According to Theorem \ref{killing}, these give a facet $\mathcal{F}$ of $\mathcal{C}$ defined by
$$
\mathcal{F}(P,u_1,\hdots,u_s,d) = \{(\lambda_1,\hdots,\lambda_s,\ell)\in \mathcal{C} : \sum_{i=1}^s \lambda_i(u_ix_P) = d\ell \frac{2}{(\alpha_P,\alpha_P)}\}.
$$

Due to $\langle \sigma_{u_1}^P,\hdots,\sigma_{u_s}^P\rangle_d^{\circledast_0}\ne 0$, it must be true that
$$
\sum_{i=1}^s \codim \sigma^P_{u_i} = \dim X +\int_{d} c_1(T_X).
$$
Fix $1\le j\le s$ and $v\in W^P$ such that either
\begin{enumerate}
\item[(A)] $v\xrightarrow{\beta} u_j$ for a simple root $\beta$, or
\item[(B)] $v = \overline{s_\theta u_j}\in W^P$ and $\ell(v)>\ell(u_j)$.
\end{enumerate}

Define $v_1,\hdots, v_s$ by setting $v_k = u_k$ for $k\ne j$ and $v_j = v$. Then
set $d'$ according to the cases above:
\begin{enumerate}
\item[(A)] $d' = d$, or
\item[(B)] $d' = d- m$ where $m = \omega_P(u_j^{-1}\theta^\vee)$ (see Lemma \ref{fraction}).
\end{enumerate}

Clearly in case (A), the collection $(v_1,\hdots, v_s, d')$ satisfy \eqref{excess}. It turns out that this also happens in case (B) (see Lemma \ref{fraction}(2)). So $(v_1,\hdots,v_s,d')$ produce a $G$-invariant divisor $D$ in $\op{Parbun}_G$ as in Definition \ref{stacky2}; let us denote this divisor by $D(v,j)$.  These divisors can be related to the semi-infinite Bruhat order (Remark \ref{LBO}).

\begin{theorem}\label{type1}
\begin{enumerate}
\item $D(v,j)$ is an irreducible divisor in $\op{Parbun}_{G}$.
\item $h^0(\op{Parbun}_{G},\mathcal{O}(m D(v,j))) = 1$ for all positive integers $m$.
\item $\mathbb{Q}_{\ge 0}\mathcal{O}(D(v,j))$ is an extremal ray of $\op{Pic}^+_{\mathbb{Q}}(\op{Parbun}_{G})$.
\item This extremal ray lies on the face $ \mathcal{F}$ (and is called a type I ray on $\mf$).
\end{enumerate}
\end{theorem}
Theorem \ref{type1} is proved in Section \ref{beginning}  (see Theorem \ref{type11}) in a generalised form applicable to arbitrary standard parabolics $P$.
As a corollary of our constructions, we obtain the quantum generalization of Fulton's conjecture (Theorem \ref{FulCon}) for arbitrary groups.

\subsection{Induction and type II rays}\label{tII-intro}
Let $D_1,\dots,D_q$ be the extremal ray generators produced as $[D(v,j)]$ in Theorem \ref{type1}. We call these type I rays. 
We will show (Corollary \ref{salvador}) that the face $\mf$ of $\mathcal{C}$ is a product
\begin{equation}\label{biject1}
\mf\leto{\sim} \prod_{i=1}^q {\Bbb{Q}_{\ge 0}} D_i \times \mf_{\op{II}}
\end{equation}
Now $\mf_{\op{II}}$ is related to the Levi subgroup by an induction apparatus. We show that it is surjected on to (Lemma \ref{pruf}) by  the effective cone of $\Parbun_L(d))$. We then show in Lemma \ref{twostepprocess} that this effective cone is isomorphic to the effective cone of $\Parbun_L(0))$, which is then related to the effective cone of 
$\Parbun_L'$ where  $L'=[L,L]$. $L'$ breaks up as a product of simple, simply-connected groups $L_k$, i.e., $L'=\prod L_k$.
We show that there is an explicit surjective map, and we thus obtain a surjection

\begin{equation}\label{inducto}
\prod\mathcal{C}_{L_k}\twoheadrightarrow \mf_{\op{II}}
\end{equation}
and thus extremal rays on $\mf_{\op{II}}$, called type II rays, are images of (some) extremal rays of the $\mathcal{C}_{L_k}$.

The induction algorithm is given in Section \ref{algorithm}, and examples are given in Section \ref{exemples}.

\begin{remark}\label{LBO}
Recall that type I rays (generated by $[D(v,j)]$ above) arise from two possible scenarios on a Weyl group element $u\in W^P$:
\begin{enumerate}
\item[(A)] $s_i u \xrightarrow{\alpha_i} u$ or

\item[(B)] $u < \overline{s_\theta u}$, where the bar indicates taking min-length representative in $W^P$.
\end{enumerate}
These two criteria can be framed in a unified way using the semi-infinite order on the affine Weyl group (see \cite[Section 2.2]{KNS}). The semi-infinite order is the transitive closure of the requirement that $s_{\beta + k\delta } w t_\xi <_{\frac{\infty}{2}} w t_\xi$ whenever the finite root $w^{-1}\beta \prec 0$. Clearly if finite Weyl group elements $v$ and $w$ satisfy $v\le w$, then $v <_{\frac{\infty}{2}}w$ also.

The affine roots for the Levi consist of all $\beta+k\delta$ such that $\beta$ is a finite Levi root: $\beta \in R_{\mathfrak{l}}$. The positive affine roots are those where $\beta\succ 0$ and $k\ge 0$ or $\beta \prec 0$ and $k>0$.

The distinguished representatives $(W^P)_{af}$ are those $x\in W_{af}$ such that $x(\gamma)\succ 0$ for every positive affine root $\gamma$ for the Levi. (These will not be ``minimum-length'' with respect to the semi-infinite length function, which can take on both positive and negative values.)

Possibility (B) above is equivalent to
$$
u^{-1} \theta \in R^+\setminus R_{\mathfrak{l}}^+.
$$

This is equivalent to the simultaneous conditions
$$
s_0u<_{\frac{\infty}{2}} u
$$
and
$$
s_0u \in (W^P)_{af}.
$$

In case (A) above, it is also immediate that $s_iu <_{\frac{\infty}{2}} u$ and $s_iu \in (W^P)_{af}$ since $W^P \subset (W^P)_{af}$. Therefore we obtain a type I ray datum whenever
$$
s_iu <_{\frac{\infty}{2}}u, s_iu\in (W^P)_{af}
$$
for $i=0, 1, \hdots, r$.
\end{remark}

\section{Some preliminaries}

\subsection{Parameter spaces and compactifications}\label{KCompact}
\begin{defi}
Let $\op{Parbun}_{G,S}$ be the stack parametrizing tuples $\overline{\mathcal{E}}=(\mathcal{E},\overline{g_1},\hdots,\overline{g_s})$ where $\mathcal{E}$ is a principal $G$-bundle on $\mathbb{P}^1$ and $\overline{g_i}\in \mathcal{E}_{p_i}/B$ for each $i$. Here $S = \{p_1,\hdots,p_s\}\subset \mathbb{P}^1$. We drop the subscript $S$ if it is clear from the context.
\end{defi}
Let $P\subseteq G$ be a standard parabolic subgroup, $w_1,\dots,w_s\in W^P$ and $d\in H_2(G/P)$.  Recall that a $P$-reduction of a principal bundle $\mathcal{E}$ is a section $F:\Bbb{P}^1\to\mathcal{E}/P$ of $\mathcal{E}/P\to\Bbb{P}^1$. 

\begin{definition}\label{stacky}
Let $\Omega^0(\vec{w},d)$ denote the set of tuples $(\mathcal{E},\overline{g_1},\hdots,\overline{g_s},F)$ where $\overline{\mathcal{E}}=(\mathcal{E},\overline{g_1},\hdots,\overline{g_s})\in \op{Parbun}_{G}$, and $F:\mathbb{P}^1\to \mathcal{E}/P$ is a section satisfying
\begin{itemize}
\item for each $i$, $\overline{g_i}$ and $F(p_i)$ are in relative position $w_i\in W^P$. This relative position is defined as follows: Pick a trivialization $e\in \mathcal{E}_{p_i}$ and write $F(p_i) = eh_iP$ and $\bar g_i = e\tilde g_iB$ for some $\tilde g_i\in G$. Then, $w_i$ is defined by the requirement $h_iP\in  \tilde g_iBw_i P\in G/P$. A different choice of $e$ acts on $h_i$ and $g_i$ by a left multiplication and therefore does not affect the relative position.
\item $F$ has degree $d=(a_i)_{\alpha_i\in S_P}$, 
meaning: the first Chern class $c_1(F^*(\mathcal{E}(\omega_i)))=a_i$ for each $\alpha_i\in S_P$. Here { $\mathcal{E}(\omega_i)=\mathcal{E}\times^P\Bbb{C}_{-\omega_i}$} is the line bundle on $\mathcal{E}/P$ associated to the fundamental weight $\omega_i:P\to \Bbb{C}^*$.
\end{itemize}
\end{definition}

There is a natural representable morphism $\Omega^0(\vec{w},d)\to \op{Parbun}_{G}$. The compactification to be described below parameterizes stable maps (with marked points) $C\to \mathcal{E}/P$. In the special case where $\mathcal{E}$ is trivial, this is the data of a map $C\to \Bbb{P}^1\times G/P$ so that the map to $\Bbb{P}^1$ has 
degree $1$: Therefore on the locus where $\mathcal{E}$ is trivial, the compactifcation is the graph space $\overline{M}_{0,n}(G/P\times\Bbb{P}^1,(d,1))$ (or rather, its fiber over a general point of $(\Bbb{P}^1)^n$ since $z_1,\dots,z_n$ are fixed).

Following \cite{Campbell,FFKM}, we introduce a relative Kontsevich compactification  $\Omega(\vec{w},d)\to \op{Parbun}_{G}$ of this morphism. A $Y$-point of $\Omega(\vec{w},d)$
is the following:
\begin{enumerate}
\item A morphism $p:C\to \Bbb{P}^1\times Y$ of degree $1$, such that  $C\to Y$ is a flat family of connected nodal projective curves of genus $0$.
\item For each $i$ a section $q_i:Y\to C$, such that the composite $Y\to C\to\Bbb{P}^1$ is the constant map to $p_i$.
\item A principal $G$-bundle $\mathcal{E}$ on $\Bbb{P}^1\times Y$.
\item A principal $P$-bundle structure on $C$ with a given isomorphism  to $p^*\me$ when the structure group is extended to $G$. This structure is equivalent to a
section $F:C\to \me/P$ over $\Bbb{P}^1\times Y$. Note that $ \me/P$ is a $G/P$ bundle over $\Bbb{P}^1\times Y$. The section $F$ must further satisfy: 
\begin{enumerate}
\item $F:(C,q_1,\dots,q_s)\to \me/P$ is Kontsevich stable over closed points $y\in Y$ (this implies in particular  that the points $q_i$ are on the smooth locus of $C$).
\item The degree of $F$ is $d$; i.e., the first Chern class $c_1(F^*(\mathcal{E}(-\omega_i)))=a_i$ for each $\alpha_i\in S_P$. Here $\mathcal{E}(\omega_i)=\mathcal{E}\times^P\Bbb{C}_{\omega_i}$ is the line bundle on $\mathcal{E}/P$ associated to the dominant integral weight $\omega_i$ as before.
\end{enumerate}
\item Sections (over $Y$) $\overline{g}_i$ of $(\mathcal{E}/B)\mid _{p_i\times Y}$ for each $i$, so that over all geometric points of $Y$, $F(q_i)$ and $\overline{g_i}$ are in relative position $\leq w_i$ (in the Bruhat order) for every $i$.
\end{enumerate}
Let $K(d)$ be the space above without the  conditions and choices (2) and (5). We claim that $K(d)$ is a smooth and irreducible stack, such that the 
fibers of $K(d)\to \Bun_G$ are proper Deligne-Mumford stacks. Recall (for motivation) first that the space of maps from $\Bbb{P}^1\to G/P$ of a given degree is  smooth and irreducible \cite{Thomsen,KP} and dense in the Kontsevich compactification. Note that the space of maps is a open subset of the space of $P$-bundles of a given degree, and is hence smooth and irreducible

  For arbitrary $P$, the methods of \cite{Campbell,Yang} show that with $\mathcal{M}_X$ as defined in \cite{Campbell}, $K(d)=\overline{Bun}^K_P$ is smooth over $\mathcal{M}_X$ (the proof of Proposition 2.4.1 uses that the stack of $B$-bundles on a possibly nodal curve is smooth. However this is the case for $B$ replaced by any algebraic group \cite{wang}.) Now  $\mathcal{M}_X$  is smooth by loc. cit. The irreducibility of $K(d)$ follows from the density of the open part where $C\to\Bbb{P}^1$ is an isomorphism (which is $\Bun_P(\Bbb{P}^1)$ plus choices of points on $\Bbb{P}^1$), see \cite[Proposition 2.6.3]{Yang}, and the connectedness of the stack  $\Bun_P(\Bbb{P}^1)$ formed by bundles of a  given degree (we reduce to the Levi subgroup). 

There is  a map $K(d)\to (\Bbb{P}^1)^n$. Let $K'(d)$ be the fibre over a general point $(p_1,\dots,p_s)\in (\Bbb{P}^1)^s$ (which incorporates the condition (2) above). By Kleiman transversality, $K'(d)$ is smooth, and has a dense open subset of points where the map $\mathcal{C}\to\Bbb{P}^1$ is an isomorphism. Clearly $\Omega(\vec{w},d)$ is a fibre bundle over $K'(d)$. 

By dimension counting over $K'(d)$, any component of $\Omega(\vec{w},d)\setminus \Omega^0(\vec{w},d)$ has dimension less than the dimension of $\Omega^0(\vec{w},d)$. Hence,
\begin{lemma}\label{relativity}
$\Omega(\vec{w},d)$ is irreducible of relative dimension (over $\Parbun_G$)
\begin{align}\label{rel-dim}
\dim G/P +\int_d c_1(T_X) -\sum_{i=1}^s \codim \sigma^P_{w_i}
\end{align}
\end{lemma}
The fibres of $\Omega(\vec{w},d)\to \op{Parbun}_{G}$ are proper Deligne-Mumford stacks, and moreover $K(d)$ is smooth as a stack. We will work for the most part over the automorphism free part of $\Omega(\vec{w},d)$ (e.g., when $F$ is an embedding). Locally, in the smooth topology of  $\op{Parbun}_{G}$, the map $\Omega(\vec{w},d)\to \op{Parbun}_{G}$  is a projective morphism of orbifolds.

\begin{remark}
The irreducibility of $\Omega^0(\vec{w},d)$ does not need the compactifications. We could avoid the assertions of irreducibilty and smoothness above for the most part by defining  $\Omega(\vec{w},d)$ to be the closure of
 $\Omega^0(\vec{w},d)$. This suffices for the proof of Theorem \ref{codim1boundary}. However for the proof of Theorem \ref{T2}, we need the fact that the curves that the (automorphism free maps) $\mathcal{C}\to \me/P$ that appear in the quantum stratum are smooth points of the moduli space of objects satisfying (1), (3)  and (4), and can be perturbed so that the map $C\to \Bbb{P}^1$ is an isomorphism. This assertion 
 is covered by \cite[Propositions 2.4.1, 4.4.1]{Campbell} (see Remark \ref{smoothening} for an example).

\end{remark}

\subsection{$T$-fixed curves in $G/P$}\label{TFIX}

The ``non-classical'' boundary of $\Omega(\vec{w},d)-\Omega^0(\vec{w},d)$ consists of points where the curve $C$ in Definition \ref{stacky} is reducible. It turns out that
in codimension one (see Theorem \ref{codim1boundary}), these curves have an extra component which maps to a $T$ fixed curve in a $G/P$. So we recall some facts about $T$-fixed curves in $G/P$ from \cite{FW}.

Let $P\supset B$ be a standard parabolic subgroup of $G$.
For any positive root $\alpha$ not in $R_{\frl}^+$, there is a unique $T$-invariant curve $C_{\alpha}$ in $X=G/P$ passing through 
$\dot{e}$ and 
$\dot{s_\alpha}$ given by  a mapping $\Bbb{P}^1\to X$ with $t\mapsto U_{-\alpha}(t)\dot{e}$ which takes $t=0$ to $\dot{e}$ and $t=\infty$ to $\dot{s_{\alpha}}$. Here $U_{-\alpha}(t)$ is the one parameter unipotent subgroup corresponding to the root $-\alpha$.

The degree of $C_{\alpha}$ is $d(\alpha)=\sum_{\alpha_i\in \Delta-\Delta_P}\omega_{i}(\alpha^{\vee})\mu(X^P_{s_{i}})$ and  by \cite[Lemma 3.5]{FW},
\begin{equation}\label{degdalpha}
\int_{d(\alpha)} c_1(T_X)=\bigl(\sum_{\gamma\in R^+\setminus R_{\frl}^+}\gamma\bigr)(\alpha^{\vee}).
\end{equation}
Any $T$-invariant curve passing through $\dot{e}$ is necessarily of this form. Note that $\int_{d(\alpha)}c_1(T_X)\geq 2$ (see \cite[Lemma 3.5]{FW}). Any $T$-invariant curve in $G/P$ is a $W$-translate of such a curve: it has fixed points $\dot{v}$ and $\dot{w}$ where $v=\overline{ws_{\alpha}}$ with $w\in W^P$ for some $\alpha \not \in R_{\frl}^+$. The degree of such a curve  $C(v,w)$ is $d(\alpha)$.

\begin{lemma}\label{eve}
$$\ell(v)-\ell(w)+1=\codim X(w) - \codim X(v)  +1 \leq \int_{d(\alpha)}c_1(TX)$$
\end{lemma}
This inequality follows from \cite[Corollary 5.4.2]{FFKM} in the case $P=B$. Lemma \ref{eve}  is proved
in Section \ref{newy}.

Following \cite{FW}, let $E(v,w)$ be the intersection 
$e_1^{-1}(Y(v))\cap e_2^{-1}(X(w))\subseteq \overline{M}_{0,3}(X,d)$
where $d=d(\alpha)$ and $Y(v)$ is the $B^-$-orbit closure of $vP$, and $e_1$, $e_2$ the evaluations at the first two of the three
 marked points. Note that $\codim Y(v) +\codim X(v)=\dim X$. The inequality in Lemma \ref{eve} is equivalent to
 $\dim E(v,w)\geq 1$, and equality would imply that $e_3(E(v,w))\subseteq X$ is the unique  $T$-fixed curve of degree $d(\alpha)$ joining $\dot{v}$ and $\dot{w}$.

\begin{lemma}\label{stringy}
Suppose equality holds in the inequality in Lemma \ref{eve}. Thus $w< v$.  Let $x\in X^P_v$ be in the Schubert cell parameterized by $v$. Then there is a unique rational curve of degree $d(\alpha)$ joining $x$ to some point (which is not a priori chosen) of $X^P_w$.
\end{lemma}
\begin{proof}
We can assume that $x=v P$ by translations. Then a point on such a curve lies in $E(v,w)$ which is $T$-fixed by the above discussion.

\end{proof}

\section{The boundary  } We want to understand the complement $\Omega(\vec{w},d)\setminus\Omega^0(\vec{w},d)$. A part of the boundary is the ``classical boundary'' formed by the $i$th marked point degenerating to the boundary of the Schubert cell parameterised by $w_i$. The other part is the ``quantum type'' boundary where the map $C\to \Bbb{P}^1$ in the definition of $\Omega(\vec{w},d)$ is not an isomorphism.

\begin{defi}\label{Hermann}
Fix $1\le j\le s$ and $v\in W^P$ such that $v\xrightarrow{\beta} w_j$. Let $v_i=w_i$ for $i\neq j$ and $v_j=v$.
Then  $\Omega(\vec{v},d)\subseteq \Omega(\vec{w},d)$ is codimension one, and is called a classical type boundary stratum
of $\Omega(\vec{w},d)$ corresponding to $(\beta,j)$.
\end{defi}

\begin{defi}\label{Maxim}
Fix $1\le j\le s$ and $\alpha\in R^+\setminus R_{\frl}^+$. Set $v=\overline{w_js_{\alpha}}\in W^P$. Assume (see \eqref{degdalpha})
$$\ell(v)-\ell(w_j)+1 =\int_{d(\alpha)}c_1(TX).$$
To this data we will associate a codimension-$1$ quantum-type stratum $\Omega_K(\vec{w},\alpha,j,d)\subset \Omega(\vec{w},d)$. The generic point of this stratum corresponds to a reducible curve $C$ with two components. The main component $C_0$ maps isomorphically to $\Bbb{P}^1$ under $p$. The other component $C_1$  attaches to $C_0$ at a node $R$, and  $C_1$ lies over $p_j$. $C_1$ carries the marked point $q_j$. In addition to the conditions in the definition of $\Omega(\vec{w},d)$, we will assume that
\begin{enumerate}
\item $F(R)$ and $\overline{g_j}$ are in relative position $v$.
\item The degree of $F:C_1\to \me_{p_j}/P$ is $d(\alpha)$.
\end{enumerate}
\end{defi}
\begin{remark}

 By Lemma \ref{stringy}, the quantum stratum  $\Omega_K(\vec{w},\alpha,j,d)\subset \Omega(\vec{w},d)$ is isomorphic to $\Omega(\vec{v},d-d(\alpha))$ over $\Parbun_G$, by the map that forgets the component $C_1$. This will be used in the computation of the push forward of the cycle class of  $\Omega_K(\vec{w},\alpha,j,d)$ to $\Parbun_G$: This class is the same as the push forward of
$\Omega(\vec{v},d-d(\alpha))$.
\end{remark}

\begin{theorem}\label{codim1boundary}  Assume that $\Omega(\vec{w},d)$ is irreducible of relative dimension $0$ over $\Parbun_G$. Let $Y$ be the union of all classical and quantum-type strata of $\Omega(\vec{w},d)$. Then  $(\Omega^0(\vec{w},d)\cup Y)\to \Parbun_G$ is surjective in codimension one (i.e., 
if $Z$ is the closure of the image of $(\Omega^0(\vec{w},d)\cup Y)\to \Parbun_G$, then each irreducible component of  $\Parbun_G-Z$ has codimension $\geq 2$ in $\Parbun_G$)
\end{theorem}

\begin{remark}
The Chevalley formulas describe the product of a cycle class $X^P_w$ with a positive divisor class $X^P_{s_i}$ where $\alpha_i\in S_P$. The terms that appear in the classical Chevalley formula in cohomology correspond to the codimension-one boundary of the Schubert intersection stack in \cite{BKR,BKiers} (here we run over choices of $w$ given by the Schubert conditions $w_1,\dots,w_s$). These correspond to the classical boundary \ref{Hermann} above and  the ordinary Bruhat order. Under the condition of Levi-movability (with corresponding intersection number $1$) only the boundary strata corresponding to simple roots map to codimension one in the product of flag varieties for $G$ (the analogue of  $\Parbun_G$) by results of \cite{BKR}.

In a similar manner, the terms that appear in the quantum Chevalley formula \cite{FW,LamS} correspond to the classical and quantum boundary above. The Bruhat order is the semi-infinite Bruhat order. In the presence of the condition  $\langle \sigma_{w_1}^P, \hdots, \sigma_{w_s}^P \rangle_d^{\circledast_0} =1$ (which implies Levi-movability), only the boundary strata corresponding to simple roots (in the semi-infinite order) map to codimension one in $\Parbun_G$ by Theorems \ref{T1} and \ref{T2} below (which generalize results of \cite{BKR}).
\end{remark}

\subsection{Proof of Theorem \ref{codim1boundary}}\label{proofcodim1boundary}
\begin{proof}This proof is modelled on the proof of \cite[Theorem 5.2]{FFKM}. Let $D$ be an irreducible component of $\Omega(\vec{w},d)\setminus\Omega^0(\vec{w},d)$ which is not of classical-type, and $(C,p,F,\me,\bar{g}_1,\dots,\bar{g}_s)$ a general point of $D$.{ If $C$ is irreducible then
it is easy to see (by dimension counting) that the point lies on one of the classical boundary strata as in Def. \ref{Hermann}. Therefore we will assume that $C$ is reducible.}

The curve $C$ has a main component $C_0$ mapping isomorphically to $\Bbb{P}^1$, and has a comb-like structure. The teeth of the comb are chains of $\Bbb{P}^1$'s with at most one marked point on each tooth.

Let $\mathcal{N}$ be the nodes on $C_0$. By our condition each $R\in\mathcal{N}$ has to lie over some marked point $p_i\in \Bbb{P}^1$. This is because otherwise the  attaching point of the tooth of the comb over $R$ can slide over on the main component, with the same image in $\Parbun_G$ (note that in this case the tooth over $R$ will not have a marked point).

For simplicity assume $\mathcal{N}=\{R_1,\dots,R_m\}$  where $R_1,\dots,R_m$ lie over $p_1,\dots,p_m\in\Bbb{P}^1$ (in that order). See Figure \ref{diag-2}.

\begin{figure}
\begin{tikzpicture}[scale=2]

\draw[-,thick] (0,0) -- (3,0);
\draw[-,thick] (0,-1.5) -- (3,-1.5);

\draw[->] (2,-0.4) -- (2,-1.1);
\node at (2.13,-0.75) {\small $p$};

\draw[->] (3.5,0.4) -- (4.2, 0.4);
\node at (3.85,0.53) {\small $F$};
\node at (4.8,0.4) {$\me/P$};

\node at (-0.2,-1.5) {$\P^1$};
\node at (-0.2,0) {$C_0$};

\draw[-,thick] (0.2,-0.4) -- (0.2,1);
\draw[-,thick] (-0.5,1.3) -- (0.3, 0.7);
\draw[-,thick] (-0.4,1.1) -- (-0.4,1.7);
\node at (-0.4, 2) {$C_1$};
\node at (0.3,0.07) {\tiny $R_1$};

\node at (-0.4,1.5) [circle,fill,inner sep=1.5pt]{};
\node at (-0.53,1.5) {\small $q_1$};
\node at (0.2,-1.5) [circle,fill,inner sep=1.5pt]{};
\node at (0.2,-1.63) {\small $p_1$};

\draw[-,thick] (0.8,-0.4) -- (0.8,1);
\draw[-,thick] (1.5,1.3) -- (0.7, 0.7);
\draw[-,thick] (1.4,1.1) -- (1.4,1.7);
\draw[-,thick] (1.5, 1.5) -- (0.7,2.1);
\node at (0.5, 2.15) {$C_2$};
\node at (0.9,0.07) {\tiny $R_2$};

\node at (0.8,2.025) [circle,fill,inner sep=1.5pt]{};
\node at (0.8,2.025-0.13) {\small $q_2$};
\node at (0.8,-1.5) [circle,fill,inner sep=1.5pt]{};
\node at (0.8,-1.63) {\small $p_2$};

\node at (1.4,0) [circle,fill,inner sep=1.5pt]{};
\node at (1.4,-0.13) {\small $q_3$};
\node at (1.4,-1.5) [circle,fill,inner sep=1.5pt]{};
\node at (1.4,-1.63) {\small $p_3$};

\node at (2.2, -0.13) {\small $\cdots$};
\node at (2.2, -1.63) {\small $\cdots$};

\node at (2.8,0) [circle,fill,inner sep=1.5pt]{};
\node at (2.8,-0.13) {\small $q_s$};
\node at (2.8,-1.5) [circle,fill,inner sep=1.5pt]{};
\node at (2.8,-1.63) {\small $p_s$};

\draw[->] (4.65,0) -- (3.2,-1.2);

\end{tikzpicture}
\caption{}\label{diag-2}
\end{figure}

Define $v_1,\dots,v_s\in W^P$ as follows:
\begin{enumerate}
\item Let $v_i=w_i$ if $i>m$.
\item  If $1\leq i\leq m$, $R_i$ lies on a tooth $C_i$ of $C$ which maps to a node $R_i$ on $C_0$. Let $F(R_i)$ be in the Schubert variety parametrized by $v_i$ with respect to the flag $\bar{g}_i.$
\end{enumerate}
Let $d_0$ be the  degree of the $P$-reduction $C_0\to \me/P$.

It is easy to see that $D\to \Parbun_G$ factors (generically) through a subset of $\Omega(\vec{v},d_0)$ and the fibres of $D\to \Omega(\vec{v},d_0)$ are products of varieties $Q_1\times \dots\times Q_m$, with one factor $Q_j$ for  every node $R_j$ on $C_0$.  Here $Q_j\subseteq \overline{M}_{0,2}(X,d_j)$, $d_j=\int_{C_j} c_1(TX)$ is the space of stable maps so that the first marked point maps to $\dot{v_j}\in X({v_j})$ and the second one goes to an arbitrary point in $X({w_j})$.

 Hence the relative dimensions over $\Parbun_G$ satisfy
$$\dim \Omega(\vec{v},d_0)-\dim \Omega(\vec{w},d) \geq -1.$$
(we put $\geq$ because $D\to \Parbun_G$ only factors through  $\Omega(\vec{v},d_0)$)
 The above equality gives
\begin{equation}\label{IQ1}
\sum_{j=1}^m (\codim X(w_j) -\codim X(v_j) -\int_{C_j} c_1(TX))\geq -1
\end{equation}

The factors $Q_j$ have therefore got to be zero-dimensional since they lie over the same point of $\Parbun_G$.
They also have to be $T$-fixed. The following picture of  $Q_j$'s emerges. For simplicity in notation we describe $Q_1$:
$C_1$ is a chain of rational curves $C(1),\dots,C(r)$.

Suppose that the image of  $C(k)$ in $G/P$ is  the (reduced) curve $E(k)$, $k=1,\dots,r$.
Let $b(k)\geq 1$ be the degree of $C(k)\to E(k)$. The curves $E(k)$ are $T$-fixed ({after translation by $\overline{g}_k$}) with $E(1)$ joining $\dot{v}_1$ to
$\dot{u}_1$, $E(2)$ joining $\dot{u}_1$ to $\dot{u}_2$, and $E(r)$ joining $\dot{u}_{r-1}$ to $\dot{u}_r)$ with $u_r=w_1$. Let $u_0=v_1$.
We have $u_{k-1}=\overline{u_k s_{\beta(k)}}\in W_P$ for $k=1,\dots,r$ where $\beta(k)$ are positive roots (use the form of $T$ fixed curves given in Section \ref{TFIX}).
\begin{remark}
If $q_1$ is not on $C(r)$, then we can deform the last component and the map to $G/P$ by varying the point where
$C(r)$ attached to $C(r-1)$. Therefore in Figure \ref{diag-2}, each point $q_i$ is on the last component.
\end{remark}

Repeated application of Lemma \ref{eve} gives
$$\codim X(w_1)-\codim X(v_1) - \int_{\sum_{k=1}^r [E(k)]}c_1(TX) \leq -r.$$
Note that $\sum [E(k)]\leq\sum b(k)[E(k)]=  d_j$, and hence the above gives for $j=1,\dots,m$,
\begin{equation}\label{IQ2}
\codim X(w_j) -\codim  X(v_j) - \int_{d_j}c_1(TX) \leq -r_j
\end{equation}
where $r_j$ is the length of the chain $C_j$.

Comparing inequalities \eqref{IQ1} and \eqref{IQ2}, we get $m=1$ and $r_1=1$, and $C_1\cong\Bbb{P}^1$ mapping isomorphically to its image. The theorem follows.
\end{proof}

\section{Proof of Lemma \ref{eve}, and consequences}\label{newy}
 Let $v,w\in W^P$ and $v = \overline{ws_\alpha}$, where $\alpha \in R^+\setminus R_{\mathfrak{l}}^+$. We need to prove
\begin{equation}\label{newEE}
\ell(v)-\ell(w)+1 \leq \bigl(\sum_{\gamma\in R^+\setminus R_{\frl}^+}\gamma\bigr)(\alpha^{\vee})
\end{equation}
The RHS is a sum of integers that can be positive, negative or zero. Note that $\ell(w)$ is equal to the number of roots $\gamma\in R^+\setminus R_{\frl}^+$ such that $w \gamma\prec 0$; likewise, $\ell(v)$ is equal to the number of such $\gamma$ satisfying $v\gamma \prec 0$. In fact $\ell(v)$ is  the number of $\gamma\in R^+\setminus R_{\frl}^+$ satisfying $ws_{\alpha}\gamma \prec 0$.

We first run through $\gamma$ so that $\gamma$ and $s_{\alpha}\gamma$ are (possibly equal and) both in $R^{+}\setminus R_{\frl}^+$.
The contribution on the RHS is $(\gamma+s_{\alpha}(\gamma))(\alpha^{\vee})=0$ (or one half of this). If $\gamma$ contributes for $w$ then $\gamma'=s_{\alpha}\gamma$ contributes for $ws_{\alpha}$. The net count is zero. That is, the contribution to $\ell(v)-\ell(w)$ from these $\gamma$ as well as their contributions to the RHS are both zero.

 Now we run through $\gamma$ so that  $s_{\alpha}\gamma$ is not in $R^{+}\setminus R_{\frl}^+$ ($\gamma=\alpha$ is one of these). It is easy to see that
$\gamma(\alpha^{\vee})>0$. These $\gamma$ can contribute at most one on the LHS (if they contribute one to $ws_{\alpha}$ and none to $w$), and on the RHS they give a positive integer $\gamma(\alpha^{\vee})$. So the RHS weighs at least as  much for such $\gamma$. In the case $\gamma=\alpha$, we get $2$ on the RHS.  Hence inequality \eqref{newEE} holds, and the proof of Lemma \ref{eve} is complete.

Now assume that equality holds in Lemma  \ref{eve}. The above computations yield the following lemmas. In the case $\gamma=\alpha$, we want the contribution of $\gamma$ in $\ell(\overline{ws_{\alpha}})-\ell(w)$ to be one (note that $\alpha(\alpha^{\vee})=2$), i.e.,
\begin{lemma}\label{above0}
$w\alpha\succ 0$.
\end{lemma}
For the remaining cases we have the following.
\begin{lemma}\label{above}
Suppose $\gamma\neq \alpha$, $\gamma\in R^+ \setminus R_{\frl}^+$, and $s_{\alpha}\gamma\not\in R^+\setminus R_{\frl}^+$. Then $\gamma(\alpha^{\vee})=1$, $w\gamma\succ 0$ and $ws_{\alpha}\gamma\prec 0$.
\end{lemma}
There are several other properties in the setting of equality in Lemma \ref{eve}.  The first one is a formula for the length of $\ell(ws_{\alpha})$ itself, not the length of the minimal length representative \cite[Lemma 10.18]{LamS} (who employ an opposite notation):
\begin{lemma}\label{equivalent}
$\ell(ws_{\alpha})=\ell(w)-1 +\bigl(\sum_{\gamma\in R^+}\gamma\bigr)(\alpha^{\vee})$.
\end{lemma}

This is implied by the following:

\begin{lemma}\label{anotherlemma?}
Suppose $\gamma\in R_{\frl}^+$, then either $ws_{\alpha}\gamma\succ 0$ and $\gamma({\alpha}^{\vee})=0$; or $ws_{\alpha}\gamma\prec  0$ and $\gamma({\alpha}^{\vee})=1$.
\end{lemma}

\begin{proof}
It is easy to see that if $\gamma\in R_{\frl}^+$, then $\gamma(\alpha^{\vee})\geq 0$. Otherwise, $\gamma'=s_{\alpha}\gamma\in R^+\setminus R_{\frl}^+$ and $s_{\alpha}\gamma'\not\in R^+\setminus R_{\frl}^+$. By Lemma \ref{above}, we get $w\gamma\prec 0$, which is a contradiction.

If $\gamma(\alpha^{\vee})=0$, then obviously $ws_{\alpha}\gamma=w\gamma\succ 0$.

So suppose $\gamma(\alpha^\vee)>0$. The root $s_\alpha\gamma = \gamma - \gamma(\alpha^\vee)\alpha$ clearly has a negative coefficient on some simple root $\alpha_i$ not in $R_{\frl}^+$, so belongs to $R^-\setminus R_{\frl}^-$. So applying Lemma \ref{above} to $-s_\alpha \gamma$ (note that $-s_\alpha \gamma \ne \alpha$), we must have $(-s_\alpha\gamma)(\alpha^\vee) = 1$, i.e., $\gamma(\alpha^\vee) = 1$, and $w (-s_\alpha \gamma) \succ 0$, i.e., $ws_\alpha \gamma \prec 0$.
\end{proof}

For us, the most prominent instance of (\ref{newEE}) will be when $w\alpha = \theta$, the highest root. In this specific instance, (\ref{newEE}) holds with equality. 

\begin{lemma}\label{old8-7}
Let $w\in W^P$ and suppose $w\alpha = \theta$, where $\alpha \in R^+\setminus R_{\frl}^+$. Set $v = \overline{ws_\alpha} = \overline{s_\theta w} \in W^P$. Then 
$$
\ell(v) - \ell(w)+1 = \left(\sum_{\gamma \in R^+ \setminus R_{\frl}^+} \gamma \right) (\alpha^\vee). 
$$
\end{lemma}

\begin{proof}
Running through $\gamma$ such that $s_\alpha \gamma\not\in R^+\setminus R_{\frl}^+$, we first have the special case $\gamma = \alpha$ and thus $\gamma(\alpha^\vee) = 2$. Note that $ws_\alpha(\gamma) = -w\alpha = -\theta$, so the contribution on the LHS of (\ref{newEE}) is $1+1=2$. 

Otherwise, $\gamma \ne \alpha$ but $\gamma(\alpha^\vee)>0$. Now, 
$$
\gamma(\alpha^\vee) = w^{-1}\gamma(\theta^\vee)
$$
which takes on only the values $-1,0,1$ except for the special cases $w^{-1}\gamma = \pm \theta$. Hence $\gamma(\alpha^\vee) = 1$, and the contribution of $\gamma$ on the RHS of (\ref{newEE}) is $1$. Due to being the highest root, the only way for $\theta$ to be expressed as the sum of two roots is for both summands to be positive. Noting that 
$$
\theta = w(-s_\alpha\gamma + \gamma) = -ws_\alpha \gamma + w \gamma,
$$
we see that $w\gamma\succ 0$ and $ws_\alpha\gamma\prec 0$, implying that the contribution from $\gamma$ to the LHS is also $1$.
\end{proof}

\section{Levi movability statements}
Assume that $\langle \sigma_{w_1}^P,\hdots,\sigma_{w_s}^P\rangle_d^{\circledast_0}=1$ (which implies Levi-movability) in this section. The first result is about the classical-type boundary and the ramification divisor of the map $\Omega(\vec{w},d)\to \Parbun_G$.
\begin{theorem}\label{T1}
Fix $1\le j\le s$ and $v\in W^P$ such that $v\xrightarrow{\beta} w_j$. Let $v_i=w_i$ for $i\neq j$ and $v_j=v$. Note that $\Omega(\vec{w},d)$ is smooth at the generic point of $\Omega(\vec{v},d)\subseteq \Omega(\vec{w},d)$.
\begin{enumerate}
\item[(a)] If $\beta$ is not a simple root, then the classical type boundary stratum (see Def. \ref{Hermann}) $\Omega(\vec{v},d)\subseteq \Omega(\vec{w},d)$ lies in the ramification divisor of the map $\Omega(\vec{w},d)\to \Parbun_G$ (See Section \ref{ramdef}).
\item[(b)] If $\beta$ is  a simple root, then the classical type boundary stratum (see Def. \ref{Hermann}) $\Omega(\vec{v},d)\subseteq \Omega(\vec{w},d)$ does not lie in the ramification divisor of the map $\Omega(\vec{w},d)\to \Parbun_G$.
\end{enumerate}
\end{theorem}

The second result is about the quantum-type boundary and the ramification divisor of the map $\Omega(\vec{w},d)\to \Parbun_G$.
\begin{theorem}\label{T2}
Fix $1\le j\le s$ and $\alpha\in R^+-R_{\frl}^+$. Set $v=\overline{w_js_{\alpha}}\in W^P$. Assume (see \eqref{degdalpha})
$$\ell(v)-\ell(w_j)+1=\int_{d(\alpha)}c_1(T_X).$$
To this data recall that we have  defined a  codimension-one quantum-type stratum $\Omega_K(\vec{w},\alpha,j,d)\subset \Omega(\vec{w},d)$ in Def.\ref{Maxim}. Let $\theta$ be the highest root.
\begin{enumerate}
\item[(a)] If $w_j \alpha\neq \theta$, then the quantum-type boundary stratum  $\Omega_K(\vec{w},\alpha,j,d)\subset \Omega(\vec{w},d)$ lies in the ramification divisor of the map $\Omega(\vec{w},d)\to \Parbun_G$ (see Section \ref{ramdef}).
\item[(b)]  If $w_j \alpha=\theta$, then the quantum-type boundary stratum  $\Omega_K(\vec{w},\alpha,j,d)\subset \Omega(\vec{w},d)$ does not lie in the ramification divisor of the map $\Omega(\vec{w},d)\to \Parbun_G$.
\end{enumerate}
\end{theorem}
\begin{remark}\label{contracted}
In Theorem \ref{T1}(a),  $\Omega(\vec{w},d)\to \Parbun_G$ is birational and hence $\Omega(\vec{v},d)\subseteq \Omega(\vec{w},d)$ is contracted under $\Omega(\vec{w},d)\to \Parbun_G$.

Similarly in Theorem \ref{T2}(a),   $\Omega(\vec{w},d)\to \Parbun_G$ is birational and hence $\Omega_K(\vec{w},\alpha,j,d)\subset \Omega(\vec{w},d)$ is contracted under $\Omega(\vec{w},d)\to \Parbun_G$.

\end{remark}

Parts (b) parts of the theorems follow from Lemma \ref{detailed} which shows that the cycle classes of the push-forwards to  $\Parbun_G$ are not zero. Parts (a) of the theorems are proved in the next two sections.

\section{Levi-movability analysis}\label{LeviM1}
\begin{defi}
Define an element $x_P\in \mathfrak{h}$ by the requirements
$$
\alpha_k(x_P) =
\left\{
\begin{array}{ll}
1, & \alpha_k\in S_P\\
0, & \text{else}.
\end{array}
\right.
$$
\end{defi}
The tangent space of $G/P$ at $\dot{e}=eP\in G/P$  admits a direct sum decomposition
$$T_{\dot{e}}(G/P)=\bigoplus_{\beta\in R^+\setminus R_{\frl}^+} T_{\dot{e}}(G/P)_{-\beta}$$

Let $$V_j=\bigoplus_{\beta\in R^+\setminus R_{\frl}^+,\beta(x_P)=j} T_{\dot{e}}(G/P)_{-\beta}$$
so that for a suitable $r$,
\begin{equation}\label{tspacedeco}
T_{\dot{e}}(G/P)=\bigoplus_{j=1}^{r} V_j
\end{equation}
The subspaces $\oplus_{j=1}^d V_d$ are $P$-stable, for $d=1,\dots,r$ and hence define a filtration of the tangent bundle of $G/P$ over $G/P$ by subbundles which we denote by
$$\mathcal{F}_1\subset\dots\subset\mathcal{F}_{r}.$$

The rank of $\mathcal{F}_j$ is the cardinality of
$$S_j=\{\gamma\in R^+\setminus R_{\frl}^+: \gamma(x_P)\leq j\}.$$

\subsection{Graded dimensions}\label{gradedDim}
We will consider on $\Bbb{P}^1$, bundles $\mathcal{G}$ with a filtration by subbundles
 $$\mathcal{G}_1\subset\dots\subset\mathcal{G}_{r}$$
  with subspaces $B_i\subseteq \mg_{p_i}$ for $i=1,\dots, s$, and corresponding maps
\begin{equation}\label{mappo}
H^0(\Bbb{P}^1,\mg)\to \oplus_{i=1}^s \mg_{p_i}/B_i.
\end{equation}
To such a filtration, we will assign the number
$$
\deg= \sum_{j}\left(\chi(\Bbb{P}^1,\mg_j)-\sum_{i=1}^s \bigl(\dim (\mathcal{G}_j)_{p_i}-\dim \bigl((\mathcal{G}_j)_{p_i}\cap B_i)\bigr)\right).
$$
Note that if this assigned number is positive, then \eqref{mappo} has a non-trivial kernel, because one of the $j$ summands in the above equation will then have to be positive, and we use $h^0\geq h^0-h^1$ and the maps
\begin{equation}\label{mappo2}
H^0(\Bbb{P}^1,\mg_j)\to \oplus_{i=1}^s (\mg_j)_{p_i}/((\mathcal{G}_j)_{p_i}\cap B_i).
\end{equation}

\subsection{Levi-Movability}\label{LeviM2}
Once again consider a datum
\begin{equation}\label{dagger2}
\langle \sigma^P_{w_1},\dots,\sigma^P_{w_{s-1}},\sigma^P_{w_s}\rangle_d^{\circledast}=1.
\end{equation}

Let $x=(\mathcal{E},\overline{g_1},\hdots,\overline{g_s},f)\in \Omega^0(\vec{w},d)$ with $\mathcal{E}$ trivial be a general point of $\Omega^0(\vec{w},d)$, and $f:
\Bbb{P}^1\to G/P$ be the corresponding map (using $\me/P=\Bbb{P}^1\times G/P$ since $\me$ is trivial).

Let $A_x$ be the fiber of $\Omega(\vec{w},d)\to \Parbun_G$ over $(\mathcal{E},\overline{g_1},\hdots,\overline{g_s})$.
The tangent space to $A_x$ is the kernel of the map to a direct sum of sky-scraper sheaves
\begin{equation}\label{tsp0}
H^0(\Bbb{P}^1,f^*T(G/P))\to \bigoplus \frac{T(G/P)_{f(p_i)}}{T\overline{g}_i X(w_i)}.
\end{equation}

The tangent space $T\overline{g}_i X(w_i)$ is taken at $f(p_i)\in \overline{g}_i X(w_i)$ which is in the open cell. The filtration on $G/P$ pulls back to a filtration on $f^*T(G/P)$ and also a corresponding filtration on  the target of \eqref{tsp0}.

The assumption of Levi-movability implies that the number $\deg$ assigned to \eqref{tsp0} as in Section \ref{gradedDim} is zero. In fact, each of the $j$-terms is separately zero (Section 3.8.1 of \cite{BKq}).
\begin{remark}
Let $\mathcal{F}_j$ denote the pull-back of $\mathcal{F}_j$ via $f:\Bbb{P}^1\to G/P$.
The rank of $\mathcal{F}_j$ is the cardinality of
$$S_j=\{\gamma\in R^+\setminus R_{\frl}^+: \gamma(x_P)\leq j\}$$
The degree is $\deg\mathcal{F}_j=\sum a_i \sum_{\gamma\in S_j}\gamma(\alpha_i^{\vee}).$

The filtered piece corresponding to  $\frac{T(G/P)}{T\overline{g}_i X(w_i)}$ has rank
$$|S_j|-|\{\gamma\in S_j:w_i\gamma\prec 0\}|.$$

The Levi-movability condition is equivalent to $\forall j, \chi(\mathcal{F}_j)=\sum_{i=1}^s (|S_j|-|\{\gamma\in S_j:w_i\gamma\prec 0\}|)$.
\end{remark}

\subsection{Proof of Theorem \ref{T1}(a)}
Continue with the assumption 
\begin{equation}\label{dagger3}
\langle \sigma^P_{w_1},\dots,\sigma^P_{w_{s-1}},\sigma^P_{w_s}\rangle_d^{\circledast_0}=1.
\end{equation}
Let $v\in W^P$ be  such that $v\xrightarrow{\beta} w_1$. Let $v_i=w_i$ for $i\neq 1$ and $v_1=v$.

Let
$y=(\mathcal{E},\overline{g_1},\hdots,\overline{g_s},f)\in \Omega^0(\vec{v},d)$ with $\mathcal{E}$ trivial be a general point of $\Omega^0(\vec{v},d)\subseteq \Omega(\vec{w},d)$, and $f:
\Bbb{P}^1\to G/P$ be the corresponding map (using $\me/P=\Bbb{P}^1\times G/P$ since $\me$ is trivial). Let $x$ be the corresponding point of $\Parbun_G$.

Let $A_x$ be the fiber of $\Omega(\vec{w},d)\to\Parbun_G$ over $x$. Note that the map $\Omega(\vec{w},d)\to \Parbun_G$ after a base change is a map of  schemes of the same dimensions, and $y$ is a smooth point of $\Omega(\vec{w},d)$ (since Schubert varieties are normal).  The tangent space to $A_x$ at $y$ is the kernel of the map to a direct sum of sky-scraper sheaves
\begin{equation}\label{tsp}
H^0(\Bbb{P}^1,f^*T(G/P))\to \bigoplus \frac{T(G/P)_{f(p_i)}}{T\overline{g}_i X(w_i)}.
\end{equation}

The tangent space $T\overline{g}_i X(w_i)$ is taken at $f(p_i)\in \overline{g}_i X(w_i)$ which is in the boundary $\overline{g}_i X^0(v)\subseteq \overline{g}_i X(w_1)$ when $i=1$, and in the open cell otherwise.

The filtration on $G/P$ pulls back to a filtration on $f^*T(G/P)$ and also a corresponding filtration on  the target of \eqref{tsp}. Assume that $\beta$ is not a simple root.

This situation gives another instance of the setting of Section \ref{gradedDim}. We note that  there exists a $j$ such that the rank of the $j$th filtered piece corresponding to  $\frac{T(G/P)}{T\overline{g}_1 X(w_1)}$ has rank strictly less than the rank  in \eqref{tsp0} (which is taken at a generic point of $\overline{g}_1 X(w_1)$). This is proved in \cite[Theorem 7.4]{BKR} (also see the semi-continuity inequality \cite[Equation (33)]{BKR}).
This means that the the $j$th-term in $\deg$ is greater than the corresponding term for \eqref{tsp0} (which is zero). It follows  that \eqref{tsp} is not an isomorphism which proves Theorem \ref{T1}(a).

\section{Proof of Theorem \ref{T2}(a)}
Consider a  codimension $1$ quantum-type stratum $\Omega_K(\vec{w},\alpha,j,d)\subset \Omega(\vec{w},d)$ as in Def. \ref{Maxim}. For simplicity assume $j=1$. Choose a general point  $y\in \Omega_K(\vec{w},\alpha,j,d)$. Assume that $y$ lies over $x=(\mathcal{E},\overline{g_1},\hdots,\overline{g_s})\in \Parbun_G$.

 The data of $y$ corresponds to a reducible curve $C$ with two components. The main component $C_0$ maps isomorphically to $\Bbb{P}^1$ under $p$. The other component $C_1$  attaches to $C_0$ at a node $R$, and  $C_1$ lies over $p_1$. $C_1$ carries the marked point $q_1$. See Figure \ref{diag-3}.

\begin{figure}
\begin{tikzpicture}[scale=2]

\draw[-,thick] (0,0) -- (3,0);
\draw[-,thick] (0,-1.5) -- (3,-1.5);

\draw[->] (2,-0.4) -- (2,-1.1);
\node at (2.13,-0.75) {\small $p$};

\draw[->] (3.5,0.4) -- (4.2, 0.4);
\node at (3.85,0.53) {\small $F$};
\node at (4.8,0.4) {$\me/P$};

\node at (-0.2,-1.5) {$\P^1$};
\node at (-0.2,0) {$C_0$};

\draw[-,thick] (0.2,-0.4) -- (0.2,1);
\node at (0.2, 1.3) {$C_1$};
\node at (0.3,0.07) {\tiny $R$};

\node at (0.2,0.9) [circle,fill,inner sep=1.5pt]{};
\node at (0.07,0.9) {\small $q_1$};
\node at (0.2,-1.5) [circle,fill,inner sep=1.5pt]{};
\node at (0.2,-1.63) {\small $p_1$};

\node at (0.8,0) [circle,fill,inner sep=1.5pt]{};
\node at (0.8,-0.13) {\small $q_2$};
\node at (0.8,-1.5) [circle,fill,inner sep=1.5pt]{};
\node at (0.8,-1.63) {\small $p_2$};

\node at (1.4,0) [circle,fill,inner sep=1.5pt]{};
\node at (1.4,-0.13) {\small $q_3$};
\node at (1.4,-1.5) [circle,fill,inner sep=1.5pt]{};
\node at (1.4,-1.63) {\small $p_3$};

\node at (2.2, -0.13) {\small $\cdots$};
\node at (2.2, -1.63) {\small $\cdots$};

\node at (2.8,0) [circle,fill,inner sep=1.5pt]{};
\node at (2.8,-0.13) {\small $q_s$};
\node at (2.8,-1.5) [circle,fill,inner sep=1.5pt]{};
\node at (2.8,-1.63) {\small $p_s$};

\draw[->] (4.65,0) -- (3.2,-1.2);

\end{tikzpicture}
\caption{}\label{diag-3}
\end{figure}

\begin{enumerate}
\item $F(R)$ and $\overline{g_1}$ are in relative position $v$.
\item The degree of $F:C_1\to \me_{p_1}/P$ is $d(\alpha)>0$. Note that the degree of $P$-reduction on the main component can be negative.
\end{enumerate}
Choose a trivialization of $\me$ at $p_1$, and assume that $F(R)$ is the point $\dot{v}$ and assume that $\bar{g}_1=\bar{e}$ by translation. Note that by Lemma \ref{stringy}, $C_1$ is the $T$-fixed curve joining $\dot{v}$ and $\dot{w}_1$ and $f(q_1)=\dot{w}_1$. Also,
\begin{lemma}\label{curve}
The Schubert variety $X(v)$ contains the curve $E=C(v,w_1)\subseteq G/P$.
Further $TE_{w_1}\cap (TX({w_1}))_{\dot{w}_1}=0$.
\end{lemma}
\begin{proof}
$C(v,w_1)$ is  the $v$-translate of the curve joining $\dot{e}$ and $\dot{s}_{\alpha}$, i.e.,
$t\mapsto vU_{-\alpha}(t)\dot{e}$ with $t=0$ corresponding to $f(R)$ and $t=\infty$ to $\dot{w}_1$. Note that the mapping can be rewritten as
$t\mapsto U_{-v\alpha}(t)\dot{e}$, but $v\alpha=-w_1\alpha\in R^-$ by Lemma \ref{above0}. Therefore the curve $f|D$ is contained in $X(v)$. 

Note further that by a calculation for $\op{SL}(2)$,
\begin{equation}
vU_{-\alpha}(t)\dot{e} =w U_{-\alpha}(-T)\dot{e}, \text{  where } T=1/t
\end{equation}
This gives the second statement, since $w\alpha\in R^{+}-R_{\frl}$.
\end{proof}

\subsubsection{Some deformation theory}\label{ramdef}
We recall some deformation theory from \cite{Behrend,BFant,Campbell,GHS}.
Under our assumptions
$F:C\to \me/P$ is an embedded local complete intersection curve. Let $\Omega^E(\vec{w},d)\subset\Omega(\vec{w},d)$ be the open substack parameterizing $F$ which are embedding and the relative position in item (5) of the definition of $\Omega(\vec{w},d)$ an equality. Clearly the point $y\in \Omega_K(\vec{w},\alpha,j,d)\subset \Omega(\vec{w},d)$
satisfies this condition.

Since $F$ is an embedding, the map $\Omega^F(\vec{w},d)\to \Parbun_G$ is a representable map of smooth Artin stacks by \cite{Campbell}: The representability is because $F$ has no automorphisms. We can consider the ramification locus of such a map via base change: If $S\to  \Parbun_G$ is an atlas, then $\Omega^F(\vec{w},d)\times_S \Parbun_G$, and the the maps $\Omega^F(\vec{w},d)\times_S \Parbun_G\to \Omega^F(\vec{w},d)$ and  $\Omega^F(\vec{w},d)\to \op{Spec}(\Bbb{C})$ are smooth. Therefore $\Omega^F(\vec{w},d)\times_S \Parbun_G\to S$ is a map between smooth schemes of the same dimension. We can consider the tangent space to the fibers to look for ramification phenomena.

Ignoring the Schubert conditions, the first-order deformations (i.e., deformations of the map $F$ keeping $\me$ fixed) are given by the hypercohomology groups $\Bbb{H}^1(C,\Bbb{R}Hom(\Omega_F^{\bullet},\mathcal{O}_C))$ (and the obstructions are given by $\Bbb{H}^2(C,\Bbb{R}Hom(\Omega_F^{\bullet},\mathcal{O}_C)$) where $\Omega_F^{\bullet}$ is the complex $F^*\Omega_{\me/P}\to \Omega_{{C}}$ placed in degrees $-1$ and $0$.

\begin{defi}
The  map $F^*\Omega_{\me/P}\to \Omega_{{C}}$ is a surjection. Since $F$ is a local complete intersection embedding, the kernel of this surjection  is locally free (see \cite[Appendix B.7]{Fulton_Int}, \cite[Section 2.1]{GHS}). Its dual is denoted  by  $N$, the ``normal bundle" (the description as normal bundle is valid at all smooth points).

Note that $\Bbb{H}^1(C,\Bbb{R}Hom(\Omega_F^{\bullet},\mathcal{O}_C))= H^0(C,N).$
\end{defi}

At the points $q_i$ with $i>1$, the fibre of the normal bundle $N$ is isomorphic to the pull back of the relative tangent bundle of $\me/P$ over $\Bbb{P}^1$. At $q_1$ the fibre of the normal bundle has a product structure: It is the product  of the tangent space of $\Bbb{P}^1$ at $p_1$ times the normal bundle of the $T$-fixed curve $C_1$ in the homogeneous space $\me_{p_1}/P=X=G/P$, since we have fixed a trivialization of $\me_{p_1}$.

\begin{defi}\label{defNi}

For each $i=1,\dots, s$, we have Schubert varieties in the fibers of $\mathcal{E}/P$ over $p_i$. Let $T_i $ be their tangent spaces at the points coming from $x$.  These subspaces inject into the fibers of the normal bundle $N$ at the points $q_i$.
(for $i>1$, the injection is immediate, for $i=1$, the injection follows from Lemma \ref{curve}.) Let $N_i\subseteq N_{q_i}$ be the images of $T_i$.
\end{defi}
Now considering Schubert conditions in $\Omega(\vec{w},d)$,
\begin{lemma}\label{newtangent}
The tangent space of the fiber of $\Omega(\vec{w},d)$ over $x$ at $y$ is given by the kernel of the following complex:
$$H^0(C,N)\to \bigoplus_{i=1}^s \frac{N_{q_i}}{N_i}$$
\end{lemma}

\begin{remark}\label{describe}
Let $N_{C_0}$ and $N_{C_1}$ be the normal bundles of the smooth curves $C_0$ and $C_1$ on $\me/P$.
We will use the following description (see \cite[Lemma 2.6]{GHS}) of $N|_{C_1}$ (similarly $N|_{C_0}$):  There is an inclusion $N_{C_1}\subseteq N|_{C_1}$. Finally,  $N|_{C_1}\subseteq N_{C_1}(R)$ is the sheaf of rational sections of $N_{C_1}$ near the node with poles of first order at the node along  normal direction determined by $T(C_0)_R\subseteq (N_{C_1})_R$. Recall that $C_0$ and $C_1$ have distinct tangent  directions at $R$ by assumption.
\end{remark}

\begin{defi}\label{trivio}
Choose a trivialization of $\me$ in an etale (or complex analytic) neighborhood $U$ of $p_1$, so that the curve $C_0\to \me/P =\Bbb{P}^1\times G/P$ is  obtained from $C_0\cong \Bbb{P}^1$ and the constant map to $G/P$ to the point $vP$ and $\bar{g}_1=1$.
\begin{enumerate}
\item[(a)]  This makes the normal bundle of $C_1$ in $\me/P$ a direct sum $N_{C_1} = (f^*TX/TC_1)\oplus \mathcal{O}$.
 \item[(b)] We get a trivialization $N|C_0=(TX)_v\otimes\mathcal{O}$ in a neighbourhood $U$ of $R$: Note that for $x\in U$, $N_x= (TX)_v$, but translation by $v^{-1}$ gives the identification with $(TX)_e$. Also recall that we have a direct sum decomposition
 $$(TX)_e=\bigoplus_{\gamma\in R^+-R_{\frl}^+}\frg_{-\gamma}$$
 The tangent space of $X(v)$ (the $B$-orbit closure of $vP$) corresponds to
 $$v^{-1}TX(v)_v=\bigoplus_{\gamma\in R^+-R_{\frl}^+,v\gamma \prec 0}\frg_{-\gamma}.$$
 which has a distinguished summand $\frg_{-\alpha}$ which corresponds to the image of $(TC_1)_R$.
 \item[(c)] $$N\mid_{C_0} =(\frg_{-\alpha}\otimes\mathcal{O}(R))\oplus \bigoplus_{\gamma\in R^+-R_{\frl}^+,\gamma\neq \alpha}\frg_{-\gamma}\otimes \mathcal{O}$$
 \end{enumerate}
\end{defi}

\subsection{Filtration strategy}
We want to use the filtration strategy to prove Theorem \ref{T2}(a).  The normal bundle agrees on $C_0-\{R\}$ with the pull back of the relative tangent bundle of ${\me/P}$ over $\Bbb{P}^1$ which carries a filtration $\mathcal{F}_i$. But $H^0(C,N)\to H^0(C_0-\{R\},N)$ is not injective, since we have  deformations  of $C_1\subseteq (\me/P)_{p_1}$ which pass through $R$. Therefore we cannot induce a filtration on $H^0(C,N)$ this way. Instead we define a modified normal bundle:

\begin{defi}
Let $N'$ be the kernel of $N\to N_{q_1}/N_1$ (see Defn. \ref{defNi})
\end{defi}
We can then express the  tangent space in Lemma \ref{newtangent} as the kernel of
$$H^0(C,N')\to \bigoplus_{i=2}^s \frac{N_{q_i}}{N_i}.$$

The following justifies the definition of $N'$.
\begin{lemma}\label{restricto1}
$$H^0(C,N')\to H^0(C_0,N'|{C_0})\subseteq H^0(C_0-\{R\},N)$$
is injective.
\end{lemma}
\begin{proof}
Let
$0\to \mathcal{I}\to \mathcal{O}_C\to i_{C_0,*}\mathcal{O}_{C_0}\to 0$ be the exact sequence coming from the inclusion $C_0\subset C$.

In an analytic neighborhood of $R$, we can write $C=\op{Spec}k[x,y]/xy\subseteq \Bbb{A}^2$ where $R=(0,0)$, and $C_0$ is given by $y=0$, and $C_1$ by $x=0$. By a simple computation: $yk[x,y]/xy$ as a $k[x,y]$ is a $k[x,y]/(xy,x)=k[y]$ module identified with $yk[y]$, therefore $I=i_{C_1,*}\mathcal{O}_{C_1}(-R)$. Hence we have an exact sequence
\begin{equation}\label{exacto}
0\to i_{C_1,*}\mathcal{O}_{C_1}(-R)\to\mathcal{O}_C\to i_{C_0,*}\mathcal{O}_{C_0}\to 0.
\end{equation}

We now begin the proof of Lemma \ref{restricto1}. Using the exact sequence \eqref{exacto},  it suffices to show that $H^0(C_1,N'|C_1(-R))=0$.  Using the description of $N$ given in Remark \ref{describe}, and Dfn. \ref{trivio}
$N_{C_1} = (f^*TX/TC_1)\oplus \mathcal{O}$, and  $(TC_0)_R$ gives a line $\Bbb{C}[0\oplus 1]$ in the fiber of $(f^*TX/TC_1)\oplus \mathcal{O}$ at $R$, and therefore
 $$N|{C_1}=(f^*TX/TC_1)\oplus \mathcal{O}(R).$$

Sections of  $H^0(C_1, N'|C_1(-R))$  are of the form $a\oplus x\subseteq N_{C_1}$ so that  $x$ vanishes at $q_1$, and $x$ is regular at $R$ making $x$ zero. Also $a\in H^0(C_1, ((f|C_1)^*TX/TC_1)(-R))$ with fiber at  $q_1$ in $(TX_w)_{q_1}$: But this is the tangent space of deformations of $f(C_1)$ which pass through $f(R)$ and intersect $X_w$. By assumption there are no such deformations, so $a=0$, and this finishes the proof.
\end{proof}

\subsection{}

\begin{defi} Define $J\subseteq  (N|C_0)_R$ to be the subspace fomed by elements $\delta$, so that $\delta$ extends to a section of $H^0(C_1,N'|C_1)$. Note that a section $s\in H^0(C_0,N|C_0)$ extends to a section of $H^0(C,N')$ if and only if $s_R\in J$.
\end{defi}
We observe that $H^0(C,N')$ is the kernel of the map
$H^0(C_0,N|C_0)\to (N\mid C_0)_R/J$. Therefore
\begin{lemma}\label{newmap}
We can express the  tangent space in Lemma \ref{newtangent} as the kernel of
\begin{equation}\label{Map0}
H^0(C,N\mid C_0)\to N_R/J\oplus \bigoplus_{i=2}^s \frac{N_{q_i}}{N_i}
\end{equation}
\end{lemma}

\begin{remark}\label{newmap1}
Next we induce filtrations on the objects appearing in Lemma \ref{newmap}:
\begin{enumerate}
\item We will induce a filtration $\mathcal{F}_{\bull}(N|C_0)$ by sub-bundles on the vector bundle $N|C_0$ (=$N'|C_0$) on $C_0$. A section is in the $j$th filtered piece, if its restriction to $C_0-\{R\}$ is in the $jth$ filtered piece as a section of $N_{C_0}$. Note that in an etale neighborhood of $R$ in  $C_0$, $N_{C_0}$ is isomorphic to $f^*TX$, using a trivialization of $\me/P$ near $p_1$ as in Dfn. \ref{trivio}.
\item The above step  induces a filtration of the fibre of $N|C_0$ at $R$.
\end{enumerate}
\end{remark}
\subsection{}
Let $J\subseteq (N|C_0)_R=(N'|C_1)_R$ be the subspace fomed by elements $\delta$, so that $\delta$ extends to a section of $H^0(C_1,N'|C_1)$. Clearly
$s\in H^0(C_0,N|C_0)$ extends to a section of $H^0(C,N')$ if and only if $s_R\in J$.
Furthermore let $\widetilde{N}\subset N|C_0$ formed by sections of $N$ whose fibre at $R$ lies in $J$. Clearly
$$H^0(C,N')=H^0(C_0,\widetilde{N})\subseteq  H^0(C_0,N'|{C_0})\subseteq H^0(C_0-\{R\},N).$$
\begin{lemma}
$\dim J=\dim X_v$.
\end{lemma}
\begin{proof}
$J$ is the set of elements in $(N'|C_1)_R$ which extend to a section of $H^0(C_1,N'|C_1)$.  The map
$H^0(C_1,N'|C_1)\to (N'|C_1)_R$ is injective by the proof of Theorem \ref{restricto1}. Therefore
$\dim J= H^0(C_1,N'|C_1)$.

Using the trivialization of $\me$ at $p_1$ in Dfn. \ref{trivio}, $N|C_1=f^*(TX/TC_1)\oplus \mathcal{O}_{C_1}(R)$.
Since $f^*(TX/TC_1)$ is globally generated on $C_1$, the map $H^0(C,N|C_1)\to N_q/N_1$ is surjective. 
We also know that $H^1(\Bbb{P}^1,N|C_1)=0$ which implies $H^1(\Bbb{P}^1,N'|C_1)=0$ and hence 
(using $\chi ({P}^1,\mv)=\rk \mv +\deg\mv $)
$$H^0(N'\mid C_1)= \chi (N'|C_1)= (\dim X -1 +1) +(\int_{d(\alpha)}c_1(TX)-2 +1) -\codim X_w.$$
But since the stratum considered is a quantum-type stratum
$\codim X_w-\codim X_v+1 =\int_{d(\alpha)}c_1(TX))$, so that $\dim H^0(C_1,N'|C_1)=\dim X_v$.
\end{proof}

 \begin{lemma}\label{node}
 $J\subseteq (N\mid_{C_0})_R$ is the fiber at $R$  of the subbundle
 $$(\frg_{-\alpha}\otimes\mathcal{O}(R))\oplus \bigoplus_{\gamma\in R^+-R_{\frl}^+,v\gamma \prec 0,\gamma\neq \alpha}\frg_{-\gamma}\otimes\mathcal{O}$$

 \end{lemma}
\begin{proof}
First note that we know lots of elements in $J$: Any deformation of $C_0$ in a neighborhood of $0$ so that $0$ stays in $X(v)$ lies in $J$ (because we can compatibly deform $C_1$ (by Lemma \ref{stringy}). This gives a map $(TX(v))_v\to J$, which has a one dimensional kernel in direction of $C_1$. So there is ``one missing direction'' in $J$ that needs to be determined.

We note that $J$ can be computed locally from a neighborhood of $C_1$ in $C_0\cup C_1$. Therefore it suffices to consider an example to prove the statement:

Let $\Bbb{A}^2=\Bbb{A}^1_t\times \Bbb{A}^1_z$ be the set of ordered pairs $(t,z)$. Blow up $(0,0)$ to get the surface $\widetilde{\Bbb{A}^2}$. Let $\tau:\widetilde{\Bbb{A}^2}\to \Bbb{A}^1_t$.
Now for $t\neq 0$, $\tau^{-1}(t)$ is $\Bbb{A}^1$ (parameterized by coordinate $z$), and $\tau^{-1}(0)=C_0\cup C_1$, where $C_1$ is the exceptional divisor and $C_0$ the strict transform of the $z$-axis.

Consider the rational map $\pi:\widetilde{\Bbb{A}^2}\to X$, given  by
$$(t,z)\mapsto wU_{-\alpha}(z/t)eP$$
We will show that $\pi$ is regular: In the coordinate patch $(a,t)$ of $\widetilde{\Bbb{A}^2}$  with $z=at$, the map is a well defined $wU_{-\alpha}(a) eP$ which lies on the curve $C(v,w)$. On the coordinate patch $(b,z)$ of $\widetilde{\Bbb{A}^2}$ with $t=bz$, the map takes the form $wU_{-\alpha}(1/b)eP$ which extends to $b=0$ with image $vP$ at $b=0$. The exceptional divisor maps to the curve $C(v,w)$.

The point $(0,t)$ for $t\neq 0$ maps to $\dot{w}$, The node point $C_0\cap C_1$ maps, under $\pi$, to $\dot{v}$. The curve $C_0$ also maps to the point $\dot{v}$, so the corresponding principal bundle on $C_0$ is trivial.

The curve $D_0$ on $\widetilde{\Bbb{A}^2}$ has two components: The exceptional divisor, and the strict transform of $t=0$. The strict transform meets the exceptional divisor at the point on the second patch with $b=z=0$, which goes to $\dot{v}$ under the map described.  Therefore $D_0$ is our curve $C$. The curves $D_t$ for $t\neq 0$ are smooth, and pass through $\dot{w}$ at $z=0$. Hence they give the missing ``deformation''. We can write the map
$\pi$ for $z\neq 0$ as $vU_{-\alpha}(t/z)eP$ and hence the vector field  produced by the deformation on $C_0-R$ has a pole at $z=0$, and gives the term $\frg_{-\alpha}\otimes\mathcal{O}(R)$ in (c) of Dfn. \ref{trivio}.
\end{proof}

\begin{figure}

\begin{tikzpicture}

\foreach \y in {5} {

\draw[-,thick,orange] (-3,\y) -- (0,\y);  
\draw[-,thick,blue] (-3,\y+1) -- (-0.5,\y+1);
\draw[-,thick] (-2.5,\y-1) -- (0,\y-1);
\draw[-,thick,blue] (-2.8,\y-0.7) -- (-2.8,\y+1.2);
\draw[-,thick] (-0.2,\y-1.2) -- (-0.2,\y+0.7);

\draw[-,thick,red] (-2.3,\y-1.1) -- (-2.9,\y-0.4);
\draw[-,thick,red] (-0.1,\y+0.4) -- (-0.7,\y+1.1);

\draw[->] (-1.5,\y-2) -- (-1.5, 2);

\draw[->] (0.7,\y) -- (2,\y);

\node at (3,\y) {$X = G/P$};

\node at (-3.1,\y+0.5) {\tiny $w$};
\node at (-1.5,\y+1.2) {\tiny $w$};

\node at (-1.5,\y-1.2) {\tiny $v$};
\node at (0.1,\y-0.5) {\tiny $v$};

\node at (-1.5,\y+0.2) {\tiny $C_{v,w}$};
\node at (-2.25,\y-0.7) {\tiny $C_{v,w}$};
\node at (-0.75,\y+0.7) {\tiny $C_{v,w}$};

\node at (-3.5,\y-1.5) {\tiny $C_1$};
\node at (-0.5,\y-1.6) {\tiny $C_0$};

\draw[->] (-3.3,\y-1.3) -- (-2.7,\y-0.8); 
\draw[->] (-0.8,\y-1.4) -- (-1.1,\y-1.2); 

\node at (-4.3,\y) {$\tilde P$};

\node at (0.3,\y) {\tiny $D_t$};

}

\draw[-,thick,orange] (-3,0) -- (0,0);
\draw[-,thick,blue] (-3,1) -- (0,1);
\draw[-,thick] (-3,-1) -- (0,-1);
\draw[-,thick,blue] (-2.8,1.2) -- (-2.8,-1.2);
\draw[-,thick] (-0.2,1.2) -- (-0.2,-1.2);

\draw[->] (0.7,0) -- (2,0);

\draw[-,thick] (2.5,-1.2) -- (2.5,1.2);
\draw[-,thick] (2.3,0) -- (2.7,0);

\node at (0.3,0) {\tiny $D_t$};
\node at (3,0) {\tiny $t$};
\node at (2.5,1.6) {$\mathbb{P}^1_t$};

\node at (-2.8,-1) {\color{red} $\times$};
\node at (-0.2,1) {\color{red} $\otimes$};

\node at (0.4,1.3) {\color{red} \tiny $(\infty, \infty)$};

\node at (-3.2,-1.3) {\color{red} \tiny $(0,0)$};

\node at (-4.4,0) {$P = \mathbb{P}^1\times \mathbb{P}^1$};

\draw[->] (-1.5,-2) -- (-1.5,-3);
\node at (-1.5,-4) {$\mathbb{P}^1_z$};

\end{tikzpicture}

\caption{A smoothing of $D_0 = C_0\cup C_1$ to $D_t\simeq \mathbb{P}^1$ for $t\ne 0$. } \label{smoothing}
\end{figure}
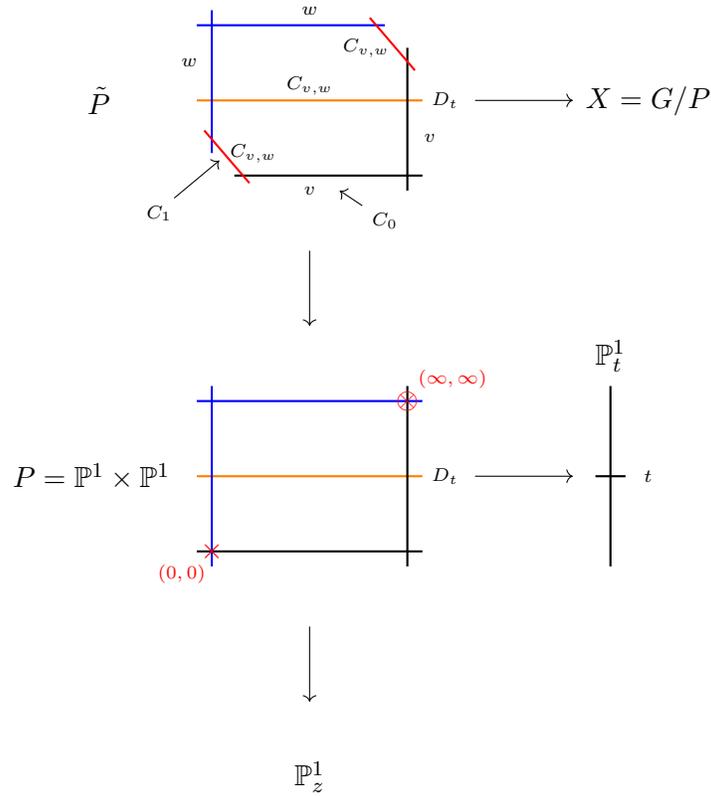
\begin{remark}
We can complete the family parametrized by $t\in \Bbb{A}^1$ to a family of complete curves parametrized by $t\in\Bbb{P}^1$ which is symmetric.
Let $P=\Bbb{P}^1_t\times \Bbb{P}^1_z$, and $\widetilde{P}$ the blow-up of $P$ at the points $(0,0)$ and $(\infty,\infty)$. 
Consider the meromorphic map from $P$ to $X$
$$(t,z)\mapsto wU_{-\alpha}(z/t)eP =vU_{-\alpha}(-t/z)eP\in X$$
This extends to a regular map $\tau:\widetilde{P}\to X$.

For $t\neq 0,\infty$, the curve $D_t=\{(z,t):z\in \Bbb{P}^1\}=\Bbb{P}^1$ maps to $C(v,w)$ with $0$ mapping to $wP$ and $\infty$ to $vP$. As we let
$t\to 0$ or $t\to \infty$ we need to form flat limits. The two limits of $D_t$ at $t=\infty$ and $t=0$ are shown in Figure \ref{smoothing}.

\end{remark}
\begin{remark}\label{smoothening}
This deformation can be used to smooth out a general point of quantum strata in Definition \ref{Maxim}. Suppose we have a principal $G$ bundle $\me$ on $\Bbb{P}^1_z$ with a given map $F:\Bbb{P}^1_z\to \me/P$. Since $P$ bundles on $\Bbb{A}^1_z$ are trivial, we may assume that we have a trivialization $e$ of $\me$ so that $F(z) = evP$ for all $z$.
Also choose a trivialization $\tilde{e}$ of $\me$ in a neighborhood $\Bbb{P}^1-\{0\}$ of $\infty$ so that $ \tilde{e}g_0(z)=e$, so that the map $F$ looks like $F(Z)=\tilde{e}P$. Hence
$g_0(z)v\in P((Z))$ where $Z=1/z$ is a local coordinate at $\infty\in\Bbb{P}^1_z$

We therefore get a  map  $\widetilde{\Bbb{A}^2}\to G/P$ with $(0,0)$ blown up as above. This is the data of a trivial principal bundle on $\widetilde{\Bbb{A}^2}$ and a $P$ reduction. We want to extend this principal bundle  to a principal $G$-bundle $\me^e$ on $\Bbb{P}^1_z\times \Bbb{A}^1_t$  and the reduction to a map 
\begin{equation}\label{mappo}
\widetilde{\Bbb{P}^1_z\times\Bbb{A}^1_t}\to \me^e/P
\end{equation}
 so that the restriction to $t=0$ coincides with the initial data on $C_0$.  The principal bundle $\me^e$ is obtained by gluing together the trivial bundle on $\Bbb{A}^1_z\times \Bbb{A}^1_t$ with the trivial bundle on $A^1_Z\times \Bbb{A}^1$. Let the transition be $g_t(z)=  g_0(z)v(vU_{-\alpha}(-t/z))^{-1}$ which at $t=0$ is $g_0(z)$. The map $F_t$ is $e vU_{-\alpha}(-t/z)P= \tilde{e} g_t(z) vU_{-\alpha}(-t/z)P= \tilde{e} g_o(z)vP=\tilde{e}P$. Thus $F_t$ extends to a map $\widetilde{\Bbb{P}^1_z\times\Bbb{A}^1_t}\to \me^e/P$. 

We could also consider this as follows. The principal bundle on $\Bbb{C}^*_z\times \Bbb{A}^1_t$ coming from \eqref{mappo} is trivial (since we have assigned an element in the $P$ coset). We also have a trivial $P$-bundle in the neighborhood of $\infty$, and at $t=0$, we are given a patching function 
$g_0(z)v\in P((Z))$. We can use the same patching function for all $t$ and hence get a $P$-bundle on $\widetilde{\Bbb{P}^1_z\times\Bbb{A}^1_t}$ as desired.
\end{remark}

\subsection{Tangent spaces}\label{class}

As in \cite{BKR}, the Levi-movability condition implies that both sides of \eqref{tsp0} are filtered, and \eqref{tsp0} is a filtered isomorphism. The tangent space $T\overline{g}_i X(w_i)$ is taken at $f(p_i)\in \overline{g}_i X(w_i)$ which is in the open cell.

The tangent bundle $f^*T(G/P)$ is filtered by subbundles coming from a canonical filtration of the tangent bundle of $G/P$ at the identity by weights of the $T$-action.
$$\mathcal{F}_1\subset\dots\subset\mathcal{F}_{r}$$

Let us compute the assigned number $\deg$ from Section \ref{gradedDim}  for the situation in Remark \ref{newmap1}.  Our aim is that the $\deg$ here is greater than the $\deg$ assigned to the map \eqref{tsp0} which is zero by the Levi-movability assumption. This will prove Theorem \ref{T2}(a).

\subsection{ Computation of $\deg$ in \eqref{tsp0}}
Let $r$ be the length of the filtration:
$$0\subseteq\mf_1\subseteq \mf_2\subseteq \dots\subseteq \mf_r=\mf$$
For simplicity introduce a weight $\tau$ so that $\tau(x_P)=r+1$.

Let us start with
$$\sum_j\chi(\Bbb{P}^1,\mf_j)=\sum_j\rk\mf_j +\sum_j\deg\mf_j$$
Now $$\sum_j\rk\mf_j=\sum_{j}|S_j|=\sum_{\gamma\in R^+-R_{\frl}^+}(\tau-\gamma)(x_P)$$
Writing $\deg\mf_j= \sum_{\gamma\in R^+-R_{\frl}^+,\gamma(x_P)\leq j} \gamma(d)$, where $d=\sum a_i \alpha_i^{\vee}$, so that
$$\sum_{j}\chi(\mf_j)=\bigl(\sum_{\gamma\in R^+-R_{\frl}^+}\gamma(d)(\tau-\gamma) +(\tau-\gamma\bigr)(x_P))$$
The key point is that there is no summation in $j$ on the right hand side.
Note that $$\sum_{i=1}^s\sum_j \bigl(\dim (\mathcal{F}_j)_{p_i}-\dim (\mathcal{F}_j)_{p_i}\cap B_i)\bigr)$$
has the form
$$\sum_{i=1}^s \sum_{\gamma\in R^+-R_{\frl}^+,w_i\gamma\succ 0} (\tau-\gamma)(x_P)$$

Therefore the assigned number $\deg$ for the situation in \eqref{tsp0}, known to be equal to zero is the following element of $\frh^*$ evaluated at $x_P$.
$$\sum_{\gamma\in R^+-R_{\frl}^+}(\gamma(d)(\tau-\gamma)+(\tau-\gamma))-(\sum_{i=1}^s \sum_{\gamma\in R^+-R_{\frl}^+,w_i\gamma\succ 0}(\tau-\gamma)).$$
By a simple computation (using Lemma \ref{node}), we get
\begin{lemma}
The assigned number $\deg$ above for the map in Lemma \ref{newmap} is
$$\tau-\alpha +\sum_{\gamma\in R^+-R_{\frl}^+}\gamma(d-\alpha^{\vee})(\tau-\gamma)-(\sum_{i=2}^s \sum_{\gamma\in R^+-R_{\frl}^+,w_i\gamma\succ 0}(\tau-\gamma) +\sum_{\gamma\in R^+-R_{\frl}^+,v\gamma\succ 0}(\tau-\gamma))$$
\end{lemma}

To prove Theorem \ref{T2}(a), by comparing the $\deg$ for Lemma \ref{newmap} with the $\deg$ from \eqref{tsp0} computed above (assumed zero), we need to show that  the following element of $\frh^*$ evaluated at $x_P$ ro be $>0$,
$$(\tau-\alpha)+\sum_{\gamma\in R^+-R_{\frl}^+}(\gamma(-\alpha^{\vee})(\tau-\gamma)-\bigl(\sum_{\gamma\in R^+-R_{\frl}^+,v\gamma\succ 0}(\tau-\gamma) -\sum_{\gamma\in R^+-R_{\frl}^+ w\gamma \succ 0}(\tau-\gamma)\bigr)$$
this is the same as
$$(\tau-\alpha)+ \sum_{\gamma\in R^+-R_{\frl}^+}(\gamma(-\alpha^{\vee})(\tau-\gamma))-\big(\sum_{\gamma\in R^+-R_{\frl}^+,w\gamma\prec 0}(\tau-\gamma) -\sum_{\gamma\in R^+-R_{\frl}^+,v\gamma \prec 0}(\tau-\gamma)\bigr)$$
The coefficient of $\tau$ is  $-\int_{d(\alpha)}c_1TX+1-\ell(w)+\ell(v)=0$, therefore we are reduced to showing that the following quantity evaluated at $x_P$ is $<0$:
$$\alpha- \sum_{\gamma\in R^+-R_{\frl}^+}\gamma(\alpha^{\vee})\gamma-(\sum_{\gamma\in R^+-R_{\frl}^+,w\gamma\prec 0}\gamma -\sum_{\gamma\in R^+-R_{\frl}^+,v\gamma \prec 0}\gamma)$$
Using Lemma \ref{anotherlemma?},  the above quantity equals the following (evaluated at $x_P$)
\begin{equation}\label{formalsum}
\alpha-\sum_{\gamma\in R^+}\gamma(\alpha^{\vee})\gamma-(\sum_{\gamma\in R^+,w\gamma\prec 0}\gamma \   \ -\  \  \  \  \sum_{\gamma\in R^+,v\gamma \prec 0}\gamma)
\end{equation}
(using that $w\in W^P$)

The quantity \eqref{formalsum} equals
$$\alpha+ w^{-1}\rho -w^{-1}\rho+(w^{-1}\rho({\alpha}^{\vee}))\alpha +\alpha -\sum_{\gamma\in R^+}\gamma(\alpha^{\vee})\rangle\gamma$$
 which equals
 \begin{equation}\label{simplified}
 (\rho(w\alpha^{\vee})+1)\alpha- \sum_{\gamma\in R^+}\gamma({\alpha}^{\vee})\gamma
 \end{equation}
 We need to  evaluate \eqref{simplified} at $x_P$, and determine it to be a quantity which  is $\leq 0$ and $0$ only if $w\alpha=\theta$. Now write
$$\rho(w{\alpha}^{\vee})=\frac{2(\rho,w\alpha)}{(w\alpha,w\alpha)}=\frac{2(\rho,w\alpha)}{(\alpha,\alpha)}$$
and using  \cite[Lemma 3.4]{BKq} (see Remark \ref{BCbelow})
\begin{equation}\label{theSum}
\sum_{\gamma\in R^+}\gamma ({\alpha}^{\vee})\gamma=\frac{2g^*}{(\alpha,\alpha)} \alpha
\end{equation}
where $g^*=\rho(\theta^{\vee})+1=(\rho,\theta)+1.$

Therefore the sum \eqref{simplified}, evaluated on $x_P$,  equals $\alpha(x_P)$ times (note that $\alpha(x_P)>0$ since $\alpha\in R^+\setminus R_{\frl}^+$)
$$
\frac{2}{(\alpha,\alpha)} (\rho,w\alpha-\theta)+(1 -\frac{(\theta,\theta)}{(\alpha,\alpha)}),
$$
 a number that is strictly negative unless $w\alpha= \theta$. Here we use that the highest root is long, and hence $\frac{(\theta,\theta)}{(\alpha,\alpha)}\geq 1$. We also know that $\theta-w\alpha$ is a sum of positive roots, and hence $(\rho,w\alpha-\theta)< 0$ unless $w\alpha=\theta$.
\begin{remark}\label{BCbelow}
Recall from \cite[Lemma 3.4]{BKq} (also \cite[Equation (7.9)]{Mac} and \cite[Corollary 1.3]{BC}) that for $y\in \frh$,
$$
\sum_{\gamma\in R^+} \gamma(\alpha^{\vee})\gamma(y) = g^*(\alpha^{\vee},y) =\frac{2g^*\alpha(y)}{(\alpha,\alpha)}
.$$
\end{remark}

\section{Cycle classes}\label{cycleclasses}
Let $v_1,\dots,v_s\in W^P$ and $d\in H_2(G/P)$ (possibly $d\leq 0$). Now assume that
$$\sum \codim \sigma^P_{v_i} = \dim X +\int_d c_1(T_X)+1.$$
We let $\mathcal{D}_{\vec v}=f_*[\Omega(\vec{v},d)]$ be the divisor class on $\op{Parbun}_{G}$ obtained by the cycle-theoretic pushforward of the cycle class of $\Omega(\vec{v},d)$.

We want to compute this cycle class $\mathcal{D}_{\vec v}$. Let the corresponding line bundle be $\mathcal{B}(\vec\lambda,\ell)=\mathcal{O}(\mathcal{D}_{\vec v})$. Then we have the following formulas for finding $\lambda_1, \hdots, \lambda_s$ and $\ell$. Write $\lambda_i = \sum c_i^b \omega_b$, where $\omega_b$ is the $b^\text{th}$ fundamental weight.

\begin{theorem}\label{firstformula}
\begin{enumerate}
\item[(a)] Set $v_{s+1} = \overline{s_{\theta}w_0}$, and $m = d(\theta)\in H_2(G/P)$. Then we have the ``level formula":
$$
\ell = \langle v_1,\hdots, v_s, v_{s+1}\rangle_{d+m}
$$
\item[(b)] $c_i^b=0$ if $s_{b}v_i$ has smaller length than $v_i$ or if $s_{b}v_i\not\in W^P$. Otherwise, $v_i\xrightarrow{\alpha_b}s_bv_i\in W^P$ and we define $w_k = v_k$ for $k\ne i$ and $w_i = s_bv_i$. Then
$$
c_i^b = \langle w_1,\hdots,w_s\rangle_{d}.
$$
\end{enumerate}
\begin{remark}
If $d< 0$, these formulas show that  $\mathcal{D}_{\vec v}$ is a multiple of the pullback of the determinant of cohomology divisor on $\Bun_G$ (with complement the points where the underlying bundle is trivial).
\end{remark}

\end{theorem}

\subsection{Proof of Theorem \ref{firstformula} (b)}\label{firstproof}

We have a map $i_{\alpha}:\Bbb{P}^1\to G/B$ corresponding to the simple root $\alpha=\alpha_i$ whose image is a $T$-fixed curve. This embedding takes $t\in \Bbb{A}^1$ to
$U_{\alpha}(t)\dot{s}_{\alpha}$. The pullback of a line bundle $\ml=\sum a_j\omega_j$ via this map is $\mathcal{O}(a_i)$.

Define the smooth map $\pi:\Bbb{P}^1\times (G/B)^{s-1}\to \Parbun_G$   by
$(t,\bar{g}_2,\dots,\bar{g}_s)\mapsto (\me,i_{\alpha}(t),\bar{g}_2,\dots,\bar{g}_s)$ where $\me$ is the
trivial $G$-bundle. Consider $\pi^{*}(\mathcal{D}_{\vec v})\subseteq \Bbb{P}^1\times (G/B)^{s-1}$. By Lemma \ref{elementary} below, $c_i^b$ equals the degree of
the map $\Omega(\vec{v},d)'\to (G/B)^{s-1}$ (where $\Omega(\vec{v},d)'$ is the base change of $\Omega(\vec{v},d)$). Fix a general point $(\bar{g}_2,\dots,\bar{g}_s)\in (G/B)^{s-1}$. Then,
\begin{itemize}
\item $c_i^b$ is the number of pairs $(f,t)$ where $f:\Bbb{P}^1\to G/P$  and $t\in \Bbb{P}^1$ so that $f(p_j)\in {g}_i C^P_{v_j}$  for $j\geq 2$ and $f(p_1)\in i(t)C^P_{v}$.
\end{itemize}

As in \cite[Section 4.3]{BKiers} if $\alpha\in \Delta(Q_u)$ (with the notation of $\Delta(Q_u)$ as in loc. cit.) the closures of $i(t)C^P_{v}=s_{\alpha}C^P_v$ are constant, and the count above is zero for codimension reasons. This shows that unless $v^{-1}\alpha\in R^+-R_{\frl}^+$, the count is zero.

Now by  \cite[Lemma 4.4]{BKiers} , $v^{-1}\alpha\in R^+-R_{\frl}^+$ if and only if  setting $w=s_{\alpha}v=s_bv$, we have $w\in W^P$ and $v\xrightarrow{\alpha}w$.  We therefore assume that this is the case. We claim that this count is then the same as  $\langle w_1,\hdots,w_s\rangle_{d}$ (evaluated for the trivial bundle with flags $(\dot{e},\bar{g}_2,\dots,\bar{g}_s)$ which follows from (see \cite[Lemma 4.2]{BKiers}):
\begin{lemma}
Let $x\in G/P$. The following are equivalent
\begin{enumerate}
\item $x\in C^P_{w}$.
\item There exists a $t\in \Bbb{P}^1$ so that $x\in  i(t)C^P_{v}$.
\end{enumerate}
(When  $t$ exists in (2), it is unique.)
\end{lemma}

\begin{lemma}\label{elementary} Let $Z$ be a smooth complex algebraic variety.
Let  $\pi:\Omega\to \Bbb{P}^1\times Z$ be a proper morphism with $\Omega$ equidimensional and generically reduced, and $\dim\Omega=\dim Z$. Let $\mathcal{O}(\pi_*[\Omega])=\mathcal{O}(\ell)\boxtimes \mathcal{L}$. Then $\ell$ is the degree of the morphism $\Omega\to Z$ (i.e., cardinality of a general fiber).
\end{lemma}

\subsection{A generalization}
We prove the following generalisation of Theorem \ref{firstformula}(a): We work with $s$ points (one of them $p_{1}=0\in\Bbb{P}^1$). Let $v_1,\dots,v_s\in W^P$ and $d\in H_2(G/P)$ be such that
$$\sum_{i=1}^{s} \codim \sigma^P_{v_i} = \dim X +\int_d c_1(T_X)+1.$$

Let $\mathcal{B}(\vec\lambda,\ell)=\mathcal{O}(\mathcal{D}_{\vec v})$, a line bundle on $\Parbun_{G,S}$.
Let $v'_{1}$ be the minimal coset representative for $s_{\theta}v_{1}W_P$.
 Then,
\begin{theorem}\label{generale}
\begin{enumerate}
\item  If $\ell(v'_{1})\geq \ell (v_{1})$, then $\ell-\lambda^{1}(\theta^{\vee})=0$.
\item  If $\ell(v'_{1})< \ell (v_{1})$, then $\alpha=-v_{1}^{-1}\theta\succ 0$. Set  $m=d(\alpha)$. Then $m> 0$ and
$$\ell-\lambda^{1}(\theta^{\vee}) = \langle v'_1,v_2,\hdots, v_s\rangle_{d+m}
$$
(the GW invariant above will be shown to have expected dimension zero).
\end{enumerate}
\end{theorem}
Note that when $v_{1}=\overline{w_0}$, $\lambda_1=0$,  and (enlarging the set $S$ by adding a new point, and relabelling so that the new point is $p_1$) we obtain the level formula in Theorem \ref{firstformula}.
\begin{remark}In case (2), 
 $\Omega^0(\vec{v},d)$
corresponds  a quantum type stratum for $\Omega((v'_1,v_2,\dots,v_s),d+m)$ (see Defn. \ref{Maxim}).
\end{remark}
To see that $-v_{1}^{-1}\theta\succ 0$, and $m>0$ in Theorem \ref{generale} (2), it suffices to observe that $\ell(s_\theta v_1)  \le \ell(v_1)$, since otherwise we have $\ell(v_1')\ge \ell(v_1)$.  Thus $v_1^{-1}\theta \prec 0$. Since $\alpha = -v_1^{-1}\theta \succ 0$, $m = d(\alpha)> 0$.  

\begin{lemma}\label{fraction}
Suppose $u\in W^P$ and let $u'$ be the minimum coset representative for $s_\theta u W_P$. 
\begin{enumerate}
\item  $\ell(u')= \ell(u)$, iff $u^{-1}\theta$ belongs to $R_{\mathfrak{l}}$. In this case $u=u'$.
\item  $\ell(u')>\ell(u)$  iff $u^{-1}\theta\in R^+-R_{\frl}^+$. In this case $u'$ corresponds to a codimension one quantum stratum of $u$ and
$ {\ell(u')-\ell(u)+1}=\int_{d(u^{-1}\theta)}c_1(TX)$.

\item $\ell(u')<\ell(u)$ iff $u^{-1}\theta\in R^- -R_{\frl}^+$. In this case $ {\ell(u)-\ell(u')+1}=\int_{d(-u^{-1}\theta)}c_1(TX)$.
\end{enumerate}
\end{lemma}
\begin{proof}
First observe the symmetry in $u$ and $u'$: Suppose $u,u'\in W^P$. Then  $u'$ is the minimum coset representative for $s_\theta u W_P$ if and only if
$u$ is the minimum coset representative for $s_\theta u' W_P$. Also note that $u^{-1}\theta\succ 0 \iff u'^{-1}\theta\prec 0$. 

By symmetry,  it suffices to prove that
\begin{enumerate}
\item[(A)] If $u^{-1}\theta\in R_{\frl}$ then $u=u'$, and hence $\ell(u)=\ell(u')$.
\item[(B)] If $u^{-1}\theta\in R^+- R_{\frl}^+$, then $ {\ell(u')-\ell(u)+1}=\int_{d(u^{-1}\theta)}c_1(TX)$ and hence $\ell(u')>\ell(u)$.
\end{enumerate}
To see (A), if $u^{-1}\theta\in R_{\frl}$  then
$$
s_\theta u W_P = uu^{-1}s_\theta u W_P = u s_{u^{-1}\theta} W_P = u W_P
$$
because $s_{u^{-1}\theta}\in W_P$  which gives $u=u'$ as desired.

For (B), assume $u^{-1}\theta\in R^+- R_{\frl}^+$ and let $\alpha=u^{-1}\theta$. The statement follows directly from Lemma \ref{old8-7}, with $v = u'$ and $w = u$. 
\end{proof}

\section{Outline of the proof of Theorem \ref{generale}}\label{outline}
The proof will be modelled after the proof of Theorem \ref{firstformula} (a) given in Section \ref{firstproof}.
\subsection{Step 1}
We construct a smooth morphism
\begin{equation}\label{etape1}
\pi:\Bbb{P}^1\times (G/B)^{s-1} \to   \Parbun_{G,S}.
\end{equation}
so that the pull back of  $\mathcal{B}(\vec{\lambda},\ell) =\mathcal{O}(\mathcal{D}_{\vec v})$ equals
$\mathcal{O}_{\Bbb{P}^1}(\ell-\lambda_{1}(\theta^{\vee}))\boxtimes \boxtimes _{i=2}^{s} L(\lambda_i)$.
Therefore by Lemma \ref{elementary}, $\ell-\lambda_{1}(\theta^{\vee})$ is the degree of the map
$\Omega(\vec{v},d)'\to (G/B)^{s-1}$, which is the same as the degree of $\Omega^0(\vec{v},d)'\to (G/B)^{s-1}$, where $\Omega(\vec{v},d)'$ (resp. $\Omega^0(\vec{v},d)'$) is the base change of $\Omega(\vec{v},d)$ (resp.  $\Omega^0(\vec{v},d)$) under \eqref{etape1}.
\subsection{Step 2}
Suppose $\ell(v'_{1})\geq \ell (v_{1})$.  This implies $v_1^{-1}\theta$ is in ${R}_{\frl}$ in which case $v_1=v_1'$ or
$v^{-1}\theta\in R^+-R_{\frl}^+$. We will show that the fiber of  $\Omega^0(\vec{v},d)'\to (G/B)^{s-1}$ over a general point of
$(G/B)^{s-1}$ is empty, and hence show that $\ell-\lambda_{1}(\theta^{\vee})=0$. This step is carried out in Section \ref{secondstep}.

\subsection{Step 3} Supposing $\ell(v'_{1})< \ell (v_{1})$, we show  that the fiber of $\Omega^0(\vec{v},d)'\to (G/B)^{s-1}$ over a general point of
$(G/B)^{s-1}$ is in bijection with the enumerative count  $\langle v'_1,v_2,\hdots, v_s\rangle_{d+m}$ (carried out at the parabolic bundle corresponding via \eqref{etape1}, to the chosen point of $(G/B)^{s-1}$ 
and $t=0\in \Bbb{P}^1$). This will conclude the proof of the theorem. This step is carried out in Section \ref{Step3proof}

\section{Moduli stacks and uniformizations}\label{Stepp1}
\begin{defi}
\begin{enumerate}
\item Let $\mathcal{Q}_G=G(\Bbb{C}((z)))/G(\Bbb{C}[[z]])$ be the affine Grassmannian parameterizing
principal $G$-bundes $E$ on a curve $C$ with a trivialization  $E\big|_U\to U\times G$ over $U=C-p$. Here $z$ is a local coordinate at $p$ (we will soon let $C=\Bbb{P}^1$ and $p=0$).

\item Let $\mathcal{I}_G=G(\Bbb{C}((z)))/I_G$ be the Iwahori Grassmannian, with $I_G$ the inverse image of $B$ under the evaluation at zero morphism
$G(\Bbb{C}[[z]])\to G$. Note that $\mathcal{I}_G$ parametrizes principal $G$ bundes on a curve $C$ with a trivialization  $E\big|_U\to U\times G$ over $U=C-p$, and an element $\bar{g}_0\in E_0/B$.

\item Let $\Parbun^0_G$ be the moduli stack of principal $G$ bundles on $\Bbb{P}^1$ with $B$-reductions at $0\in \Bbb{P}^1$.
\end{enumerate}
\end{defi}
Consider the embedding  $SL_2\to G$ given by the highest root $\theta$. Let
$\phi:SL_2(\Bbb{C})\to SL_2(\Bbb{C}((z)))$ be the homomorphism
$$\begin{pmatrix}
a & b \\
c & d
\end{pmatrix}\longmapsto \begin{pmatrix}
d & cz^{-1} \\
bz & a
\end{pmatrix}$$
Letting $B_{SL_2}$ be the group of upper triangular matrices in $SL_2(\Bbb{C})$, we get morphisms
\begin{equation}\label{long}
\Bbb{P}^1= SL_2(\Bbb{C})/B_{SL_2}\to SL_2(\Bbb{C}((z)))/I_{SL_2}=\mathcal{I}_{SL_2}\to \mathcal{I}_G
\end{equation}
Since principal bundles parametrized by $\mathcal{I}_G$ are trivialized outside of zero, we have a map
$\mathcal{I}_G\times (G/B)^{s-1}\to \Parbun_G$, which when composed with \eqref{long} gives (this gives the morphism \eqref{etape1}):
\begin{equation}\label{mapPar}
\Bbb{P}^1\times (G/B)^{s-1}\to \Parbun_G.
\end{equation}
\begin{proposition}\label{properties}
\begin{enumerate}
\item The map $\eqref{mapPar}$ above pulls back $\mathcal{B}(\vec{\lambda},\ell)$ to $\mathcal{O}_{\Bbb{P}^1}(\ell-\lambda_{1}(\theta^{\vee}))\boxtimes\boxtimes _{i=2}^{s} L(\lambda_i)$.
\item The map $\eqref{mapPar}$ can be extended to an atlas of  $\Parbun_{G,S}$. Therefore any codimension two substack of $\Parbun_{G,S}$ pulls back to a codimension two subscheme of $\Bbb{P}^1\times (G/B)^{s-1}$.
\end{enumerate}
\end{proposition}

For part (1) we are  reduced to showing that the map $\Bbb{P}^1\to \Parbun^0_G$ pulls back $L(\lambda)\tensor \mathcal{L}^{\ell}$   to  $\mathcal{O}(\ell-\lambda (\theta^{\vee}))$ where $\mathcal{L}$ is the positive generator of the Picard group of $\Bun_{SL_2}$. For the factor $L(\lambda)$  we use the map $B_{SL_2}\to I_G$  and compose with the corresponding character of $I_G$. For the statement for the factor $\ml$, see  \cite[Section 10.2.2]{SorgerLectures}.

\subsection{Proof of  Proposition \ref{properties} (2)}
Note that $\eqref{mapPar}$ is a base change of $\Bbb{P}^1\to \Parbun^0_G$. It therefore suffices to show that $\Bbb{P}^1\to \Parbun^0_G$ is smooth. To analyse this consider the two coordinate patches of $\Bbb{P}^1$.

\begin{equation}\label{la_map}
p:\Bbb{A}^1\to SL_2(\Bbb{C})/B_{SL_2}, t\in \Bbb{A}^1\mapsto \begin{pmatrix}
1 & 0 \\
 t & 1
\end{pmatrix}\in SL_2(\Bbb{C}).
\end{equation}
Under \eqref{la_map}, $t=\infty$ goes to $\begin{pmatrix}
0 & -1 \\
 1 & 0
\end{pmatrix}$. Since we want a parameterization of \eqref{la_map} at $t=\infty$, we set $T=1/t$ and compute
\begin{lemma}\label{lem2}
$$p(T)=\begin{pmatrix}
1 & T \\
 0& 1
\end{pmatrix}\begin{pmatrix}
0 & -1 \\
 1 & 0
\end{pmatrix}B_{SL_2}\in \op{SL}_2(\Bbb{C})/B_{SL_2}$$

\end{lemma}

We now prove that $\Bbb{P}^1\to \Parbun^0_G$ is smooth. Any point $T\neq 0$ maps to a parabolic bundle with trivial underlying bundle. The map $\Bbb{P}^1\to \Parbun^0$ is smooth at such points: Note that any map $S\to [pt/A]$ where $A$ is an algebraic group, and $S$ is a smooth variety, is smooth, because such a map consists of the data of a principal $A$ space, $E\to S$. When we base change this map to the atlas $pt \to pt/A$, we get $E\to pt$. Now  $E$ is a smooth variety, because it is smooth over a smooth variety $S$. For the problem at hand $A$ is the subgroup of the  automorphism group of the trivial $G$ bundle on $\Bbb{P}^1$, formed by automorphisms which preserve the standard flag at $0$.

Consider the point $T=0$, i.e., $t=\infty$. The local deformation space of a principal $G$ bundle $E$ on $\Bbb{P}^1$ is $H^1(\Bbb{P}^1,E\times_G\frg)$ \cite{Ram}. If we fix the choice of a section $s\in E_0$ (i.e., a trivialization of the fiber over $0$), the deformation space is $H^1(\Bbb{P}^1,E\times_G\frg(-1))$. At infinity $E$ comes from the patching data $\theta:\Bbb{C}^*\to G$, and the $H^1$ above is a direct sum over root spaces, and the Cartan. Note that the line bundles look like $\mathcal{O}(\alpha(H_{\theta})-1))$ where $\alpha$ runs through all roots. If we choose a $B$-structure instead of a trivialization of the fiber, the direct sum looks like $H^1(\Bbb{P}^1,\mathcal{O}(\alpha(H_{\theta})-1)))$ when $\alpha$ is a positive root, and $H^1(\Bbb{P}^1,\mathcal{O}(\alpha(H_{\theta})))$ when $\alpha$ is a negative root. There is only one contribution of the latter type, which is for $\alpha=-\theta$, and is one dimensional. The transition functions of our family definitely hit this.

\begin{lemma}
Let $\overline{g(z)}\in \mathcal{I}_G=G(\Bbb{C}((z)))/I_G$. Let $E$ be its image in  $\Bun_G$. Then the automorphism group of $E$ (without parabolic structures) is the set of all $A(z)\in G(\Bbb{C}[[z]])$ such that $g(z)A(z)g(z)^{-1}$ is holomorphic on the set $z\neq 0$, i.e, on $\Bbb{P}^1-\{0\}$.

If $g(z)=s_{\theta}z^{\theta^{\vee}}$, then the automorphism group acting on the space $E_0/B$ has the following property:  The given element in $E_0/B$ has a dense orbit. This is because the automorphism group then contains $B^-\subseteq G[[z]]$, since $\alpha(\theta^{\vee})\geq 0$ for positive roots $\alpha$.

\end{lemma}

\begin{proof}
Recall that $B^-$ is generated by negative root subgroups $U_{-\gamma}$ and the maximal torus $H$.
Let $a\in \C$ and consider the element $U_{-\gamma}(a)$. With $g(z) = s_\theta z^{\theta^\vee}$, we have
$$
g(z)U_{-\gamma}(a)g(z)^{-1} = s_\theta U_{-\gamma}(z^{-\langle \gamma, \theta^\vee\rangle}a)s_\theta^{-1},
$$
which is clearly regular at $z = \infty$ given that $\langle \gamma,\theta^\vee\rangle \ge 0$.

For $h\in H$, we know that $z^{\theta^\vee}$ commutes with $h$, so
$$
g(z) h g(z)^{-1} = s_\theta h s_\theta^{-1},
$$
which is also regular at $z = \infty$.
\end{proof}

\section{A basic computation}
\

In Section \ref{Stepp1}, we completed the first step in the proof of Proposition \ref{generale} as outlined in Section \ref{outline}.
Assume $\ell(v'_{1})\geq \ell (v_{1})$.  The second  step is completed in Section \ref{secondstep}.

\subsubsection{ A Basic Computation}\label{notate}
We have a family of $G$-bundles $E(T)$ with $T\in \Bbb{P}^1$. These bundles are trivialized outside of zero. Fix $(\bar{g}_2,\dots,\bar{g}_s)\in (G/B)^{s-1}$.  The map \eqref{mapPar} gives a parabolic structure at $(p_1=0,p_2,\dots,p_s)$ on $E(T)$ for all $T\in \Bbb{P}^1$.

Suppose we have a map $f:\Bbb{P}^1\to E(T)/P$ for $T\in \Bbb{A}^1$. Since $E(T)$ are trivialized outside of $0$, $f$ gives a map $g:\Bbb{P}^1\to E(T=\infty)/P$. We can then ask  to compare the Schubert state of $f$ and $g$ at $p_1=0$ (at other points they are obviously the same). For convenience we will assume that the parabolic bundle at $T=\infty$ (with trivial underlying bundle) is general.

Let $U_{-\theta}(Tz)\in G((z))$ be the image of the matrix $\begin{pmatrix}
1 & T \\
 0 & 1
\end{pmatrix}$
under the composite $$SL_2\to SL_2(\Bbb{C}((z)))\to G((z)).$$ By construction $E(T)$ is trivialized near $0$ and also outside of $0$. The section ``1'' defined near zero goes to the section  ``$U_{-\theta}(Tz)s_\theta^{-1}z^{\theta^\vee}$'' outside of zero. Consider a $P$ reduction of $E_T$ defined near $z=0$ by a matrix $A(z)P$. Here $A(z)$ is
holomorphic at $z=0$. The corresponding $P$ reduction of $E(0)$ in the coordinate system near $z=0$ is given by $U_{-\theta}(Tz)s_\theta^{-1}z^{\theta^\vee} A(z)$. Write $A(z)= I(z) bw$, where $I(z)$ is in the kernel of $G[[z]]\to G$, $b\in B$, and $w\in W$. We therefore want to understand the limit
\begin{equation}\label{limitt}
\lim_{z\to 0}U_{-\theta}(Tz)s_\theta^{-1}z^{\theta^\vee} I(z)bwP
\end{equation}
\begin{remark}
We have therefore assumed the form of the $P$-reduction of $E(T)$ at $z=0$, and use it to deduce information about the $P$-reduction of
$E(t=0)(=E(T=\infty))$ at $0$.
\end{remark}


The limit \eqref{limitt}
exists in $G/P$; therefore we have an equation
$$
U_{-\theta}(Tz) s_{\theta}^{-1}z^{\theta^\vee} I(z)bw = g(z) p(z)
$$
for some $g(z)\in G(\mathcal{O}), p(z)\in P(\mathcal{K})$.

\begin{remark}\label{degreecount}
Let $$d_T=(a_i)_{\alpha_i\in S_P}=\sum_{\alpha_i \in S_P} a_i\mu(X^P_{s_i})$$ be the degree of the $P$-reduction of $E(T)$ given by $f:\Bbb{P}^1\to E(T)/P$, and
$$d_{t=0}=(b_i)_{\alpha_i\in S_P}=\sum_{\alpha_i \in S_P} b_i\mu(X^P_{s_i})$$ be the degree of the  corresponding $P$ reduction of
$E(t=0)$ given by $g:\Bbb{P}^1\to E(t=0)/P$. Recall that $a_i =\deg f^* E\times_P \Bbb{C}_{-\omega_i}$ and $b_i =\deg g^* E\times_P \Bbb{C}_{-\omega_i}$.
Note that $f^* E\times_P \Bbb{C}_{-\omega_i}$ and $ g^* E\times_P \Bbb{C}_{-\omega_i}$ are identified outside of $0$ on $\Bbb{P}^1$. It follows that
$a_i-b_i$ equals the order of vanishing of $\omega_{i}(p(z))$ at $z=0$.

\end{remark}

We wish to find
\begin{enumerate}
\item
The order of vanishing of $\omega_i(p(z))$ for $\alpha_i\in S_P$, as a function of the initial data $b,w,I,T$.
\item The Schubert cell in which $g(0)\in G/P$ lies, with respect to the standard flag.
\end{enumerate}

Immediately we observe that the leftmost factor of $U_{-\theta}(Tz)$ has no effect on the answers to (1) and (2). In general, in answering these questions, we can and will left-multiply $g(z)$ by some element of the Iwahori subgroup $I_G$ with no change to (1) or (2).

\begin{proposition}\label{muncie}
Suppose $I(z)$ is in the kernel of  the $G(\Bbb{C}[[z]])\to G$ given by evaluation at $z=0$. Write a factorization of $I(z)$ into root subgroups in a particular order (the order within the parentheses does not affect our approach, so we leave it unspecified):
$$
I(z) = \left(\prod_{\gamma\ne -\theta }U_\gamma(a_\gamma)\right) U_{-\theta}(a_{-\theta}),
$$
where each $a_\gamma$ (including $a_{-\theta}$) is a power series in $z$ divisible by $z$.

Suppose
$$
s_{\theta}^{-1}z^{\theta^\vee} I(z)bw = g(z) p(z)
$$
where $g(z)\in G(\mathcal{O}), p(z)\in P(\mathcal{K})$.

Then, if $w^{-1}\theta^\vee\prec 0$ or $\op{val}(a_{-\theta}) \ge 2$, then $g(0)\in Bs_\theta wP$ and for $i \in S_P$, $\op{ord}(\omega_i(p(z))) = \omega_i(w^{-1}\theta^\vee)$.
Otherwise, $g(0)\in BwP$ and $\op{ord}(\omega_i(p(z))) = 0$.
\end{proposition}
\subsection{Proof of Proposition \ref{generale}, Step 2}\label{secondstep}
Fix a general point of $(G/B)^{s-1}$. The map $\pi:\Bbb{P}^1\times (G/B)^{s-1} \to   \Parbun_{G,S}$ gives a family of $G$-bundles $E(T)$ on $\Bbb{P}^1$ with $T\in \Bbb{P}^1$. These bundles are trivialized outside of zero. Fix $(\bar{g}_2,\dots,\bar{g}_s)\in (G/B)^{s-1}$.  We also get parabolic structures at $(p_1=0,p_2,\dots,p_s)$ on $E(T)$ for all $T\in \Bbb{P}^1$. We may assume that the bundle at $T=\infty$ (i.e., $t=0$) does not lie in $\Omega(\vec{v},d)$. 
We now begin the proof of Step 2. Suppose to the contrary that $E(T)$ is in $\Omega^0(\vec{v},d)$. We have a map
$f:\Bbb{P}^1\to E(T)/P$ which gives rise to a map $g:\Bbb{P}^1\to E(T=\infty)/P$ as before.  

We claim that $g$ gives a point of $\Omega^0(\vec{v},d)$ a contradiction. 
To do this we divide into two cases:
\begin{enumerate}
\item  $v_1^{-1}\theta\in R_{\frl}$. In this case $v'_1=v_1$ by Lemma \ref{fraction}, and $\omega_i(v_1^{-1}\theta)=0$, and so the degrees of the maps $g$ and $f$ are same. Therefore $g$ gives a point of $\Omega^0(\vec{v},d)$.
\item $v_1^{-1}\theta\in R^{+}-R^+_{\frl}$. If $\op{val}(a_{-\theta}) \ge 2$ in Proposition \ref{muncie}, then degree of the $P$ reduction corresponding to $g$ equal to $d-d(\alpha)$, where $\alpha$ is the positive root $v^{-1}\theta$. The expected dimension of such $g$ is one less than the expected dimension for $f$, hence negative (in fact, $g$ after attaching a $T$ fixed curve lies in a quantum  stratum (Defn. \ref{Maxim}) of $\Omega(\vec{v},d)$ and therefore contributes to $\Omega(\vec{v},d)$). In the remaining case of Proposition \ref{muncie}, $g$ is in $\Omega^0(\vec{v},d)$.
\end{enumerate}

\subsection{Proof of Proposition \ref{muncie}}
To help break up the proof, we first record the following useful lemmas.
\begin{lemma}\label{sanders}
Suppose $\gamma$ is a root not equal to $-\theta$, and let $a\in \C[[z]]$. If either
\begin{enumerate}
\item $\langle \gamma,\theta^\vee\rangle \ge 1$, or
\item $\op{val}(a)\ge 1$,
\end{enumerate}
then $s_\theta^{-1}z^{\theta^\vee} U_\gamma(a) z^{-\theta^\vee} s_\theta$ belongs to $I_G$.
\end{lemma}

\begin{proof}
We observe that
$$s_\theta^{-1}z^{\theta^\vee} U_\gamma(a) z^{-\theta^\vee} s_\theta = U_{s_\theta \gamma} \left(z^{\langle \gamma, \theta^\vee\rangle} a\right).$$
In case (1), this expression goes to the identity as $z\to 0$.

For case (2), we note that $\langle \gamma,\theta^\vee\rangle \ge -1$ since $\gamma\ne -\theta$. Moreover, we may as well assume $\langle \gamma, \theta^\vee\rangle \in \{-1,0\}$ since otherwise (1) applies. If $\langle \gamma,\theta^\vee\rangle = 0$, then since $a\to 0$ as $z\to 0$,
$$
U_{s_\theta \gamma} \left(z^{\langle \gamma, \theta^\vee\rangle} a\right) = U_{\gamma}(a)
$$
once again goes to the identity as $z \to 0$.

Otherwise, $\langle \gamma,\theta^\vee\rangle  = -1$ and $\gamma+\theta$ is a root, necessarily positive. So
$$
U_{s_\theta \gamma} \left(z^{\langle \gamma, \theta^\vee\rangle} a\right) = U_{\gamma+\theta}(z^{-1}a),
$$
which belongs to $B(\C[[z]])\subset I_G$.
\end{proof}

It is also helpful to recall the formula for the commutator of two root subgroups:
\begin{lemma}{\cite{Stein}*{Chapter 3}}\label{sberg}
For arbitrary roots $\gamma_1, \gamma_2$ such that $\gamma_1+\gamma_2 \neq 0$, we have
$$
U_{\gamma_1}(a)U_{\gamma_2}(b)U_{\gamma_1}(a)^{-1} U_{\gamma_2}(b)^{-1}=\prod_{i \gamma_1+j \gamma_2 \in \Phi} U_{i\gamma_1+j\gamma_2}(c_{\gamma_1,\gamma_2}^{i,j}a^ib^j)
$$
for $i, j \in \mathbb{Z}^+$ and integer constants $c_{\gamma_1,\gamma_2}^{i,j}$ which depend on $i,j,\gamma_1,$ and $\gamma_2$.
\end{lemma}

When we wish to use this commutator formula, we will be in the specific situation where $\gamma_1 = \theta$ and $\langle \gamma_2,\theta^\vee\rangle = -1$. In that context, there are the following constraints on $i$ and $j$. A short proof can be found at the end of this section. 

\begin{lemma}\label{greater-or-one}
Suppose $\gamma$ is any root such that $\langle \gamma, \theta^\vee\rangle = -1$. Then for all $i,j\in \Z^+$ such that
$$
i\theta + j\gamma
$$
is a root, we have $i\le j$. Moreover, if $i=j$ then $i=j=1$.
\end{lemma}


\begin{proof}[Proof of Proposition \ref{muncie}]
As already noted, we can ignore the factor of $U_{-\theta}(Tz)$ and focus solely on
$$
s_\theta^{-1} z^{\theta^\vee} I(z) bw = g(z) p(z).
$$

By Lemma \ref{sanders}(2), we can reduce to considering just
$$
s_\theta^{-1} z^{\theta^\vee} U_{-\theta}( a_{-\theta}) b w =  U_{\theta}(z^{-2} a_{-\theta}) s_\theta^{-1} z^{\theta^\vee}b w.
$$

Write
$$
b = \left(\prod_{\tiny \begin{array}{c} \alpha\succ 0\\ \langle \alpha,\theta^\vee\rangle = 0\end{array}} U_\alpha(b_\alpha)\right)
\left(\prod_{\tiny \begin{array}{c} \alpha\succ 0\\ \langle \alpha,\theta^\vee\rangle = 1\end{array}} U_\alpha(b_\alpha)\right) U_\theta(b_\theta).
$$

Then we move the $s_\theta^{-1}z^{\theta^\vee}$ past these roots subgroups as follows:
\begin{align}\label{yikes}
s_\theta^{-1} z^{\theta^\vee} b = \left(\prod_{\tiny \begin{array}{c} \alpha\succ 0\\ \langle \alpha,\theta^\vee\rangle = 0\end{array}} U_\alpha(b_\alpha)\right)
\left(\prod_{\tiny \begin{array}{c} \alpha\succ 0\\ \langle \alpha,\theta^\vee\rangle = 1\end{array}} U_{\alpha-\theta}(zb_\alpha)\right) U_{-\theta}(z^2b_\theta)s_\theta^{-1} z^{\theta^\vee}.
\end{align}

Next we want to move the $U_{\theta}(z^{-2}a_{-\theta})$ past these root subgroups, picking up some extra terms along the way as in Lemma \ref{sberg}. First we examine the root subgroups $U_\alpha(b_\alpha)$ such that $\langle \alpha,\theta^\vee\rangle = 0$:
\begin{align}\label{commuters}
U_{\theta}(z^{-2}a_{-\theta}) U_\alpha(b_\alpha) = U_\alpha(b_\alpha) U_{\theta}(z^{-2}a_{-\theta});
\end{align}
for future use set $X_\alpha = U_\alpha(b_\alpha)$ when $\langle \alpha,\theta^\vee\rangle = 0$; note that $X_\alpha\in I_G$ since $X_\alpha\in B$.

The more interesting case is when $\langle \alpha,\theta^\vee \rangle = 1$ and we want to move the $U_{\theta}$ past $U_{\alpha-\theta}$. Set $\eta_{i,j} = i\theta+j(\alpha-\theta)$. By Lemmas \ref{sberg} and \ref{greater-or-one} (here $\gamma = \alpha-\theta$),
\begin{align}\label{stickers}
U_\theta(z^{-2}a_{-\theta})&U_{\alpha-\theta}(zb_\alpha) =  U_{\alpha}\left(c_{\theta,\alpha-\theta}^{1,1}z^{-1}a_{-\theta}b_\alpha\right)\times \\\nonumber
&\left(\prod_
{\tiny
\begin{array}{c}
i<j \\
\eta_{i,j}\in \Phi
\end{array}
}
U_{\eta_{i,j}}\left(c_{\theta,\alpha-\theta}^{i,j}(z^{-2}a_{-\theta})^i(zb_\alpha)^j\right)
\right)U_{\alpha-\theta}(zb_\alpha)U_\theta(z^{-2}a_{-\theta}).
\end{align}
Therefore set
$$
X_\alpha:= \left(U_{\alpha}\left(c_{\theta,\alpha-\theta}^{1,1}z^{-1}a_{-\theta}b_\alpha\right)\prod_
{\tiny
\begin{array}{c}
i<j \\
\eta_{i,j}\in \Phi
\end{array}
}
U_{\eta_{i,j}}\left(c_{\theta,\alpha-\theta}^{i,j}(z^{-2}a_{-\theta})^i(zb_\alpha)^j\right)
\right)U_{\alpha-\theta}(zb_\alpha)
$$
for those roots such that $\langle \alpha,\theta^\vee\rangle = 1$. Note that $X_\alpha\in I_G$, since $\alpha\succ 0$ and since when $i<j$, the valuation of
$$
(z^{-2}a_{-\theta})^i(zb_\alpha)^j
$$
is at least $j-i$.

To summarize, by combining (\ref{yikes}), (\ref{commuters}), and (\ref{stickers}), we have managed to write
$$
U_\theta(z^{-2}a_{-\theta}) s_\theta^{-1} z^{\theta^\vee}b = \left(\prod_{\alpha\succ 0} X_\alpha\right) U_\theta(z^{-2}a_{-\theta})U_{-\theta}(z^2b_\theta)s_\theta^{-1} z^{\theta^\vee}
$$
Since everything in the parentheses belongs to $I_G$, we have reduced the problem to studying
$$
\bar g(z) p(z) = U_\theta(z^{-2}a_{-\theta})U_{-\theta}(z^2b_\theta)s_\theta^{-1} z^{\theta^\vee}w.
$$
Here $\bar g(z)\in G(\mathcal{O})$ and $p(z)\in P(\mathcal{K})$, and we want to calculate the Schubert cell containing $\bar g(0)P$ and the order of vanishing $\op{ord}(\omega_i(p(z)))$ when $\alpha_i\in S_P$.

From the $SL_2$ theory, we find that
$$
U_\theta(z^{-2} a_{-\theta}) U_{-\theta}(z^2 b_\theta) = U_{-\theta}(u^{-1}z^2 b_\theta) u^{\theta^\vee} U_\theta(z^{-2}a_{-\theta} u^{-1}),
$$
where $u\in \mathcal{O}$ is the unit $u = 1+a_{-\theta}b_\theta$.

Up to left multiplication by $I_G$, we just examine
\begin{align}\label{remainder}
U_\theta(z^{-2}a_{-\theta} u^{-1})s_\theta^{-1} z^{\theta^\vee}w  = U_\theta(z^{-2}a_{-\theta}u^{-1}) s_\theta^{-1} w z^{w^{-1}\theta^\vee}.
\end{align}
If $a_{-\theta}$ is divisible by $z^2$, then this gives $g(0) \in Bs_\theta wP$ and $\op{ord}(\omega_i(p(z))) = \omega_i(w^{-1}\theta^\vee)$ for $\alpha_i\in S_P$.

If $-w^{-1}\theta\succ 0$, then (\ref{remainder}) can be rewritten as
$$
s_\theta^{-1} w U_{-w^{-1}\theta} (z^{-2} a_{-\theta} u^{-1})z^{w^{-1}\theta^\vee},
$$
which once again yields $g(0)\in B s_\theta w P$ and $\op{ord}(\omega_i(p(z)) = \omega_i(w^{-1}\theta^\vee)$ for $\alpha_i\in S_P$.

Otherwise, for the remainder of the proof, we assume:
\begin{itemize}
\item $\op{ord}(a_{-\theta}) = 1$ and
\item $w^{-1}\theta\succ 0$.
\end{itemize}

Thus $f :=a_{-\theta}/z$ is a unit in $\mathcal{O}$ and $z^{-2}a_{-\theta} = z^{-1}f$.
Once again from the $SL_2$ theory, we have an identity
$$
U_\theta(h^{-1}) = U_{-\theta}(h) s_\theta^{-1}U_{-\theta}(h^{-1}) h^{\theta^\vee}
$$
with $h = uf^{-1}z$. Substituting for $U_\theta(h^{-1})$ in (\ref{remainder}), and moving the $s_\theta w$ term to the left, we get
$$
U_{-\theta}(uf^{-1}z) s_\theta^{-1} U_{-\theta}(uf^{-1}z) h^{\theta^\vee}s_\theta w z^{w^{-1}\theta^\vee} = U_{-\theta}(uf^{-1}z)
wU_{w^{-1}\theta}(u f^{-1} z)h^{-w^{-1}\theta^\vee}z^{w^{-1}\theta^\vee}.
$$
Set $p(z) = U_{w^{-1}\theta}(u f^{-1} z)h^{-w^{-1}\theta^\vee}z^{w^{-1}\theta^\vee}$ and $g(z) = U_{-\theta}(uf^{-1}z)w$. Applying $\omega_P$ to $p(z)$ yields the same result as $\omega_P(h^{-w^{-1}\theta^\vee}z^{w^{-1}\theta^\vee})$, which is just $1$ since
$$
h^{-w^{-1}\theta^\vee}z^{w^{-1}\theta^\vee} = (f^{-1} u)^{-w^{-1}\theta^\vee} z^{-w^{-1}\theta^\vee} z^{w^{-1}\theta^\vee} \in T(\mathcal{O}).
$$
So we conclude that $g(0)\in BwP$ and $\op{ord}(\omega_i(p(z))) = 0$ for $\alpha_i\in S_P$.
\end{proof}

\subsection{Conclusion of proof of Proposition \ref{generale}}\label{Step3proof} 
We  complete Step 3 of the outline. Assume $\ell(v'_{1})< \ell (v_{1})$, and we show then that the fiber of $\pi^{-1}(\mathcal{D}_{\vec v})\to (G/B)^{s-1}$ over a general point  $(\bar{g_2},\dots,\bar{g}_s)$ of
$(G/B)^{s-1}$ is in bijection with the enumerative count  $\langle v'_1,v_2,\hdots, v_s\rangle_{d+m}$.

We assume that the parabolic bundle at $T=\infty$ (with trivial underlying bundle) is general.
We use the notation of Proposition \ref{muncie}. If the Schubert position of $f(0)$ is $\tilde{w}P$, then the Schubert position of $g(0)$ is $s_{\theta}\tilde{w}P$ or $\tilde{w}P$ (with corresponding changes in degrees). 

Therefore if $w=v_1$ and the degree is $d$ then (using Lemma \ref{fraction} and Remark \ref{degreecount}), the Schubert position  of $g$ is $v_1'$ since $v_1^{-1} \theta \prec 0$ by Lemma \ref{fraction}. The degree of  $g$ is $d+d(-v_1^{-1}\theta)$. This implies that 
$g$ contributes to the count $\langle v'_1,v_2,\hdots, v_s\rangle_{d+m}$ for the trivial bundle with flags $(\dot{e},\bar{g_2},\dots,\bar{g}_s)$.

Finally  fixing the point of $(G/B)^s$ we need to show that there is exactly one value of $T$ which via $f$ gives rise to a given $g$. This means that the coefficient $a_{-\theta}$ varies linearly in $T$ and vanishes to order $\geq 2$ exactly once:

We get reductions of $E(T)/P$ for all values of $T\neq 0$, and also to $t=0$. Fix a value of $T$.
Suppose the  $P$ reduction of $E(T)$ is defined near $z=0$ by a matrix $A(z)P$. Here $A(z)$ is
holomorphic at $z=0$. Write  $A(z)= I(z) bw$.

The corresponding $P$ reduction of $E(0)$  near $z=0$ is given by $U_{-\theta}(Tz)s_\theta^{-1}z^{\theta^\vee} A(z)$.
Suppose this equals $I_0(z)b v'_1P$ (since we want $g(0)$ to correspond to the enumerative count $\langle v'_1,v_2,\hdots, v_s\rangle_{d+m}$).

Therefore we have an equation
  $$U_{-\theta}(Tz) s_{\theta}^{-1}z^{\theta^\vee} I(z)bwP= I_0(z)bv'_1P.$$

 and hence
 $$ I(z)wP= z^{-\theta^{\vee}} s_{\theta} U_{-\theta}(-Tz) I_0(z)bv'_1P.$$
 Note that $z^{-\theta^{\vee}} s_{\theta}=s_{\theta}^{-1}z^{\theta^{\vee}}$, and thus
 $$ I(z)bwP= s_{\theta}^{-1} z^{\theta^{\vee}}  U_{-\theta}(-Tz) I_0(z)bv'_1P.$$

 We break up $I_0(z)$ into a product
 $$
I_0(z) = U_{-\theta}(a_{-\theta})\cdot \left(\prod_{\gamma\ne -\theta }U_\gamma(a_\gamma)\right) ,
$$
where each $a_\gamma$ is a power series in $z$ divisible by $z$ (the order in the second product is unspecified). 
And hence
$$ I(z)b'wP= s_{\theta}^{-1} z^{\theta^{\vee}}   U_{-\theta}(a_{-\theta}-Tz)\cdot \left(\prod_{\gamma\ne -\theta }U_\gamma(a_\gamma)\right) bv'_1P$$

By the main computation in the proof of Proposition  \ref{muncie} (and using Lemma \ref{sberg} to move the $U_{-\theta}$ term) , if $\op{val}(a_{-\theta}-Tz)>1$, then $w=s_{\theta}v'_1=v_1$ (note that $(v'_1)^{-1}\theta\succ 0$) and the
degree rise is as expected for $f$ to have degree $d$). There is a unique value of $T$ when this happens. For other values of $T$, $w=v'_1$ with no drop in degree. The value of $T$ could be zero (for example when $d<0$).

\begin{proof}[Proof of Lemma \ref{greater-or-one}]
Given that $\langle \gamma,\theta^\vee\rangle = -1$ and $(\theta,\theta) = 2$, we have $(\gamma,\theta) = -1$.
If $\eta = i\theta+j\gamma$ is a root, then $(\eta,\eta)\le 2$. Consider
\begin{align*}
(i\theta+j \gamma, i\theta+j\gamma) &= 2i^2-2ij+j^2(\gamma,\gamma) \\
&\ge 2i^2-2ij+2/3 j^2\\
&\ge 2i^2-2ij+1/2 j^2 \\
&= 1/2(2i-j)^2,
\end{align*}
so $4\ge (2i-j)^2$ and hence $2i-j\le 2$.

If $j\ge 2$, then we obtain $2i\le 2+j\le 2j$, so that $i\le j$ as desired.

Otherwise, $j=1$. Since $\eta = i\theta+\gamma$ is a root, we must have $(\eta,\theta)\in \{-2,-1,0,1,2\}$. But
$$
(\eta,\theta) = 2i-1,
$$
which is odd and positive. Hence $(\eta,\theta) = 1$ and $i=1$; $\eta = \theta+\gamma$.

Finally, if $i=j\ge 0$ and $\eta = i(\theta+\gamma)$ is a root, then $i=1$ is forced, since $\theta+\gamma$ is also a root.
\end{proof}

\section{A basic diagram, Levification and the ramification divisor}

Recall that  (Defn. \ref{stacky}) $\Omega^0(\vec{u},d)$ parametrises principal bundles $\me$ together with elements $\me_{p_i}/B$, and a map $F:\Bbb{P}^1\to\me/P$ of degree $d$ so that $f(p_i)$ are in relative position $u_i$ with the flags $\me_{p_i}/B$. The  principal bundle and the map $F:\Bbb{P}^1\to \me/P$ are equivalent to the data of a $P$-bundle. The parabolic data at the points is captured by some additional data as in \cite[Lemma 3.3]{BKiers}.
\begin{lemma}\label{cantata}
The stack $\Omega^0(\vec{u},d)$ parametrises principal $P$-bundles $\mathcal{P}$ together with elements $\overline{z}_i\in\mathcal{P}/(u_i^{-1}Bu_i\cap P)$ such that
$\deg(\mathcal{P}\times_P \Bbb{C}_{-\omega_i})=a_i$ for each $\alpha_i\in S_P$.
\end{lemma}
\begin{defi}
 $\Parbun_L(d)$ parametrises data of an $L$-bundle $\mathcal{E}'$ and $\alpha_i\in \me'_{p_i}/B_L$, such that $\deg(\mathcal{E'}\times_L \Bbb{C}_{-\omega_i})=a_i$ for each $\alpha_i\in S_P$.
\end{defi}

We want to define maps to obtain the following basic diagram
\begin{equation}\label{Dbasic}
\xymatrix{
 & \Omega^0(\vec{u},d)\ar[dl]^{\pi}\ar[dr]_{\xi}\\
\Parbun_G   &  & \Parbun_L(d)\ar@/_/[ul]_i\ar[ll]^{i'=\pi\circ i}
}
\end{equation}
\begin{remark}
As in the related contexts of \cite{BRays,BKiers,KBranch, BRigid}, the use of this diagram to pull back line bundles from $\Parbun_L(d)$ and pushing them forward (suitably) to $\Parbun_G$ is inspired by \cite[Section 4.1]{ress}.
\end{remark}

\begin{itemize}
\item
The map $\pi$ has already been defined.
\item The map $\xi$ arises from Lemma \ref{cantata}, the morphism $P\to L$ which takes $u_i^{-1}Bu_i\cap P$ onto $B_L$ (cf. \cite[Lemma 3.2]{BKiers}).
\item The map $i$ takes a pair $\mathcal{E}'$  of a $L$ bundle of degree $d$ and data $\alpha_i\in \me'_{p_i}/B_L$ to the $P$ bundle $\mathcal{P}=\mathcal{E}'\times_L P$ and $\alpha_i$ give elements in $\mathcal{P}/(u_i^{-1}Bu_i\cap P)$ since $u_iB_Lu_i^{-1}\subseteq B$.
\item The map $i'$ has the following description: It takes a pair $\mathcal{E}'$  of a $L$ bundle of degree $d$ and data $\alpha_i\in \me'_{p_i}/B_L$ to the $G$-bundle
$\mathcal{P}=\mathcal{E}'\times_L G $ and data $\alpha_i\times u_i^{-1}$.
\end{itemize}

\subsection{Levification}\label{secLev}

The space $\Omega^0(\vec{u},d)$ retracts to $\Parbun_L(d)$ by a process called Levification (cf. \cite[Section 3.8]{BKq} also \cite[Section 3.6]{BKiers}). Consider an element $x_L=\sum'_{k} N_kx_k$ where the sum is over $k$ such that $\alpha_k\not\in \Delta(P)$, with $N_k$ such that $N_kx_k$ is in the co-root lattice. Then $t^{x_L}=\exp((\ln t) x_L)$ topologically generates $Z^0(L)$ if $P$ is maximal. Note that $P$ is arbitrary in this section.

For $t\in\Bbb{C}^*$  consider (cf. \cite[Section 3.8]{BKq}) $\phi_t(p)=t^{x_L}pt^{-x_L}$, with  $\phi_1$ the identity on $P$.
This extends to a group homomorphism $\phi_0:P\to L$ which coincides with the standard projection of $P$ to $L$ giving rise to a morphism $\hat{\phi}:P\times \Bbb{A}^1\to P$. Clearly, $\phi_t:L\to L$ is the identity on $L$ for all $t$.

\begin{defi}\label{levity}
Let $(\mathcal{P},\bar{z}_1,\dots,\bar{z}_s)$ be a point of $\Omega^0(\vec{u},d)$ where $\mathcal{P}$ is a principal $P$-bundle and $\bar{z}_i\in
\mathcal{P}_{p_i}/(w_i^{-1}Bw_i)\cap P$. Define the Levification family $\mathcal{P}_t=\mathcal{P}'\times_{\phi_t} P$ for $t\in \Bbb{A}^1$, and $\bar{z}_i(t)=\bar{z}_i\times_ {\phi_t} e$. Clearly, at $t=0$, $(\mathcal{P}_t,\bar{z}_1(t),\dots,\bar{z}_s(t))$ is in the image of $i:\Parbun_L(d)\to \Omega^0(\vec{u},d)$, and equals $i\circ\xi(\mathcal{P},\bar{z}_1,\dots,\bar{z}_s)$.
\end{defi}

\subsection{Sections of line bundles}
We recall some definitions generalizing \cite[Section 3.7]{BKiers}
\begin{defi}\label{indexx}
Let $\mathcal{M}$ be a line bundle on $\Parbun_L(d)$.  $Z(L)$ acts on fibers  of $\mathcal{M}$ and gives rise to a (multiplicative) character ``index" $\gamma_{\mathcal{M}}:Z^0(L)\to \Bbb{C}^*$.  More generally, this character can be defined if $\mathcal{M}$ is defined over an open substack  of $\Parbun_L(d)$ since $\gamma_{\mathcal{M}}$ is constant over connected families.
\end{defi}
\begin{proposition}\label{comparison}
Let $U$ be a non-empty open substack of $\Parbun_L(d)$,  $\ml$ be a line bundle on $\xi^{-1}(U)$ and $\mathcal{M}=i^*\mathcal{L}$, a line bundle on $U$, where $i:\Parbun_L(d)\to \Omega^0(\vec{u},d)$ and $\xi:\Omega^0(\vec{u},d)\to \Parbun_L(d)$.
Then
\begin{enumerate}
\item $\ml=\xi^*\mathcal{M}$. Therefore $\xi^*$ and $i^*$ set up isomorphisms $\Pic(U)\leto{\sim}\Pic(\tau^{-1}(U))$.
\item If $\gamma_{\mathcal{M}}$ is trivial then
$H^0(\xi^{-1}(U),\ml)\to H^0(U,\mathcal{M})$ is an isomorphism.
\end{enumerate}
\end{proposition}
\begin{proof}
The second part is essentially  \cite[Theorem 15 and Remark 31(a)]{BKq}, and the main point is that if $\mathcal{P}_t$ is a Levification family then a section of $\ml$ (under the assumption of (2)) at $\mathcal{P}_1$ can be propagated in a unique way to all $\mathcal{P}_t$, $t\neq 0$ (since $\mathcal{P}_1$ is isomorphic to $\mathcal{P}_t$ for $t\neq 0$), and there are no poles or zeroes of this extended section at $t=0$.  For the surjectivity we can extend any section of $\ml$ at $\mathcal{P}_0$ to all  $\mathcal{P}_t$ since the corresponding $\Bbb{C}^*$-equivariant line bundle on $\Bbb{A}^1$ is trivial.

For the first part consider $\ml'=\ml\tensor \xi^*{\mathcal{M}}^{-1}$. Note that $\mathcal{M}'= i^* \mathcal{L}'$ is trivial, and
$\gamma_{\mathcal{M}'}$ is trivial. We can apply (2) to $(\ml',\mathcal{M}')$. The nowhere vanishing global section of $H^0(U,\mathcal{M}')$ gives a global section of $H^0(\xi^{-1}(U),\ml)$. It can be seen that this is nowhere vanishing as  well (see \cite[Lemma 3.17]{BKq}: Consider the corresponding Levification family (Definition \ref{levity}), if a global section vanishes for $\mathcal{P}_t$ then it will also vanish for $\mathcal{P}_0$.)
\end{proof}
\subsection{Ramification divisors}
We  recall from \cite[Section 3.7]{BKq} the definition of the ramification divisor $\mathcal{R}$ of the map $\Omega^0(\vec{u},d)\to \Parbun_G$ and show that  it is the pullback of a divisor $\mathcal{R}_L$ in $\Parbun_L(d)$ under the condition $\langle \sigma_{w_1}^P, \hdots, \sigma_{w_s}^P \rangle_d^{\circledast_0} =1$ (which implies Levi-movability).

Consider a point of $\Omega^0(\vec{u},d)$, a $P$-bundle $\mathcal{P}$ together with elements $\overline{z}_i\in\mathcal{P}/(u_i^{-1}Bu_i\cap P), 1\leq i\leq s$. Let
$T_{\dot{e}}=T(G/P)_{\dot{e}}$. Consider the vector bundle $\mathcal{K}$ on $\Bbb{P}^1$ given by the exact sequence
$$0\to\mathcal{K}\to \mathcal{P}\times^P T_{\dot{e}}\to \bigoplus_{i=1}^s (i_{p_i})_*\frac{\mathcal{P}\times^P T_{\dot{e}}}{T(\overline{z}_i,u_i,p_i)} \to 0$$
Here $T(\overline{z}_i,u_i,p_i)=z_i\times T(w_i^{-1}C_{u_i})_{\dot{e}}$.

By assumption on codimensions, $\mathcal{K}$ has zero Euler characteristic, and hence the determinant of cohomology $D(\mk)$ carries a canonical section $\theta$. This canonical section vanishes if and only if the point $\Omega^0(\vec{u},d)$ lies on the ramification divisor of $\pi: \Omega^0(\vec{u},d)\to\Parbun_G$ (since $\pi$ is a representable morphism of smooth stacks, it is meaningful to speak of its ramification divisor). The condition of Levi-movability is that the character $\gamma_{D(\mk)}$ is zero \cite[Proposition 3.14]{BKq}. This means, by Lemma \ref{comparison} that $D(\mathcal{K})$ is the pull back of a line bundle on $\Parbun_L$, and $\theta$ is the pullback of  a section of a corresponding section on $\Parbun_L$, Let $\mathcal{R}\subseteq \Omega^0(\vec{u},d)$ be the divisor of $\theta$; it is the pullback of a divisor $\mathcal{R}_L\subseteq\Parbun_L(d)$, which gives the following.

\begin{lemma}\label{imagemap}
The map $i':\Parbun_L(d)-\mathcal{R}_L\to \Omega^0(\vec{u},d)$ has image in  $\Omega^0(\vec{u},d)-\mathcal{R}_L$.
\end{lemma}

\section{Proof of Theorem \ref{type1}}\label{beginning}
The strategy of the proof of Theorem \ref{type1} is different from the proofs of analogous statements in \cite{BRigid,BRays,BKiers,KBranch}.
For one thing, the relevant generalization of Fulton's conjecture is not available yet (in fact it will follow from our work, see Corollary \ref{FulCon}).
We also do not try to extend the map $\Omega^0(\vec{u},d)-\mathcal{R}\to \Parbun_L(d)$ over a larger open subset of $\Parbun_G$ (see Lemma \ref{zussamen}). Such an
extension can probably be obtained by considering Drinfeld compactifications. But we instead argue in a more computational manner by first
considering a subset $\Pic'\subseteq \Pic(\Parbun_G)$  which we show (later, in Lemma \ref{zussamen2}) maps isomorphically to the Picard group of $\Omega^0(\vec{u},d)\setminus \mathcal{R}$, so that
$\Pic'$ and the basic divisors $D(v,j)$ generate  $\Pic(\Parbun_G)$.

Fix a parabolic $P$, Weyl group elements $u_1,\hdots,u_s\in W^P$, and degree $d$ such that $\langle \sigma_{u_1}^P,\hdots,\sigma_{u_s}^P\rangle_d^{\circledast_0}=1$. According to Theorem \ref{killing}, if $P$ is maximal, these give a facet of $\mathcal{C}$ defined by
$$
\mathcal{F}(P,u_1,\hdots,u_s,d) = \{(\lambda_1,\hdots,\lambda_s,\ell)\in \mathcal{C} : \sum_{i=1}^s \lambda_i(u_ix_P) = d\ell \frac{2}{(\alpha_P,\alpha_P)}\}.
$$

Suppose  $P$ is not necessarily maximal. The system of equalities for $\alpha_k\in S_P=\Delta\setminus \Delta_P$, and $(\lambda_1,\hdots,\lambda_s,\ell)\in \mathcal{C}$
\begin{equation}\label{system}
 \sum_{i=1}^s \lambda_i(u_ix_k) = \omega_k(d)\ell\frac{2}{(\alpha_k,\alpha_k)}
\end{equation}
defines a face $\mathcal{F}$ of $\mathcal{C}$.
\begin{remark}
We will show that this face has the expected dimension $|\Delta\setminus \Delta_P|$ even when $P$ is not maximal.
\end{remark}

\begin{remark}\label{beforelemma}
A line bundle on $\Parbun_G$ with a non-zero global section which does not vanish identically on $\Parbun_L(d)$ under  $\pi\circ i$ is on the face $\mathcal{F}$.

In fact a line bundle on  $\Parbun_G$ satisfies the equations \ref{system} if and only if $\gamma_{i'^*\ml}$ is trivial (cf. Proof of \cite[Theorem 5.2]{BKq}).
\end{remark}

\begin{defi}
The abelian group $\Pic'\subseteq \Pic(\Parbun_G)$ is defined as follows: It is the set of all $(\vec{\lambda},\ell)$ such that for all $1\leq j\leq s$
\begin{enumerate}
\item $\lambda_j({\beta}^{\vee})=0$ whenever $v\xrightarrow{\beta} u_j$ for a simple root $\beta$ and $v\in W^P$.
\item $\lambda_j(\theta^{\vee})=\ell$ whenever letting $v = \overline{s_\theta u_j}\in W^P$, we have  $\ell(v)>\ell(u_j)$.
\end{enumerate}

\end{defi}

\begin{defi}
Fix $1\le j\le s$ and $v\in W^P$ such that either
\begin{enumerate}
\item $v\xrightarrow{\beta} u_j$ for a simple root $\beta$, or
\item $v = \overline{s_\theta u_j}\in W^P$ and $\ell(v)>\ell(u_j)$.
\end{enumerate}

Define $v_1,\hdots, v_s$ by setting $v_k = u_k$ for $k\ne j$ and $v_j = v$. Then
set $d'$ according to the cases above:
\begin{enumerate}
\item $d' = d$, or
\item $d' = d- d(\alpha)$ where $\alpha= u_j^{-1}\theta^\vee\succ 0$.
\end{enumerate}

Clearly in case (1), the collection $(v_1,\hdots, v_s, d')$ satisfy the codimension one condition \eqref{excess}. It turns out that this also happens in case (2) {{by Lemma \ref{fraction}(2)}}. So $(v_1,\hdots,v_s,d')$ produce a $G$-invariant divisor $D$ in $\op{Parbun}_G$ as in Definition \ref{stacky2}; let us denote this divisor by $D(v,j)$.
\end{defi}
Note that the equalities in the definition of $\Pic'$ may be indexed by the basic divisors $D(v,j)$.
\begin{lemma}\label{detailed}
Let $D(v_0,j_0)$ be one of the divisors. The line bundle $\mathcal{O}(D(v_0,j_0))$ fails to be on $\Pic'$ for exactly one reason each. It violates the  equality indexed by $D(v_0,j_0)$. The cycle class of $D(v_0,j_0)$ is therefore non-zero, and $D(v_0,j_0)$ is an irreducible divisor in $\op{Parbun}_{G}$ (using Lemma \ref{relativity} for the irreducibility and codimension).
\end{lemma}
\subsection{Proof of Lemma \ref{detailed}}
\subsubsection{The case when $[D(v_0,j_0)]$ is a divisor of the first type}
 For simplicity assume $j_0=1$ and set $w=v_0$ and $w\xrightarrow{\beta_0} u_{j_0}$.
To verify the first condition in the definition of $\Pic'$, let $1\le j\leq s$ and $v\xrightarrow{\beta} u_j$ for a simple root $\beta$.

If $j=1$  then, by Theorem \ref{firstformula} to  get a non-zero value of $\lambda_j({\beta}^{\vee})$ for $[D(w,1)]$ we need $w\xrightarrow{\beta} v'$. So we have $w\xrightarrow{\beta} v'$, $v\xrightarrow{\beta} u_1$ and $w\xrightarrow{\beta_0} u_{j}$ which we show is  incompatible if $\beta\neq \beta_0$. The argument below if from \cite[Corollary 5.1]{BKiers}, which we repeat for the convenience of the reader.
\begin{enumerate}
\item If $\beta=\beta_0$, then by Theorem \ref{firstformula}, the value of $\lambda_j({\beta}^{\vee})$ is the intersection number $\langle \sigma_{u_1}^P,\hdots,\sigma_{u_s}^P\rangle_d^{\circledast_0}=1$.
\item If $\beta\neq \beta_0$, $u_1^{-1}\beta\prec 0$, $u_1^{-1}\beta_0\prec 0$,  $w^{-1}\beta\succ 0$ with $w=s_{\beta_0} u_1$, so that $u_1^{-1} s_{\beta_0}\beta\succ 0$. Now
$s_{\beta_0}\beta=\beta +m \beta_0$ with $m\geq 0$. Hence  $u_1^{-1} s_{\beta_0}\beta =  u_1^{-1}\beta + m u_1^{-1}\beta_0\prec 0$ a contradiction.
\end{enumerate}
If $j\neq 1$, then to get a non-zero value of $\lambda_j({\beta}^{\vee})$ for $[D(w,1)]$, we need $u_j\xrightarrow{\beta} v'$ for some $v'\in W^P$. This cannot happen because $v\xrightarrow{\beta} u_j$ by assumption.

We now consider the second type of conditions in the definition of $\Pic'$. Let  $v = \overline{s_\theta u_j}\in W^P$ with $\ell(v)>\ell(u_j)$. This implies
$u_j^{-1}\theta\succ 0$. Let $v=\overline{s_{\theta}u_1}\in W^P$.

We need to show that $\lambda_j({\beta}^{\vee})-\ell$ for for $[D(w,1)]$ is zero.
\begin{enumerate}
\item If $j\neq 1$, $\lambda_j({\beta}^{\vee})-\ell$ is zero because $u_j^{-1}\theta\succ 0$.
\item If $j=1$, to get a non-zero value for $\lambda_j({\beta}^{\vee})-\ell$, we need $w^{-1}\theta \prec 0$. We also have $u_1^{-1}\theta\succ 0$
and $u_1^{-1}\beta_0\prec 0$ with $w=s_{\beta_0}u_1$. Therefore $w^{-1}\theta=u_1^{-1}s_{\beta_0}\theta =u_1^{-1}\theta -\theta(\beta_0^{\vee})u_1^{-1}\beta_0$.
But $\theta(\beta_0^{\vee})\geq 0$ since $\theta$ is a dominant integral weight. Now $u_1^{-1}\theta \succ0$ and $-u_1^{-1}\beta_0\succ 0$ and hence $w^{-1}\theta\succ 0$ which contradicts  $w^{-1}\theta \prec 0$.
\end{enumerate}

\subsubsection{The case when $[D(w,1)]$ is a divisor of the second type}
For simplicity assume $j_0=1$ and set $w=v_0$ and $w=\overline{s_{\theta}u_{j_0}}$ in $W_P$. Note that $\ell(w)>\ell(u_j)$.
To verify the first condition in the definition of $\Pic'$, let $1\le j\leq s$ and $v\xrightarrow{\beta} u_j$ for a simple root $\beta$, so that $u_j^{-1}\beta\prec 0$. If $j\neq 1$, the argument is the same as the earlier case, so we assume $j=1$.  To  get a non-zero value of $\lambda_j({\beta}^{\vee})$ for $[D(w,1)]$ we need $w\xrightarrow{\beta} v'$. This implies $w^{-1}\beta\succ 0$. But $w^{-1}\beta = u_1^{-1} s_{\theta}\beta= u_1^{-1}\beta - \beta(\theta^{\vee})u_1^{-1}\theta$. Now $\beta(\theta^{\vee})\geq 0$ and $u_1^{-1}\beta\prec 0$ therefore  $w^{-1}\beta\prec 0$ a contradiction.

For the second type of condition, note that they fail for $j=1$ and are true for other $j$.

\begin{remark}
In fact the divisor $D(v_0,j_0)=\mathcal{B}(\vec{\lambda},\ell)$ in case
\begin{itemize}
\item $v_0\xrightarrow{\beta_0} u_j$ for a simple root $\beta_0$, satisfies all conditions above except for $\lambda_j({\beta_0}^{\vee})=0$, in fact, $\lambda_j({\beta_0}^{\vee})=1$.
\item $v_0 = \overline{s_\theta u_j}\in W^P$ and $\ell(v_0)>\ell(u_j)$ satisfies $\lambda_j(\theta^{\vee})=\ell-1$.
\end{itemize}
\end{remark}

\begin{corollary}\label{sparrow}
Let $D_1,\dots,D_q$ be the extremal ray generators produced as $[D(v,j)]$ in Theorem \ref{type1}.
\begin{equation}\label{biject}
 \prod_{i=1}^q \Bbb{Q} D_i \times \Pic'_{\Bbb{Q}}\leto{\sim} \Pic_{\Bbb{Q}}(\Parbun_{G})
\end{equation}
\end{corollary}
\subsection{Proof of Theorem \ref{type1}}
Theorem \ref{type1} is the special case in the following theorem with $P$ maximal parabolic.
\begin{theorem}\label{type11}
With assumptions as in the beginning of Section \ref{beginning},
\begin{enumerate}
\item[(a)] $D(v,j)$ is an irreducible divisor in $\op{Parbun}_{G}$.
\item[(b)] $h^0(\op{Parbun}_{G},\mathcal{O}(m D(v,j))) = 1$ for all positive integers $m$.
\item[(c)] $\mathbb{Q}_{\ge 0}\mathcal{O}(D(v,j))$ is an extremal ray of $\op{Pic}^+_{\mathbb{Q}}(\op{Parbun}_{G})$.
\item[(d)]  This extremal ray lies on the face $ \mathcal{F}$.
\end{enumerate}
\end{theorem}
\begin{proof}
Part (a) is proved in Lemma \ref{detailed} (we also use the computation of the cycle class which shows it is non-zero). For part (b) stronger statements hold:
\begin{equation}\label{before}
h^0(\op{Parbun}_{G},\mathcal{O}(mD)) = 1, \forall m\geq 0
\end{equation}
 (where $D=\sum D(v,j)$.)

 To prove this statement: We begin by noting that by Lemma \ref{imagemap} and Remark \ref{beforelemma}, $\mathcal{O}(mD)$ lies on $\mathcal{F}$.
 Now if $m>0$,  $\mathcal{O}(mD)$ violates the stability inequality corresponding to any $D(v,j)$ by Lemma \ref{detailed2} below.

 The argument for (c) is the same as in \cite{BKiers}: No multiple  $\mathcal{O}(mD(v,j)), m>0$ is a non-trivial sum of
efffective divisors since it has only one global section up to scalars. Part (d) follows from  Remark \ref{beforelemma}.
\end{proof}
\begin{lemma}\label{detailed2}
Let $D(v,j)$be a basic divisor and $\ml=\mathcal{B}(\vec{\lambda},\ell)\in \Pic(\Parbun_G)$ satisfy the system of equalities \ref{system} (which defines the face $\mathcal{F}$)  with $\lambda_i$ dominant at level $\ell$.
Suppose $\ml$ fails the equality in the definition of $\Pic'$ corresponding to $D(v,j)$. Then the points of $D(v,j)$ are destabilising for  $\ml$, and hence $H^0(\Parbun_G,\ml)= H^0(\Parbun_G, \ml(-D(v,j)))$.
\end{lemma}
\begin{proof}
Assume $j=1$ and assume the contrary. Hence for $\alpha_k\not\in S_P$
$$
\sum_{i=1}^s \lambda_i(u_ix_jk) = \omega_k(d)\frac{2\ell}{(\alpha_k,\alpha_k)}
$$
as well as the system
$$\sum_{i=1}^s \lambda_i(v_ix_k) \leq  \omega_k(d') \frac{2\ell}{(\alpha_k,\alpha_k)}.$$
Therefore $\lambda_1(u_1x_k)-\lambda_1(v x_j)\geq \ell \omega_k(d-d')$ for each $k$. This can be rewritten as
\begin{equation}\label{listo}
(u_1^{-1}\lambda_1-v^{-1}\lambda_1)(x_k)\geq  \omega_k(d-d')\frac{2\ell}{(\alpha_k,\alpha_k)}
\end{equation}

We divide the proof into two cases depending on the type of the divisor $D(v,1)$. We start with the first type: $v\xrightarrow{\beta} u_j$ for a simple root $\beta$ and
$v^{-1}\beta\succ 0$ and $u_1=s_{\beta}v$ and $d=d'$. This argument is the same as in \cite[Proposition 5.2]{BKiers}. Now,
$$(u_1^{-1}\lambda_1-v^{-1}\lambda_1)(x_k)= \lambda_1(\beta^{\vee})(-v^{-1}\beta(x_k))$$
We will show that this number is strictly less than zero for some $k$. Note that $v^{-1}\beta\in R^+-R_{\frl}^+$. Let $\alpha_k$  appear in $\beta$ with a
positive coefficient and hence $\alpha(x_k)\ne 0$. By assumption $\lambda_1(\beta^{\vee}) >0$ and \eqref{listo} fails.

Now we deal with the second type of divisors $D(v,1)$.  Let  $v = \overline{s_\theta u_j}\in W^P$ with $\ell(v)>\ell(u_j)$. This implies
$\alpha=u_j^{-1}\theta^{\vee}\in R^+\setminus R_{\frl}^+$. Let $v=\overline{s_{\theta}u_1}\in W^P$.

The inequality\eqref{listo} now becomes
$$\lambda_1(\theta^{\vee})(u_1^{-1}\theta(x_k))\geq  \omega_k(d-d')\frac{2\ell}{(\alpha_k,\alpha_k)}= \omega_k(d(\alpha))\frac{2\ell}{(\alpha_k,\alpha_k)}=\omega_{k}(\alpha^{\vee})\frac{2\ell}{(\alpha_k,\alpha_k)}$$

That is,
$$\lambda_1(\theta^{\vee})\alpha(x_k)\geq  \omega_{k}(\alpha^{\vee})\frac{2\ell}{(\alpha_k,\alpha_k)}$$ since $(\alpha,\alpha)=2$.
But  $\alpha(x_k)= \omega_{k}(\alpha^{\vee})\frac{2}{(\alpha_k,\alpha_k)}$ since  $(\alpha,\alpha)=2$.
\begin{remark}
To see this write $\alpha =\sum c_j\alpha_j$, then
$$\alpha^{\vee}= \alpha =\sum c_j\alpha_j= \sum c_j \frac{(\alpha_j,\alpha_j)}{2}\alpha_j^{\vee}$$
\end{remark}

The inequality\eqref{listo}  is therefore equivalent to $\lambda_1(\theta^{\vee})-\ell \geq 0$ which is false because $\lambda_1(\theta^{\vee})<\ell$ under our assumptions (given that   the equality corresponding to $D(v,1) $ fails).

\end{proof}

\subsection{}
 Recall that there is a natural map $\Omega^0(\vec{u},d)\to \Parbun_L(d)$ which restricts to  $\Omega^0(\vec{u},d)\setminus \mathcal{R}\to \Parbun_L(d) \setminus \mathcal{R}_L$. There is also a section $\Parbun_L(d)\to  \Omega^0(\vec{u},d)$ mapping $\mathcal{R}_L$ to $\mathcal{R}$.

By Zariski's main theorem applied to the birational map $\pi$ the closure of $\pi(\mathcal{R})$ has codimension $\geq 2$ in $\Parbun_G$. Using Theorems \ref{Hermann}, \ref{T1}, \ref{T2}, and Remark \ref{contracted} we get
\begin{lemma}\label{zussamen}
$\Omega^0(\vec{u},d)\setminus \mathcal{R}\to \Parbun_{G,S}$ is a open substack, and the complement is a union of the basic divisors $D(v,j)$ (and perhaps a codimension two substack).
\end{lemma}


The following is the quantum generalization of Fulton's conjecture  for arbitrary groups (the classical version for arbitrary groups was proved in \cite{BKR}):
Recall that we are assuming that $P$ is an arbitrary standard parabolic and $\langle \sigma_{u_1}^P,\hdots,\sigma_{u_s}^P\rangle_d^{\circledast_0}=1$.
 It is a corollary of Theorem \ref{type11}.
\begin{corollary}\label{FulCon} 
$h^0(\Parbun_L(d), \mathcal{O}(m\mathcal{R}_L))=1$ for all $m\geq 0$, equivalently,
$h^0(\Parbun_L(d)-\mathcal{R}_L, \mathcal{O})=1$ for all $m\ge 0$.
\end{corollary}
\begin{proof}
Using and the diagram \eqref{Dbasic}, it suffices to show that
$h^0(\Omega^0(\vec{u},d)- \mathcal{R}, \mathcal{O})=1$ (since we have a section $i$ in Diagram \ref{Dbasic}). But
$$\Omega^0(\vec{u},d)- \mathcal{R}=\Parbun_G- D$$
 up to codimension two substacks, where $D$ is the union of the basic divisors $D(v,j)$, by Lemma \ref{zussamen}, and therefore the result follows from \eqref{before}.
\end{proof}
We refer the reader to \cite[Remark 8.6]{BKq} for a history of earlier work on Fulton's conjecture.

\section{Induction}

\begin{defi}\label{stom}
Le $\Pic(\Parbun_L(d) \setminus \mathcal{R}_L)^{\deg =0}$ (similarly $\Pic(\Parbun_L(d))^{\deg =0}$) be the set of line bundles $\mathcal{M}$ such that $\gamma_{\mathcal{M}}=0$ (see Defn.\ref{indexx}).

Let  $\Pic^{\deg=0}(\Omega^0(\vec{u},d)\setminus \mathcal{R})$ be the set of line bundles whose pullbacks to $\Parbun_L(d) \setminus \mathcal{R}_L$ have degree $0$ (see Remark \ref{beforelemma}).
\end{defi}

We begin by noting a corollary of Lemma \ref{detailed2} and Corollary \ref{sparrow}:
\begin{corollary}\label{salvador}
The face $\mf$ of $\mathcal{C}$ is a product, with $\Pic'^{+,\deg=0}_{\Bbb{Q}}=\Pic'\cap\mf$:
\begin{equation}\label{biject1}
\mf\leto{\sim} \prod_{i=1}^q \Bbb{Q}_{\geq 0} D_i \times \Pic'^{+,\deg=0}_{\Bbb{Q}}
\end{equation}
\end{corollary}

\begin{lemma}\label{zussamen2}
\begin{enumerate}
\item The restriction mapping $\Pic'\to \Pic(\Omega^0(\vec{u},d)\setminus \mathcal{R})$ is an isomorphism.
\item $\Pic^{+,\deg=0}(\Omega^0(\vec{u},d)\setminus \mathcal{R})$  is in bijection with $\Pic'^{+,\deg=0}$. See Defn. \ref{stom} for the definition of
$\Pic^{+,\deg=0}(\Omega^0(\vec{u},d)\setminus \mathcal{R})$.
\end{enumerate}
\end{lemma}
\begin{proof}
By Lemma \ref{zussamen}, any element in the kernel is linear combination of the $D(v,j)$, Evaluating such a linear dependence on suitable conditions in the definition of $\Pic'$ we find that the linear combination but be trivial (i.e., use \eqref{biject}). For the surjection, we can  choose an extension as a line bundle \cite[Corollary 3.4]{Heinloth} to the smooth stack $\Parbun_{G,S}$ and then use \eqref{biject}. This proves the first part.

For the second part, we clearly have a injective map $\Pic'^{+,\deg=0}\to \Pic^{+,\deg=0}(\Omega^0(\vec{u},d)\setminus \mathcal{R})$. We show this map is surjective. Let $\ml$ be a line bundle on $\Omega^0(\vec{u},d)\setminus \mathcal{R}$ which is effective with a section $s$ and satisfies $\deg=0$ (i.e., its restriction to $\Parbun_L(d)-\mathcal{R}_L$ has this property). By considering the closure of the Weil divisor of the section $s$ we can extend $\ml$ to an effective line bundle on $\Parbun_{G,S}$ (Lemma \ref{zussamen}), but this line bundle may not be in $\Pic'$. It still satisfies $\deg=0$ by Remark \ref{beforelemma}.  Using Lemma \ref{detailed2} we can subtract a suitable linear combination of the boundary divisors and get an element of  $\Pic'^{+,\deg=0}$ with
the desired properties.

\end{proof}
\begin{lemma}\label{pruf}
\begin{enumerate}
\item[(A)] $\Pic'^{+,\deg=0}$ is isomorphic to $\Pic^+(\Parbun_L(d) \setminus \mathcal{R}_L)$.
\item[(B)]  The natural map $\Pic^+(\Parbun_L(d))\to\Pic^+(\Parbun_L(d) \setminus \mathcal{R}_L)$ is a surjection. If $D_0$ is any divisor on $\Parbun_L$ whose support is contained in $\mathcal{R}_L$ then the image of $\mathcal{O}(D_0)$ in $\Pic^+(\Parbun_L(d) \setminus \mathcal{R}_L)$ is the trivial line bundle.

\end{enumerate}
\end{lemma}
Part (A) follows from Lemma \ref{zussamen2} and  Proposition \ref{comparison} which gives the following
\begin{lemma}
\begin{enumerate}
\item $\Pic(\Omega^0(\vec{u},d)\setminus \mathcal{R})$ is isomorphic to $\Pic(\Parbun_L(d) \setminus \mathcal{R}_L)$.
\item  $\Pic^{+,\deg =0}(\Omega^0(\vec{u},d)\setminus \mathcal{R}) $ is isomorphic to $\Pic^+(\Parbun_L(d) \setminus \mathcal{R}_L)$.
\end{enumerate}
\end{lemma}

For part (B), we use that For (2) we use that there is a factoring $$\Pic'\to \Pic(\Parbun_L(d)) \to \Pic(\Parbun_L(d) \setminus \mathcal{R}_L).$$

Finally note that we get a surjection
\begin{align}\label{surj}
\Pic^+(\Parbun_L(d)) \twoheadrightarrow \Pic'^{+,\deg=0}=\mf_{\op{II}}.
\end{align}

\begin{lemma}\label{twostepprocess}
\begin{enumerate}
\item[(1)]   $\Pic^+(\Parbun_{L}(d))$ is in bijection with $\Pic^+(\Parbun_{L}(0))$.
\item[(2)]  $\Pic^+_{\Bbb{Q}}(\Parbun_{L}(0))\leto{\sim}\Pic^+_{\Bbb{Q}}(\Parbun_{L'})$, where
$L'=[L,L]$ (which will be shown to be simply connected and semisimple, and hence breaks up as product over the simple components).
\end{enumerate}
\end{lemma}
These statements are proved in Sections \ref{partie1} and \ref{partie2}. We make the maps explicit in Section \ref{algorithm}, so that we get an explicit mapping \eqref{inducto} as in the introduction with $\mf_{\op{II}}=\Pic'^{+,\deg=0}_{\Bbb{Q}}$.

\subsection{Changing the degree}\label{partie1}
 We now show that $\Pic^+(\Parbun_{L}(d))$ is in bijection with $\Pic^+(\Parbun_{L}(0))$ by a sequence of degree change isomorphisms and hence prove part (1) of Lemma \ref{twostepprocess}. The following is a generalization of  \cite{BKq}*{Lemma 7.2} to parabolics that are not necessarily maximal.

\begin{lemma}\label{dynkin}
Let $P$ be any standard parabolic, and $S_P = \Delta \setminus \Delta(P)$. For each $\alpha_i\in S_P$, there exists a $\mu_i$ belonging to the coroot lattice $Q^\vee$ such that
\begin{enumerate}
\item[($a$)] $0\le \alpha(\mu_i)\le 1$ for all roots $\alpha\in R_{\mathfrak{l}}^+$,
\item[($b$)] $|\omega_i(\mu_i)| = 1$, and
\item[($c$)] $\omega_j(\mu_i)=0$ for all $\alpha_j\in S_P$ such that $j\ne i$.
\end{enumerate}
\end{lemma}
This lemma is proved in Section \ref{dynkinproof}. Using this lemma we produce some isomorphisms $\Parbun_L(d)\to \Parbun_L(d')$ for lower degrees $d$; a composition will give us the desired isomorphism
$\Parbun_L(d)\to \Parbun_L(d')$.

Let $\mu\in Q^\vee$ satisfy condition (a) in Lemma \ref{dynkin}. Let $d=(a_i)_{\alpha_i\in S_P}\in H_2(G/P)$. Let $d'= (a_i-\omega_i(\mu))_{\alpha_i\in S_P}$. We  define a map
\begin{equation}\label{changeD}
\tau_{\mu}:\op{Parbun}_L(d)\to\op{Parbun}_L(d')
\end{equation} generalizing the construction in \cite[Lemma 7.3]{BKq}. Note that \cite[Lemma 7.3]{BKq} has a sign typo in the degree computation.

For simplicity assume $p_1=0$. Let $\ell_{\mu}:\Bbb{C}^*\to H\subseteq L$ be the map corresponding to $\mu$. Pick a point $\tilde{\ml}=(\mathcal{L}; \bar l_1,\hdots, \bar l_s)\in \op{Parbun}_L(d)$. Fix a trivialization $s$ of $\ml$ in a neighborhood of $0$ so that $s(0)B_L=\bar{l}_1.$

We define a new principal bundle $\ml_{\mu}$ as follows: It coincides with $\ml$ outside of zero. Sections of $\ml_{\mu}$ in a neighborhood of $0$ are sections $sa(z)$ of $\ml$ so that $\ell_{\mu}(z)a(z)$ is regular at $z=0$. The section
$s'= s\ell_{\mu}(z)^{-1}$ is a local trivialization of $\ml_{\mu}$ in a neighborhood of $0$.

If we had chosen a different trivialization $sc(z)$ of $\ml$ with $c(0)\in B_L$, then the meromorphic sections do not change: $sa(z)=sc(z)^{-1}a(z)$, so we would examine if $\ell_{\mu}(z)c(z)^{-1}a(z)$  is holomorphic at zero. But  $\ell_{\mu}(z)c(z)^{-1}\ell_{\mu}(z)^{-1}$ is holomorphic at the origin because of assumption (a) in Lemma \ref{dynkin}. But  the limit of $\ell_{\mu}(z)c(z)^{-1}\ell_{\mu}(z)^{-1}$ as $z\to 0$ is not necessarily an element of $B_L$.
The new  $s'\times^L \C_{-\omega_i}$ as a meromorphic section of $\ml_{\mu}\times^L \C_{-\omega_i}$ equals the old $s\times^L \C_{-\omega_i}$ times $z^{\omega_i(\omega)}$. This gives that the degree of $\ml_{\mu}$ is $d'$.

We still to give a well defined parabolic structure to $\ml_{\mu}$ at $p_1$: Let
$E_{\mu}\subseteq L$ be the set of limits
$$\lim_{z\to 0}\ell_{\mu}(z)c(z)\ell_{\mu}(z)^{-1},\ c(z)\in L[[z]], c(0)\in B_L.$$

By construction we get a well defined element in $({\ml}_{\mu})_{p_0}/E_{\mu}$. As in \cite[Page 1348]{BKq}, $E_{\mu}$ is a Borel subgroup of $L$ with Lie algebra spanned by
$\frh$, and root spaces $\frg_{\beta}$ when $\beta\in R_{\ell}^{+}$ with $\beta(\omega)=0$, and root spaces $\frg_{-\beta}$ when $\beta\in R_{\ell}^{+}$ with $\beta(\omega)=1$ (adding any other $\frg_{\beta}$ will produce non-solvable factors, check that this is solvable). Write $E_{\mu}=w_{\mu}B_Lw_{\mu}^{-1}$ for $w_{\mu}\in W_L$. Therefore the element $\tilde{l}_1=s'w_{\mu}\in (\ml_{\mu})_0/B_L$ is canonically defined.

\begin{defi}
The mapping \eqref{changeD} takes  $\tilde{\ml}=(\mathcal{L}; \bar l_1,\hdots, \bar l_s)\in \op{Parbun}_L(d)$ to $\tilde{\ml}_{\mu}=(\mathcal{L}_{\mu}; \overline{\tilde{l}}_1, \bar l_2, \hdots, \bar l_s)\in \op{Parbun}_L(d')$
\end{defi}
It can be verified that \eqref{changeD} is an isomorphism of stacks, hence is an isomorphism on Picard groups.

\subsection{Proof of Lemma \ref{dynkin}}\label{dynkinproof}
The simple factors of the standard Levi $L\subset P$ are indexed by the connected components of the Dynkin diagram of $L$ as a sub-diagram of that for $G$. These connected components partition the roots of the Levi, $\Phi_L$, into a maximal disjoint union of orthogonal root systems:
$$
\Phi_L = \bigsqcup_{k=1}^m \Phi_k,
$$
where $\alpha\in \Phi_k$ and $\beta\in \Phi_\ell$ implies $\alpha(\beta^\vee) = 0$ for all $k\ne \ell$, and
such that each $\Phi_k$ cannot further be so decomposed.

Now consider $\alpha_i\in S_P$. Though it belongs to none of the $\Phi_k$, it might be adjacent (in the Dynkin diagram) to several of the $\Phi_k$. Set $\Delta_k = \Delta\cap \Phi_k$. Define
$$
\Delta_L^{(i)} =  \bigsqcup_{k \text{ s.t. }\exists \beta\in \Phi_k \text{ s.t. } \beta(\alpha_i^\vee)\ne0} \Delta_k \sqcup \{\alpha_i\}.
$$
Then $\Delta_L^{(i)}$ is the base of a Levi sub-root system of $G$; it is the minimal standard Levi containing the $\alpha_i$ root subgroup $G_{\alpha_i}$ and any simple factors of $L$ that do not commute with $G_{\alpha_i}$. It is never a subgroup of $L$, and it often does not contain all of $L$. Let us call this new Levi $M_i$. $M_i$ is always simple.

\begin{figure}
\begin{tikzpicture}
\draw[-] (0,0) -- (5,0) -- (5.7,0.7);
\draw[-] (5,0) -- (5.7,-0.7);

\draw[red,-,thick] (0,0.05) -- (1,0.05);
\draw[red,-,thick] (3,0.05) -- (4,0.05);

\draw[blue,-,thick] (3,-0.05) -- (5,-0.05);
\draw[blue,-,thick] (5,-0.07) -- (5.7,0.63);

\fill[color=black] (0,0) circle (3pt);
\node at (0,0) [circle,red,fill,inner sep=1.9pt]{};
\fill[color=black] (1,0) circle (3pt);
\node at (1,0) [circle,red,fill,inner sep=1.9pt]{};
\fill[color=black] (2,0) circle (3pt);
\node at (2,0) [circle,white,fill,inner sep=1.9pt]{};
\fill[color=blue] (3,0) circle (5pt);
\node at (3,0) [circle,red,fill,inner sep=1.9pt]{};
\fill[color=blue] (4,0) circle (5pt);
\node at (4,0) [circle,red,fill,inner sep=1.9pt]{};
\fill[color=blue] (5,0) circle (5pt);
\node at (5,0) [circle,green,fill,inner sep=1.9pt]{};

\fill[color=blue] (5.7,0.7) circle (5pt);
\node at (5.7,0.7) [circle,red,fill,inner sep=1.9pt]{};
\fill[color=black] (5.7,-0.7) circle (3pt);
\node at (5.7,-0.7) [circle,white,fill,inner sep=1.9pt]{};

\end{tikzpicture}
\caption{Red nodes depict the original simple roots for the Levi $L$. The green node is $\alpha_i$. The blue nodes comprise the simple roots for $M_i$. The remaining nodes make up $S_P\setminus \{\alpha_i\}$.}\label{dynk}
\end{figure}
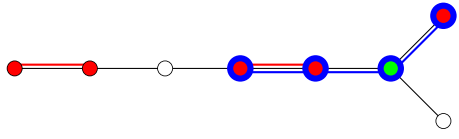

Inside $M_i$ there is a \emph{maximal} Levi subgroup $L_i\subset M_i$ obtained by removing the simple root $\alpha_i$. See Figure \ref{dynk} for an example of the different Levi subgroups.

We may apply \cite{BKq}*{Lemma 7.2} to the pair $L_i\subset M_i$ and obtain a $\mu_i\in Q_{M_i}^\vee$ satisfying
\begin{enumerate}
\item[($a'$)] $0\le \alpha(\mu_i)\le 1$ for all positive roots $\alpha$ of $M_i$, and
\item[($b'$)] $|\omega_i^{M_i}(\mu_i)| = 1$.
\end{enumerate}
Viewing $Q_{M_i}^\vee\subset Q^\vee$, we may write $\mu_i = \sum_{\alpha_k\in \Delta_L^{(i)}} c_k\alpha_k^\vee$. Then property ($b'$) says that $c_i = \pm 1$, which is exactly what is required to satisfy ($b$).

Let us now verify that $\mu_i$ satisfies both ($a$) and ($c$). For ($a$), let $\beta$ be any positive root in $\Phi_L$. If $\beta$ is not a root of $M_i$, then $\beta\in \Phi_k$ for some $\Phi_k$ orthogonal to both $\alpha_i$ and all roots of $L\cap M_i$. Therefore $\beta(\mu_i) = 0$. Otherwise, $\beta$ is a root of $M_i$, and ($a'$) implies ($a$) for the root $\beta$.

Finally, if $\alpha_j\in S_P$ and $\alpha_j\ne \alpha_i$, then clearly $\alpha_j\not\in \Delta_L^{(i)}$. So $\omega_j(\mu_i)=0$.

\begin{remark}
Here is another way to think about $\mu_i$. Condition ($a$) says that the restriction of $\mu_i$ to $\mathfrak{h}_{\mathfrak{l}}$ must live in the Weyl alcove (really the product of Weyl alcoves, one for each simple factor of $L$). Adjusting $\mu_i$ by $Q^\vee_L$ translations or by the $W_L$ action will not affect ($b$) or ($c$), so first set $\mathring \mu_i = \alpha_i^\vee$, then use the action of the Levi affine Weyl group on $\mathfrak{h}_{\mathfrak{l}}$ to obtain $\mu_i$, the unique alcove representative for $\mathring \mu_i$. If one starts with $\mathring \mu_i = -\alpha_i^\vee$ instead, the $\mu_i$ is clearly different also. However, these are the only two possibilities. 
\end{remark}

\subsection{Proof of (2) of Lemma \ref{twostepprocess}}\label{partie2}

\begin{defi}\label{def:n(L)}
The  map $Z^0(L)\to L/L'$ is an isogeny, and the kernel $Z^0(L)\cap L'\subseteq L'$ is a finite abelian group.
Let $n(L)$ be the order of the kernel.
\end{defi}
Since $L/B_L=L'/B_{L'}$, we have a morphism
  $i:\Parbun_{L'}\to \Parbun_L(0)$, which gives rise to $\Pic(\Parbun_L(0))\to \Pic(\Parbun_{L'})$. This is surjective on objects (\cite[Lemma 7.1]{BKq}).
 \begin{lemma}\label{together}
 There is a bijection of cones
 $$\Pic_{\Bbb{Q}}^{\op{index}=0}(\Parbun_L(0))\to \Pic_{\Bbb{Q}}(\Parbun_{L'})$$

Where the index of a line bundle $\ml$ on  $\Parbun_L(0)$ equals $\gamma_{\ml}$ (Defn. \ref{indexx}).
 \end{lemma}

 \begin{proof}
 We begin with the injection.

 Suppose $\ml$ is a  line bundle on $\Parbun_L(0)$ whose restriction to $\Parbun_{L'}$ is trivial, and let $s$ be a trivialization.
 Therefore for all principal $L'$ bundles $x$, we have $s(x)\in \ml_{i(x)}$. Now every point in $\Parbun_L(0)$ is of the form $i(x)$ for some $x$.
 If $i(x)\to i(x')$ is an isomorphism, then we want $s(x)$ to be carried to $s(x')$ under this isomorphism. Note that this would be obvious if the map $i(x)\to i(x')$
 came from a morphism $x\to x'$. However it can be seen (see \cite[Lemma 7.1]{BKq}) that there is a $z_o\in Z^0(L)$ such that  $i(x)\to i(x')$ composed with the action of $z_o$ comes from a map $x\to x'$.
 This is sufficient because $z_o\in Z^0(L)$ acts trivially by the index criterion.

 For the surjection, we show that if $\mk$ is in $\Pic(\Parbun_{L'})$, then $\mk^{\tensor n(L)}$ (see Dfn. \ref{def:n(L)}) is the restriction of a index zero line bundle on $\Parbun_L(0)$.
  Define a line bundle on $\Parbun_L(0)$ whose fiber at $i(x)$ is $\mk_x$. If $i(x)\to i(x')$ is a map, then we want a map $\mk_x\to\mk_{x'}$ which would again be given if the map $i(x)\to i(x')$
 came from a morphism $x\to x'$. Again, there is  $z_o\in Z^0(L)$ such that  $i(x)\to i(x')$ composed with the action of $z_o$ comes from a map $x\to x'$. This $z_0$ is well defined up-to elements
 in $Z_0(L)\cap L'$ in  \cite[Lemma 7.1]{BKq}, which acts trivially on $\mk^{\tensor n(L)}$. Therefore  $\mk^{\tensor n(L)}$ lifts to a line bundle on $\Parbun_L(0)$, which is verified to have index $0$.
  \end{proof}

The following implies (2) from Lemma \ref{twostepprocess}. It follows from \cite[Lemma 7.1]{BKq}.
\begin{lemma}
 The map
 $\Pic_{\op{index}=0}(\Parbun_L(0))\to \Pic(\Parbun_{L'})$ preserves global sections. That is, if $\mn\in \Pic_{\op{index}=0}(\Parbun_L(0))$, then
 $H^0(\Parbun_L(0),\mn)\leto{\sim} H^0(\Parbun_{L'},i^*\mn)$. This (together with Lemma \ref{together}) shows that
 $\Pic^+_{\Bbb{Q}}(\Parbun_{L}(0))\leto{\sim}\Pic^+_{\Bbb{Q}}(\Parbun_{L'})$.
 \end{lemma}

\subsection{Formulas for pullback of line bundles}

\begin{defi}\label{parfactors}
Given a weight $\lambda:B_L\to \Bbb{C}^*$ (or just $T_L\to \Bbb{C}^*$) and a $p_k$, we can form the line bundle $\ml_k(\lambda)$ on ${\Parbun}_L$ whose fiber over $\tilde{\ml}=(\mathcal{L}; \bar l_1,\hdots, \bar l_s)$ equals the fiber of
$\ml_{p_k}\times^{B_L}\Bbb{C}_{-\lambda}\to \ml_{p_k}/B_L$ over $\bar{l}_k$.

Similarly if given $\lambda:B_{L'}\to \Bbb{C}^*$, we have a corresponding line bundle  $\ml_k(\lambda')$ on ${\Parbun}_L$.
\end{defi}
 We now study the pullbacks of line bundles $\ml_1(\lambda)$ (see Defn. \ref{parfactors}) under the map \eqref{changeD}.
\begin{lemma}\label{step1}
$\tau_{\mu}^*\ml_1(\lambda)=\ml_1(w_{\mu}\lambda)$
\end{lemma}
\begin{proof}
This isomorphism sends the element $s\times 1$ where $1\in \Bbb{C}_{-w_{\mu}\lambda}$  to $s\ell_{\mu}^{-1}w_{\mu}\times 1$, where $1\in \Bbb{C}_{-\lambda}$. We just need to show that this is well defined, i.e., we want this map to also take
$sb\times 1$ to $sb\ell_{\mu}^{-1}w_{\mu}\times 1= s\ell_{\mu}^{-1}w_{\mu}\cdot (w_{\mu}^{-1}\ell_{\mu}(z)b\ell_{\mu}^{-1}(z)w_{\mu}) \times 1$. Therefore we need
$\lambda(w_{\mu}^{-1}\ell_{\mu}(z)b\ell_{\mu}^{-1}(z)w_{\mu})=w_{\mu}\lambda(b)$, which is clear.
\end{proof}

\begin{remark}\label{step0}
Let $\lambda:T_{L'}\to \Bbb{C}^*$, we can extend (after scaling $\lambda$ by $\mu=\lambda^{n(L)}$) to a $\widetilde{\mu}:T\to\Bbb{C}^*$, so that $\lambda$ is trivial on $Z^0(L)$. The map in Lemma \ref{together} pairs the line bundle $\ml_k(\mu)$ on ${\Parbun}_L'$ with the line bundle $\ml_k(\widetilde{\mu})$ on $\Pic_{\Bbb{Q}}^{\op{index}=0}(\Parbun_L(0))$.
\end{remark}

\subsubsection{}Let $V=T(G/P)$ as a representation of $L$.
Then we get determinant of cohomology line bundles $D(V)$ on $\Bun_L$ as follows: The fiber of $D(V)$ over $\ml$ is the determinant of cohomology of the vector bundle $\ml\times^L V$ on $\Bbb{P}^1$, i.e.,
$$\det H^0(\Bbb{P}^1, \ml\times_L V)^*\tensor \det H^1(\Bbb{P}^1, \ml\times_L V).$$ By pullback we get line bundles $D(V)$ on ${\Parbun}_L$ and $\Parbun_{L'}$.

\begin{defi}
Let $L'_k$ be the simple factors of $L'$, and let $\mathcal{M}_k$ be the positive generators of the Picard groups of $\Bun_{L'_k}$. These give line bundles $\mathcal{M}_k$ on $\Parbun_{L'}$. Let $d_k$ be the Dynkin indices of the maps $L'_k\to \SL(V)$.
\end{defi}

\begin{lemma}\label{step4}
The line bundle $D(V)$ on $\Parbun_{L'}$  equals $\prod_k \mathcal{M}_k^{d_k}$.
\end{lemma}

Under the isomorphism from Lemma \ref{together}
$$\Pic_{\Bbb{Q}}^{\op{index}=0}(\Parbun_L(0))\to \Pic_{\Bbb{Q}}(\Parbun_{L'})$$
induced from a map $\Parbun_{L'}\to\Parbun_L(0)$. Note first that we have a line bundle
$D(V)$ on ${\Parbun}_{L'}$, and also on $\Pic_{\Bbb{Q}}(\Parbun_L(0))$. The index of the latter bundle is given by $x_i$ in the Lie algebra of $Z^0$ ($\alpha_i\in S_P$) acting by multiplication by $\rk V$ (since the degree is zero).
\begin{lemma}\label{step3}
The line bundle $D(V)$ on $\Parbun_{L'}$ corresponds to the line bundle
$D(V)\tensor \mathcal{L}_1(\chi_e)$ on $\Parbun_L(0)$, where $\chi_e=2\rho-2\rho^L$.
\end{lemma}
We only need that the index of $D(V)\tensor \mathcal{L}_1(\chi_e)$ is zero. This follows from the computation in \cite[Prop 3.15]{BKq}. The last step in the computational side is a formula for the pullback $\tau_{\mu}^*D(V)$ of $D(V)$ via the map \eqref{changeD}.

Generalizing \cite[Lemma 7.5]{BKq}, we have, by the same proof,
\begin{lemma}\label{step2}
 $\tau_{\mu}^*D(V)\tensor D(V)^{-1}= \mathcal{L}_1(\nu)$ where (note that $\mathcal{L}_1$ has a sign)
 $$\nu= \sum_{\alpha\in R^+-R_{\ell}^+} \alpha(\mu)\alpha.$$
\end{lemma}

\subsection{Putting the explicit maps together}
We have seen that $\Pic^+_{\Bbb{Q}}(\Parbun_{L}(d))$ is in bijection with $\Pic^+_{\Bbb{Q}}(\Parbun_{L}(0))$ which is in bijection with  $\Pic^+_{\Bbb{Q}}(\Parbun_{L'})$, which breaks up as a product.
We also have a  map from $\Pic^+_{\Bbb{Q}}(\Parbun_{L}(d))$ to $\mf$ which maps to the $\Pic'^{+,\deg=0}$ factor.

To make the induction map explicit we need to describe the composite map $\Pic^+_{\Bbb{Q}}(\Parbun_{L'})\to \Pic'^{+,\deg=0}$. We will do this in several steps.
\begin{enumerate}
\item The line bundles $\ml_k(\lambda)$ on $\Parbun_{L'}$ are first induced to elements of $\Pic^{\op{index}=0}_{\Bbb{Q}}(\Parbun_{L}(0))$ using Remark \ref{step0}. The corresponding elements of
$\Pic(\Parbun_{L}(d))$ are identified using Lemma \ref{step1} (it is again of the form $\ml_k(\lambda)$). The element of $\Pic'^{+,\deg=0}$ is then computed by identifying the corresponding line bundle on  $\Omega^0(\vec{u},d)\setminus \mathcal{R}$ and adding the right correction of boundary divisors to have it lie on the face $\mf$.
\item We need to induce the line bundles $\mathcal{M}_k$ on $\Parbun_{L'}$. To do this use the fact that all irreducible components of the ramification divisors induce to zero. We first write the ramification line bundle on $\Parbun_{L}(d)$
in the form given in \cite[Prop 3.15]{BKq}. Note that the determinant of cohomology factors in this expression are of the form $D(V)$ for $V=T(G/P)_e$, and we know how to write the corresponding expression in $\Pic^+_{\Bbb{Q}}(\Parbun_{L'})$
using Lemmas \ref{step2} and \ref{step3}. Using Lemma \ref{step4} this line bundle becomes a product over the simple factors of $L'$. Each of the factors features a determinant of cohomology element (raised to a factor) and finite parabolic contributions. Each of the factors induces to zero, and hence we obtain formulas for inductions of $\mathcal{M}_k$.
\end{enumerate}
\section{The induction algorithm}\label{algorithm}

As before, we fix the data of a face $\mathcal{F}$: a maximal parabolic $P$, Weyl group elements $u_1,\hdots,u_s\in W^P$, and degree $d$ such that $\langle \sigma_{u_1}^P,\hdots,\sigma_{u_s}^P\rangle_d^{\circledast_0} = 1$. Let $L_1,\hdots, L_m$ be the simple factors of $L' = [L,L]$.

The induction map $\op{Ind}$ ultimately goes from
$$
\pic^+_\mathbb{Q}(\op{Parbun}_{L'}) \simeq \pic^+_\mathbb{Q}(\op{Parbun}_{L_1}) \times \cdots \times \pic^+_\mathbb{Q}(\op{Parbun}_{L_m})
$$
to
$$
\mathcal{F}\subset \pic^+_\mathbb{Q}(\op{Parbun}_G),
$$
and the type II extremal rays of $\mathcal{F}$ are included in the images of the extremal rays of the separate factors $\pic^+_\mathbb{Q}(\op{Parbun}_{L_k})$.

Each cone $\pic^+_\mathbb{Q}(\op{Parbun}_{L_k})$ is defined inside the vector space
$$
(X^*(T_{L_k})_\mathbb{Q})^s\times \mathbb{Q}\mathcal{M}_k.
$$
Rather than describe $\op{Ind}$ on all of $\pic^+_\mathbb{Q}(\op{Parbun}_{L_k})$, we first describe the operation on $(X^*(T_{L_k})_\mathbb{Q})^s$, then indicate how to find $\op{Ind}(\mathcal{M}_k)$.

\subsection{Induction of weights}

$\op{Ind}$ is the composition of several linear maps:
$$
\pic^+_\mathbb{Q}(\op{Parbun}_{L'}) \xrightarrow{\op{ext}}
\pic^+_\mathbb{Q}(\op{Parbun}_L(0)) \xrightarrow{(\tau_\mu^*)^{\pm d}}
\pic^+_\mathbb{Q}(\op{Parbun}_L(d)) \xrightarrow{p}\pic^{'+,\deg=0}=\mf_{\op{II}},
$$
where the last map $p$ is the surjection (\ref{surj}). Restricted to $(X^*(T_{L'})_\mathbb{Q})^s$, the map $\op{ext}$ operates the same in each factor by taking a tuple of weights
$$
(\mu_1, \hdots, \mu_s)
$$
to $(\lambda_1,\hdots,\lambda_s)\in (X^*(T)_\mathbb{Q})^s$ where $\lambda_i$ is the unique weight satisfying $\lambda_i|_{T_{L'}} \equiv \mu_i$ and $\lambda_i(x_j)=0$ for every $\alpha_j \in S_P = \Delta \setminus \Delta_P$.

To describe the linear map $\tau_\mu^*$, first identify $\mu\in Q^\vee$ as in Lemma \ref{dynkin}. Then there is an associated Weyl group element $w_\mu=w_0^\mu w_0^L\in W_L$, where $w_0^L$ is the longest element of the Levi Weyl subgroup $W_L$, and $w_0^\mu$ is the longest element of the smaller subgroup $W_\mu$ that stabilizes $\mu$. The action of $\tau_\mu^*$ on weights is to twist by $w_\mu$ in only the first factor (cf. Lemma \ref{step1}):
\begin{align}\label{tau-mu-weights}
\tau_\mu^* : (\lambda_1,\lambda_2,\hdots,\lambda_s) \mapsto (w_\mu \lambda_1, \lambda_2,\hdots,\lambda_s)
\end{align}

Finally, for computations it is helpful to factor $p$ as follows:
\begin{align}\label{verdant}
\pic^+_\mathbb{Q}(\op{Parbun}_L(d)) \xrightarrow{\vec{u}} \mathbb{Q} \mathcal{F} \xrightarrow{p_2} \mathbb{Q}\mf_{\op{II}},
\end{align}
where $\vec{u}$ is the map twisting each weight by the corresponding $u_i$, and $p_2$ is the second projection map from (\ref{biject1}). The explicit formula for $p$, acting on a tuple of weights, is therefore
\begin{align}\label{twist-and-correct}
(\lambda_1,\hdots,\lambda_s) \mapsto (u_1\lambda_1, \hdots, u_s\lambda_s; 0) - \sum_{j=1}^s \sum_{v\xrightarrow{\alpha_i}u_j}  u_j\lambda_j(\alpha_i^\vee)\vec{\mu}(D(v,j)) - \sum_{j=1}^s \sum_{v=\overline{s_\theta u_j}> u_j} -u_j\lambda_j(\theta^\vee) \vec{\mu}(D(v,j)).
\end{align}

\subsection{Induction of determinants of cohomology}

Every component of the ramification divisor $\mathcal{O}(\mathcal{R}_L)$ on $\Parbun_L(d)$ induces to $0$, and this yields enough information to calculate each $\op{Ind}(\mathcal{M}_k)$.

Inside $\pic^+(\Parbun_L(d))$, $\mathcal{O}(\mathcal{R}_L)$ is identified with
$$
(\chi_{u_1}-\chi_{e},\chi_{u_2}, \hdots, \chi_{u_s}) \otimes D(V),
$$
where $V = \mathfrak{g}/\mathfrak{p}$. First we use the above formulas to find its isomorphic preimage in $\pic^+_\mathbb{Q}(\Parbun_{L'})$. The only new ingredient needed here is the formula for $\tau_\mu^*(D(V))$ (cf. Lemma \ref{step2}):
\begin{align}\label{dv}
\tau_\mu^*(D(V)) = (\nu,0,\hdots,0)\otimes D(V),
\end{align}
where $\nu = \sum_{\alpha\in R^+\setminus R_L^+} \alpha(\mu)\alpha$. Possibly the inverse formula will be required:
\begin{align}\label{inverse-dv}
(\tau_\mu^*)^{-1} (D(V)) = (-w_\mu^{-1}\nu,0,\hdots,0) \otimes D(V)
\end{align}

 The inverse of $\op{ext}$ naturally restricts the weights of $T$ as weights of $T_{L'}$; for this part of the process it is helpful to know $\nu|_{T_{L'}}$, see Proposition \ref{nu-restrict}. Moreover,
 $$
 \op{ext}^{-1} D(V) = \mathcal{M}_1^{d_1}\otimes \hdots \otimes \mathcal{M}_m^{d_m},
 $$
 where each $d_k$ is the Dynkin index of the representation $\mathfrak{g}/\mathfrak{p}$ for the simple factor $L_k$.

Finally, if $(\mu_1,\hdots,\mu_s)\in X^*(T_{L_k})$ are the weights of the component of $\mathcal{O}(\mathcal{R}_L)$ restricted to $L_k$, then
$$
\op{Ind}(\mathcal{M}_k) =- \frac{1}{d_k} \op{Ind}(\mu_1,\hdots,\mu_s).
$$

\subsection{Summary}

Here is a step-by-step summary of the algorithm:

\begin{enumerate}

\item Identify $\mu\in Q^\vee$ as in Lemma \ref{dynkin}.

\item Calculate $k_0 = \omega_P(\mu)$, which is $\pm1$ by design. The map $\tau_\mu^*$ will be applied a total of $dk_0$ times.

\item Calculate $w_\mu$: this is the product $w_0^\mu w_0^L$, where $w_0^L$ is the longest element of the Levi Weyl subgroup, and $w_0^\mu$ is the longest element of the smaller subgroup that fixes $\mu$.

\item Calculate Dynkin indices (any of several equivalent ways, see e.g. \cite[Appendix A]{KumarConformalBlocks}): for each simple Levi factor $L_k\subset L$, the Dynkin index $d_k$ of the representation $\mathfrak{g}/\mathfrak{p}$ can be calculated as
$$
d_k = \frac{1}{2} \sum_{\alpha\in R^+\setminus R_L^+} \alpha(\theta_k^\vee)^2,
$$
where $\theta_k$ is the highest root in the simple root system associated with $L_k$. Alternatively, see Corollary \ref{dual-cox}.

\item Calculate $\nu = \sum_{\alpha\in R^+\setminus R_L^+} \alpha(\mu)\alpha$, or at least the restriction $\nu|_{L'}$. If $\alpha_j$ is a simple root for the simple component $L_k$, then
$$
\nu(\alpha_j^\vee) = \left\{
\begin{array}{ll}
d_k, & \alpha_j(\mu)=1\\
0, & \alpha_j(\mu)=0
\end{array}
\right.
$$
Moreover, the restriction of $\nu$ to a simple factor of $L'$ is either identically $0$ or a multiple $(d_k)$ of a fundamental weight.

\item Calculate the weights of the ramification divisor $\mathcal{O}(\mathcal{R}_L)$:
$$
(\chi_{u_1},\hdots,\chi_{u_s})
$$

\item Apply $(\tau_\mu^*)^{-k_0d}$ to  $(\chi_{u_1},\hdots,\chi_{u_s})\otimes D(V)$ using (\ref{dv}) or (\ref{inverse-dv}).

\item For each simple factor $L_k$, restrict the output of step (7) to obtain weights $(\mu_1,\hdots,\mu_s)$.
The generator $\mathcal{M}_k$ induces to
$$
- \frac{1}{d_k} \op{Ind}(\mu_1,\hdots,\mu_s),
$$
and $\op{Ind}$ is the composition $p\circ (\tau_\mu^*)^{k_0d}\circ \op{ext}$ (use (\ref{tau-mu-weights}) and (\ref{twist-and-correct})).

\item For each simple factor $L_k$, now apply $\op{Ind}$ to the extremal rays of $\pic^+(\Parbun_{L_k})$.

\end{enumerate}

\subsection{Shortcuts for $\nu$ and $d_k$}

\begin{proposition}\label{nu-restrict}
Fix a standard Levi subgroup $L$ and a simple component $L_1$ of $L$. Let $\mu\in Q^\vee$ be such that $\alpha \in R_{L_1}^+\implies \alpha(\mu) = 0$ or $1$. Let $\theta_1\in R_{L_1}^+$ denote the highest root for $L_1$, and set
$$
\nu =  \sum_{\alpha\in R^+\setminus R_L^+} \alpha(\mu)\alpha.
$$
Then for any simple root $\alpha_j \in R_1^+$, we have
$$
\nu(\alpha_j^\vee) = \left\{
\begin{array}{ll}
\displaystyle \frac{1}{2} \sum_{\alpha\in R^+\setminus R_L^+} \alpha(\theta_1^\vee)\alpha(\theta_1^\vee), & \alpha_j(\mu) = 1\\
0, & \text{else}
\end{array}
\right.
$$
\end{proposition}

\begin{proof}
Recall from \cite[Lemma 3.4]{BKq} that
$$
\sum_{\alpha\in R^+} \alpha(x)\alpha(y) = (1+\rho(\theta^\vee))(x,y),
$$
where $(\cdot,\cdot)$ is the Killing form for $G$ normalized so that $(\theta,\theta)=2$.

In like manner, for $x,y$ belonging to $\mathfrak{h}_{L_1}$,  we have
$$
\sum_{\alpha\in R_{L_1}^+} \alpha(x)\alpha(y) = (1+\rho_1(\theta_1^\vee))(x,y)_1,
$$
where $\rho_1$ denotes the half-sum of positive roots for $L_1$ and $(\cdot,\cdot)_1$ is the Killing form for $L_1$, likewise normalized so that $(\theta_1,\theta_1)_1 = 2$.

Combining these two formulas, and observing that roots of any other simple components of $L$ will pair trivially on $\mathfrak{h}_{L_1}$, we obtain, for any $x\in \mathfrak{h}_{L_1}$,
$$
\nu(x) = (1+\rho(\theta^\vee))(\mu, x) - (1+\rho_1(\theta_1^\vee))(\mu_1,x)_1,
$$
where $\mu_1$ is any projection of $\mu$ onto $\mathfrak{h}_{L_1}$ such that $\alpha(\mu) = \alpha(\mu_1)$ for all $\alpha\in R_{L_1}$. Therefore if $\alpha_j(\mu)=0$, we have $(\mu,\alpha_j^\vee)=(\mu_1,\alpha_j^\vee)_1=0$ and thus
$$
\nu(\alpha_j^\vee) = (1+\rho(\theta^\vee))(\mu, \alpha_j^\vee) - (1+\rho_1(\theta_1^\vee))(\mu_1,\alpha_j^\vee)_1 = 0.
$$

Now, we can write $\theta_1 = \sum_{\alpha_j \in R_{L_1}^+} c_j \alpha_j$, with each $c_j\in \Z_{>0}$. Moreover,
$$
\theta_1(\mu) = \sum_{\alpha_j \in R_{L_1}^+} c_j \alpha_j(\mu).
$$
As $\theta_1(\mu)$ is constrained to be $0$ or $1$, there exists at most one $\alpha_j$ such that $\alpha_j(\mu)=1$, for which $c_j$ must also equal $1$, and $\alpha_k(\mu)=0$ for all other simple roots $\alpha_k$ belonging to $R_{L_1}$.
Supposing now that such an $\alpha_j$ does exist having $\alpha_j(\mu)=1$, it is apparent that
$$
\alpha_j(\mu) = \theta_1(\mu)=1
$$
Due to the coefficient of $1$ on $\alpha_j$ inside $\theta_1$, $\alpha_j$ is a long root and has the same length as $\theta_1$ with respect to both $(\cdot,\cdot)$ and $(\cdot,\cdot)_1$. So we can calculate
$$
(\mu,\alpha_j^\vee) = \frac{2}{(\alpha_j,\alpha_j)} \alpha_j(\mu) = \frac{2}{(\alpha_j,\alpha_j)} = \frac{2}{(\theta_1,\theta_1)}  = \frac{2}{(\theta_1,\theta_1)} \theta_1(\mu) = (\mu,\theta_1^\vee)
$$
and likewise
$$
(\mu_1,\alpha_j^\vee)_1 = \frac{2}{(\alpha_j,\alpha_j)_1} = \frac{2}{(\theta_1,\theta_1)_1}  =1= (\mu_1,\theta_1^\vee)_1.
$$

Putting it all together, we obtain
$$
\nu(\alpha_j^\vee) = (1+\rho(\theta^\vee))\frac{2}{(\theta_1,\theta_1)} - (1+\rho_1(\theta_1^\vee)).
$$
On the other hand,
$$
\frac{1}{2}\sum_{R^+\setminus R_L^+} \alpha(\theta_1^\vee)\alpha(\theta_1^\vee) = \frac{1}{2}(1+\rho(\theta^\vee))(\theta_1^\vee,\theta_1^\vee) - \frac{1}{2}(1+\rho_1(\theta_1^\vee))(\theta_1^\vee,\theta_1^\vee)_1.
$$
Note that $(\theta_1^\vee,\theta_1^\vee)_1=2$ due to the normalization, while
$$
\frac{1}{2}(\theta_1^\vee,\theta_1^\vee)  = \frac{1}{2}\left(\frac{2\theta_1}{(\theta_1,\theta_1)},\frac{2\theta_1}{(\theta_1,\theta_1)}\right) = \frac{2}{(\theta_1,\theta_1)},
$$
completing the proof.
\end{proof}

\begin{corollary}\label{dual-cox}
The Dynkin index of $\mathfrak{g}/\mathfrak{p}$ as an $L_1$-representation is precisely
$$
\frac{1}{2}\sum_{\alpha\in R^+\setminus R_L^+} \alpha(\theta_1^\vee)^2.
$$
The proof showed that this can be calculated more efficiently as
$$
\frac{2g^*}{(\theta_1,\theta_1)} - g_1^*,
$$
where $g^* = 1+\rho(\theta^\vee)$ and $g_1^* = 1+\rho_1(\theta_1^\vee)$ are the dual Coxeter numbers of $G$ and of $L_1$, respectively.
\end{corollary}

\subsection{Observations on $\mu$}

The ``$\mu$'' in \cite[Lemma 7.2]{BKq} is not unique, but it is unique up to a duality.

\begin{lemma}
Let $P$ be a maximal parabolic, $\mu$ satisfying (i) $0 \le \alpha(\mu) \le 1$ on the Levi and (ii) $|\omega_P(\mu)| =1$. Set $w_\mu = w_0^\mu w_0^L$. Then
$$
-w_\mu^{-1} \mu
$$
also satisfies (i) and (ii).
\end{lemma}

\begin{proof}
Paired with $\omega_P$, $w_\mu^{-1}$ has no effect so the sign simply changes on the $\pm1$.

Let $\beta\in R_L^+$. Either $w_\mu(\beta)$ belongs to $R_L^+$ and pairs with $\mu$ to $0$ or $w_\mu(\beta) \in R_L^-$ and $\langle w_\mu \beta, \mu\rangle = -1$. So we have $\beta(-w_\mu^{-1}\mu) \in \{0,1\}$.
\end{proof}

Note that $-w_\mu^{-1}\mu$ is also $-w_0^L \mu$, which can be viewed as the dual coweight for $\mu$ with respect to the Levi root system.

Besides the dual pair $\mu, -w_0^L \mu$, one might wonder whether there are any others that satisfy the requirements (i) and (ii).

\begin{proposition}
Let $P$ be a maximal parabolic, and let $\mu$ satisfy (i) and (ii) of the previous Lemma. Then $\mu$ and $-w_\mu^{-1}\mu$ are the only elements of $Q^\vee$ satisfying (i) and (ii).
\end{proposition}

\begin{proof}
By the duality statement, we can freely choose the sign on $\omega_P(\mu)$ to be $= 1$.

If there existed two such $\mu$’s, say $\mu_1$ and $\mu_2$, then their difference is in the coroot lattice for the Levi, as $\omega_P(\mu_1-\mu_2)=0$.

Moreover, by definition, they are both contained in the fundamental Weyl alcove (for each connected component of the Levi).

But the Weyl alcove is fundamental for the affine Weyl group action, including coroot translations, so no two points in the alcove differ by something in the coroot lattice.
\end{proof}

\section{Examples}\label{exemples}

These examples all take $s=3$. We follow the Bourbaki conventions for labelling Dynkin diagrams. Our quantum cohomology calculations for type $D_4$ were accomplished using A. Buch's Maple program {\tt qcalc}, available at \href{https://sites.math.rutgers.edu/~asbuch/qcalc/qcalc-manual.txt}{https://sites.math.rutgers.edu/$\sim$asbuch/qcalc/qcalc-manual.txt}. For type $G_2$, we used the tables in \cite[\S 5.2]{TW} (as well as the deformed product tables in \cite[\S 9]{BKq}). For deducing the extremal rays of rational cones from their defining inequalities, we used {\tt Normaliz} \cite{nmz}. To calculate dimensions of spaces of conformal blocks, we used the Macaulay2 package {\tt ConformalBlocks} by D. Swinarski. 

\subsection{Some calculations in type $D_4$}

For $G$ of type $D_4$, the cone $\mathcal{C}_G$ has $771$ facets. Of these, $12$ facets are dominant chamber walls and $3$ are the ``alcove walls'' $\lambda_i(\theta^\vee)\le \ell$. The remaining $756$ facets are regular facets with an associated Levi-movable Gromov-Witten invariant equal to $1$. For the purpose of illustration, we will examine just two 
of these regular facets.

For this root system, the highest root is $\theta = \alpha_1+2\alpha_2+\alpha_3+\alpha_4$, and
$$
s_\theta = s_2s_1s_3s_4s_2s_4s_3s_1s_2.
$$

\subsubsection{A regular facet with maximal parabolic index $1$}

The degree $d=1$ deformed Gromov-Witten invariant
$$
\langle \sigma_{s_1}^P,    \sigma_{s_2s_1}^P,     \sigma_{s_4s_2s_1 }^P  \rangle_1^{\circledast_0}
$$
is equal to $1$.  
Here $P = P_1$ is the maximal parabolic whose negative simple roots consist of $\{-\alpha_2,-\alpha_3,-\alpha_4\}$ with Levi Dynkin diagram highlighted in blue:

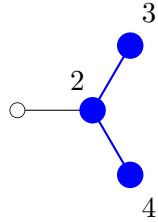
\begin{center}
\begin{tikzpicture}
\draw[-] (4,0) -- (5,0) -- (5.5,0.86);
\draw[-] (5,0) -- (5.5,-0.86);

\draw[blue,thick,-] (5.5,0.86) -- (5,0) -- (5.5,-0.86);

\fill[color=blue] (5,0) circle (5pt);
\fill[color=black] (4,0) circle (3pt);
\node at (4,0) [circle,white,fill,inner sep=1.9pt]{};

\fill[color=blue] (5.5,0.86) circle (5pt);
\fill[color=blue] (5.5,-0.86) circle (5pt);

\node at (4.8,0.4) {$2$};
\node at (5.75,1.3) {$3$};
\node at (5.75,-1.3) {$4$};

\end{tikzpicture}
\end{center}

The four type $I$ rays are listed below, along with the data from which they are derived:

$$
\begin{array}{rrr}
j=1: ~~ & e \xrightarrow{\alpha_1} s_1: & (\omega_1, \omega_4, \omega_3, 1)  \\
  & s_1 < \overline{s_\theta s_1}: 
   &  (0, \omega_3, \omega_3, 1)   \\
j=2: ~~ & s_1 \xrightarrow{\alpha_2} s_2s_1: & (\omega_2, \omega_2, 2\omega_3, 2)  \\
j=3: ~~ & s_2s_1 \xrightarrow{\alpha_4} s_4s_2s_1:  &  ~~~  (\omega_2, \omega_3+\omega_4, \omega_3+\omega_4, 2)
\end{array}
$$~\\

Let us verify, for example, the ray obtained from the choice $(j=1, v = \overline{s_\theta s_1})$. After exchanging $u_1$ with $v$ and shifting the degree via $d' = d - \omega_P(s_1^{-1}\theta^\vee) = 1 - 1 = 0$, we calculate the pushforward of the cycle class
$$
[\Omega(s_1s_2s_3s_4s_2s_1,  ~ s_2s_1,  ~ s_4s_2s_1,  ~ 0)]
$$
in the first position, we note that
$$
\overline{s_\theta s_1} = s_1s_2s_3s_4s_2s_1 = \overline{w_0},
$$
so none of $s_1, \hdots, s_4$ will increase the coset. Thus $\lambda_1 = 0$.

In the second position, only $s_3$ and $s_4$ will increase the coset. One finds that
$$
\langle s_1s_2s_3s_4s_2s_1,  ~ s_3s_2s_1,  ~ s_4s_2s_1 \rangle_0 = 1
$$
and
$$
\langle s_1s_2s_3s_4s_2s_1,  ~ s_4s_2s_1,  ~ s_4s_2s_1 \rangle_0 = 0,
$$
a direct consequence of the fact that the Schubert classes
$$
[X_{s_3s_2s_1}^P] ~~ \text{  and  } ~~ [X_{s_4s_2s_1}^P]
$$
are Poincar\'e dual to one another. Thus $\lambda_2 = \omega_3$.

In the last position, the only possible increase is via $s_3$. As
$$
\langle s_1s_2s_3s_4s_2s_1,  ~ s_2s_1,  ~ s_3s_4s_2s_1 \rangle_0 = 1
$$
(another Poincar\'e duality statement), we find that $\lambda_3 = \omega_3$.

Finally, we calculate the level (Theorem \ref{formula11}(1)). Finding that $\overline{s_\theta w_0} = s_1$ and $\omega_P(\theta^\vee)=1$, we evaluate
$$
\ell  = \langle s_1s_2s_3s_4s_2s_1, ~ s_2s_1, ~ s_4 s_2 s_1, ~ s_1 \rangle_{0+1} = 1.
$$

Putting it all together,
$$
[\Omega(s_1s_2s_3s_4s_2s_1,  ~ s_2s_1,  ~ s_4s_2s_1,  ~ 0)]
 = (0, \omega_3, \omega_3, 1).
 $$


Now let us illustrate the induction algorithm. As $\alpha_1$ is a long root (all roots are the same length in type $D_4$), we may take $\mu = -\alpha_1^\vee$ (see \cite[Lemma 7.2]{BKq}).


This clearly has $\omega_P(\mu) = \omega_1(-\alpha_1^\vee) = -1$, so $k_0 = -1$.

Among the simple roots of $L$, only $\alpha_2$ pairs nontrivially with $\mu$: $\alpha_2(\mu) = 1$. Therefore
$$
w_\mu = (s_3s_4)(s_3s_2s_4s_2s_3s_2) = s_2s_3s_4s_2.
$$

The Dynkin index for the restriction of $\mathfrak{g}/\mathfrak{p}$ to $L$ comes out to $2$:
$$
\begin{array}{r|c}
\alpha \in R^+\setminus R_{\mathfrak{l}}^+ &  \alpha(\alpha_2^\vee+\alpha_3^\vee +\alpha_4^\vee) \\\hline
\alpha_1 & -1 \\
\alpha_1+\alpha_2 & -1 \\
\alpha_1+\alpha_2+\alpha_3 & 0 \\
\alpha_1+\alpha_2+\alpha_4 & 0 \\
\alpha_1+\alpha_2+\alpha_3+\alpha_4 & 1 \\
\alpha_1+2\alpha_2+\alpha_3+\alpha_4 & 1 \\
\end{array}
$$
$$
\frac{1}{2}\left((-1)^2+ (-1)^2+ (0)^2+ (0)^2+ (1)^2+ (1)^2 \right) = 2
$$
Alternatively, we could use Corollary \ref{dual-cox}:

$$
\frac{2 g^*}{(\theta_L, \theta_L)} - g_L^* = \frac{2\cdot 5}{2} - 3 = 2
$$

This directly leads to the calculation of $\nu|_{L'}$:
$$
\nu|_{L'} = 2\omega_2^L
$$

Next we calculate the (restrictions of the) $\chi_w$'s:
\begin{align*}
\chi_{s_1} &= \rho - 2\rho^L + s_1 \rho  \xrightarrow{|_{L'}} \omega_2^L\\
\chi_{s_2s_1}& = \rho - 2\rho^L + s_1s_2 \rho \xrightarrow{|_{L'}} \omega_3^L + \omega_4^L \\
\chi_{s_4s_2s_1}& = \rho - 2 \rho^L + s_1s_2s_4 \rho \xrightarrow{|_{L'}} 2\omega_3^L
\end{align*}

We must apply $(\tau_\mu^*)^{-k_0d} = \tau_\mu^*$ to $(\chi_{s_1}, \chi_{s_2s_1} , \chi_{s_4s_2s_1} )\otimes D(V)$ and then restrict to $L'$, resulting in:
\begin{align*}
(w_\mu \omega_2^L+\nu|_{L'}, \omega_3^L + \omega_4^L, 2 \omega_3^L) \otimes \mathcal{M}^2 &= (-\omega_2^L+2\omega_2^L, \omega_3^L + \omega_4^L, 2 \omega_3^L) \otimes \mathcal{M}^2 \\
&= (\omega_2^L, \omega_3^L + \omega_4^L, 2 \omega_3^L) \otimes \mathcal{M}^2 \\
\end{align*}

As the latter will induce to $0$, we can find $\op{Ind}(\mathcal{M})$ by first applying $\op{Ind}$ to the finite part
$$
(\omega_2^L, \omega_3^L + \omega_4^L, 2 \omega_3^L).
$$

Recall that $\op{Ind}$ is the composition $p\circ (\tau_\mu^*)^{-1} \circ \op{ext}$. One can check that
\begin{align*}
(\omega_2 -  \omega_1)(x_1) &= 0 \\
(\omega_3 - \frac{1}{2} \omega_1)(x_1) &=0 \\
(\omega_4 - \frac{1}{2} \omega_1)(x_1) &= 0;
\end{align*}
therefore under $\op{ext}$,
\begin{align*}
\omega_2^L &\mapsto \omega_2 - \omega_1 \\
\omega_3^L &\mapsto \omega_3 - \frac{1}{2}\omega_1 \\
\omega_4^L &\mapsto \omega_4 - \frac{1}{2}\omega_1.
\end{align*}
So
$$
\op{ext}(\omega_2^L, \omega_3^L + \omega_4^L, 2 \omega_3^L) = ( -\omega_1 + \omega_2, -\omega_1 + \omega_3 + \omega_4, -\omega_1+2\omega_3 ).
$$
Applying $(\tau_\mu^*)^{-1}$, which twists only the first entry by $s_2s_3s_4s_2$, we obtain
$$
(\tau_\mu^*)^{-1}\circ \op{ext}(\omega_2^L, \omega_3^L + \omega_4^L, 2 \omega_3^L) = (\omega_1 - \omega_2, -\omega_1 + \omega_3 + \omega_4, -\omega_1+2\omega_3 ).
$$
Finally, we project what we have so far onto $\mathcal{F}_{\op{II}}$ by twisting the three entries by $s_1, s_2s_1$, and $s_4s_2s_1$, respectively, and then making a correction for each type $I$ ray datum (recall $p  = p_2 \circ \vec u $ from (\ref{verdant})).

First the twisting ($\vec u$):
\begin{align*}
s_1(\omega_1-\omega_2) &= -\omega_1 \\
s_2s_1(-\omega_1+\omega_3+\omega_4) &= \omega_2 \\
s_4s_2s_1(-\omega_1+2\omega_3) &= \omega_3+\omega_4
\end{align*}
Now the corrections ($p_2$):
$$
\begin{array}{rrrr}
& \text{type $I$ datum} &~~~ ~~~ \text{coefficient} & \text{correction} \\
j=1: ~~ & e \xrightarrow{\alpha_1} s_1: &  -\omega_1(\alpha_1^\vee)  & +1(\omega_1, \omega_4, \omega_3, 1)  \\
  & s_1 < \overline{s_\theta s_1}:
   & 0-(-\omega_1)(\theta^\vee) & -1 (0, \omega_3, \omega_3, 1)   \\
j=2: ~~ & s_1 \xrightarrow{\alpha_2} s_2s_1: & \omega_2(\alpha_2^\vee)  &  -1 (\omega_2, \omega_2, 2\omega_3, 2)  \\
j=3: ~~ & s_2s_1 \xrightarrow{\alpha_4} s_4s_2s_1:  & (\omega_3+\omega_4)(\alpha_4^\vee)  & ~~~  -1(\omega_2, \omega_3+\omega_4, \omega_3+\omega_4, 2)
\end{array}
$$
for a total of
$$
(-2\omega_2, -2\omega_3, -2\omega_3, -4).
$$
With the Dynkin index of $2$, we conclude that
$$
\op{Ind}(\mathcal{M}) = -\frac{1}{2}(-2\omega_2, -2\omega_3, -2\omega_3, -4) = (\omega_2, \omega_3, \omega_3, 2).
$$
One can double-check that this generates an extremal ray of $\mathcal{F}$.

There are $20$ rays of the type $A_3$ cone $\op{Pic}^+(\op{Parbun}_{L'})$. As one example, we will induce the ray
$$
(\omega_4^L, \omega_3^L, 0, 1)
$$
Starting with the finite part,
\begin{align*}
(\omega_4^L, \omega_3^L, 0) & \xrightarrow{\op{ext}} (-\frac{1}{2}\omega_1 + \omega_4, -\frac{1}{2}\omega_1 + \omega_3, 0)  \xrightarrow{(\tau_\mu^*)^{-1}} (\frac{1}{2}\omega_1-\omega_2+\omega_3, -\frac{1}{2} \omega_1 + \omega_3,  0 ) \\
& \xrightarrow{\vec u}
(-\frac{1}{2} \omega_1 - \frac{1}{2} \omega_2 + \omega_3, ~~ \frac{1}{2} \omega_2 + \frac{1}{2} \omega_3 - \frac{1}{2} \omega_4, ~~ 0  ) \\
&\xrightarrow{p_2}
(-\omega_2+\omega_3, 0, -\omega_3, -1)
\end{align*}
Add this to $\op{Ind}(\mathcal{M})$ to obtain
$$
\op{Ind}(\omega_4^L, \omega_3^L, 0, 1) = (\omega_3, \omega_3, 0, 1),
$$
which is clearly an extremal ray of $\mathcal{C}_G$.

We should point out that $\op{Ind}$ of an extremal ray need not be extremal (this was also the case in \cite{BRays,BKiers,KBranch}). For instance, one can check that
$$
\op{Ind}(0,\omega_4^L, \omega_3^L,1) = (0,0,0,0).
$$

It can even happen that $\op{Ind}$ returns a nonzero, but non-extremal, vector:
$$
\op{Ind}(\omega_3^L+\omega_4^L,\omega_2^L, \omega_2^L,2) = (\omega_2+\omega_3+\omega_4, \omega_1+\omega_3+\omega_4, \omega_1+2\omega_3, 4),
$$
which is the sum of nonparallel elements of $\mathcal{C}_G$:
$$
(\omega_3,0,\omega_3,1) + (\omega_4,\omega_1,\omega_3,1) + (\omega_2,\omega_3+\omega_4,\omega_1,2).
$$



\subsubsection{A regular facet with maximal parabolic index $2$}

There are some minor additional considerations when the Levi has a reducible root system. Let $P = P_2$, the maximal parabolic whose negative simple roots are $\{-\alpha_1,-\alpha_3,-\alpha_4\}$ with Levi Dynkin diagram highlighted in blue:

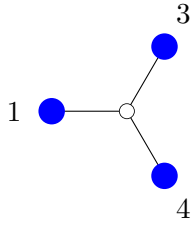
\begin{center}
\begin{tikzpicture}
\draw[-] (4,0) -- (5,0) -- (5.5,0.86);
\draw[-] (5,0) -- (5.5,-0.86);

\fill[color=blue] (4,0) circle (5pt);
\fill[color=black] (5,0) circle (3pt);
\node at (5,0) [circle,white,fill,inner sep=1.9pt]{};

\fill[color=blue] (5.5,0.86) circle (5pt);
\fill[color=blue] (5.5,-0.86) circle (5pt);

\node at (3.5,0) {$1$};
\node at (5.75,1.3) {$3$};
\node at (5.75,-1.3) {$4$};

\end{tikzpicture}
\end{center}

The degree $d=2$ deformed Gromov-Witten invariant
$$
\langle \sigma_{s_2}^P,    \sigma_{s_3s_1s_2}^P,     \sigma_{s_4s_3s_1s_2}^P  \rangle_2^{\circledast_0}
$$
is equal to $1$.

There are $9$ type $I$ rays, listed below:
$$
\begin{array}{rrr}
j=1: ~~ & e \xrightarrow{\alpha_2} s_2: & (\omega_2, 2\omega_4, \omega_2, 2)  \\
  & s_2 < \overline{s_\theta s_2}: 
   &  (\omega_1+\omega_3+\omega_4, \omega_2+2\omega_4, 2\omega_2, 4)   \\
j=2: ~~ & s_3s_2 \xrightarrow{\alpha_1} s_3s_1s_2 : & (\omega_1+\omega_4, \omega_1+\omega_4, \omega_2, 2)  \\
     & s_1s_2 \xrightarrow{\alpha_3} s_3s_1s_2 : & (\omega_3+\omega_4, \omega_3+\omega_4, \omega_2,2)  \\
       & s_3s_1s_2 < \overline{s_\theta s_3s_1s_2} : & (\omega_1+\omega_3, \omega_4, \omega_2, 2)  \\
j=3: ~~ & s_4s_3s_2 \xrightarrow{\alpha_1} s_4s_3s_1s_2:  &  ~~~  (\omega_3, \omega_4, \omega_1, 1) \\
            & s_4s_1s_2 \xrightarrow{\alpha_3} s_4s_3s_1s_2:  &  ~~~  (\omega_1, \omega_4, \omega_3, 1) \\
            & s_3s_1s_2 \xrightarrow{\alpha_4} s_4s_3s_1s_2:  &  ~~~  (\omega_1+\omega_3+\omega_4, \omega_2+\omega_4, \omega_2+\omega_4, 3) \\
            &  s_4s_3s_1s_2 < \overline{s_\theta s_4s_3s_1s_2}:  &  ~~~  (\omega_4, \omega_4, 0, 1) \\
\end{array}
$$~\\

The derived subgroup of the Levi breaks up into three simply-connected simple groups of type $A_1$; that is,
$$
L' = [L,L] \simeq \underbrace{SL_2}_{\alpha_1} \times \underbrace{SL_2}_{\alpha_3} \times \underbrace{SL_2}_{\alpha_4}.
$$

By symmetry, the Dynkin indices of $\mathfrak{g}/\mathfrak{p}$ for each $SL_2$ factor are the same and come out to $4$. This also implies that $\nu|_{L'} = 4\omega_1^L + 4\omega_3^L + 4\omega_4^L$.

Once again we may choose $\mu = -\alpha_2^\vee$, which yields $k_0 = -1$, and one finds that
$
w_\mu = s_1s_3s_4.
$

The class of $\mathcal{R}_{L}$, shifted to degree $0$ and restricted to $L'$, corresponds to
\begin{align*}
(\tau_\mu^*)^2 &\left((\chi_{s_2},\chi_{s_3s_1s_2}, \chi_{s_4s_3s_1s_2})\otimes D(V) \right) |_{L'} \\
&= (\omega_1^L+ \omega_3^L + \omega_4^L, \omega_1^L+ \omega_3^L + 3\omega_4^L, 2\omega_1^L+ 2\omega_3^L + 2\omega_4^L)\otimes \mathcal{M}_1^4 \otimes \mathcal{M}_2^4 \otimes \mathcal{M}_3^4.
\end{align*}

Crucially, the restriction of this class to each simple factor of $L'$ will induce to $0$. So for example,
$$
(\omega_1^L, \omega_1^L, 2\omega_1^L)\otimes \mathcal{M}_1^4
$$
will itself induce to $0$, allowing one to solve for $\op{Ind}(\mathcal{M}_1)$.

In this particular example, it so happens that
$$
\op{Ind}(\mathcal{M}_1) = \op{Ind}(\mathcal{M}_2) = \op{Ind}(\mathcal{M}_3) = 0.
$$

The type $A_1\times A_1 \times A_1$ cone $\op{Pic}^+(\op{Parbun}_{L'})$ has $12$ extremal rays. We list them below together with their images under $\op{Ind}$. It happens in this example that each nonzero entry on the right is an extremal ray of $\mathcal{F}_{\op{II}}$.

$$
\begin{array}{r|r}
\vec r & \op{Ind}(\vec r)  \\ \hline
 \mathcal{M}_1 & 0 \\
(\omega_1^L, \omega_1^L, 0) \otimes \mathcal{M}_1 & (2 \omega_1, \omega_2, \omega_2, 2) \\
(0, \omega_1^L, \omega_1^L)\otimes \mathcal{M}_1 & 0 \\
(\omega_1^L, 0,  \omega_1^L)\otimes \mathcal{M}_1 & 0 \\
\mathcal{M}_2 & 0 \\
(\omega_3^L, \omega_3^L, 0) \otimes \mathcal{M}_2 & (2\omega_3, \omega_2, \omega_2, 2)  \\
(0, \omega_3^L, \omega_3^L) \otimes \mathcal{M}_2 & 0 \\
(\omega_3^L, 0, \omega_3^L) \otimes \mathcal{M}_2 & 0 \\
\mathcal{M}_3 & 0 \\
(\omega_4^L, \omega_4^L, 0) \otimes \mathcal{M}_3 & 0 \\
(0, \omega_4^L, \omega_4^L) \otimes \mathcal{M}_3 & 0 \\
(\omega_4^L, 0, \omega_4^L) \otimes \mathcal{M}_3 &  (2\omega_4, \omega_2, \omega_2, 2)
\end{array}
$$



\subsubsection{Failure of saturation in a simply-laced root system}

The Saturation Conjecture of \cite{KapMill} concerns the tensor product decomposition problem for a complex reductive Lie group $G$. Specifically, they conjectured that if $G$ is simple of simply-laced type, then for dominant weights $\lambda, \mu, \nu$ such that $\lambda+\mu+\nu$ is in the root lattice,
\begin{align}\label{sat-conj}
\left(V_\lambda \otimes V_\mu \otimes V_\nu \right)^G \ne 0 \iff \left(V_{N\lambda} \otimes V_{N\mu} \otimes V_{N\nu}\right)^G \ne 0
\end{align}
for some positive integer $N$. Clearly $(\Leftarrow)$ is the interesting direction. If $G$ is of type $A$, this conjecture was already proven by Knuston and Tao \cite{KT}. It has been checked by computer in types $D_4$, \cite{KKM} $D_5, D_6,$ \cite{KSat} and $E_6$ \cite{Kthesis}.

One can ask a similar saturation question in the setting of conformal blocks. Recall
$$
V(\lambda, \mu, \nu,\ell) = H^0(\op{Parbun}_G,\mathcal{B}(\lambda,\mu, \nu, \ell))^*.
$$
Then the question is whether the following implication holds, for simple and simply-laced $G$:
\begin{align}\label{sat-conj-2}
V(\lambda, \mu, \nu,\ell) \ne 0 \iff V(N\lambda, N\mu, N\nu, N\ell) \ne 0
\end{align}
for some $N>0$. Once again in type $A$, this has been proven \cite{BQHorn}. However, in the course of our investigations into the cone $\mathcal{C}_{\op{Spin}(8)}$, we found a counterexample to (\ref{sat-conj-2}).

\begin{proposition}\label{non-sat}
The na\"ive generalization of the Saturation Conjecture to the conformal blocks setting (\ref{sat-conj-2}) fails in type $D_4$.
\end{proposition}

\begin{proof}
One can check that
$$
V(\omega_2, \omega_2, \omega_2,2) = 0
$$
while
$$
V(2\omega_2, 2\omega_2, 2\omega_2,4) \ne 0,
$$
and $3\omega_2$ belongs to the root lattice. 

In fact, we have checked,  for $n\le 10$, that 
$$
V(\omega_2, \omega_2, \omega_2,2) = 0
$$
\end{proof}

The vector $(\omega_2, \omega_2, \omega_2,2)\in \mathcal{C}_{\op{Spin}(8)}$ is an extremal ray but is not on any regular facet. This is a new phenomenon compared to the saturated tensor cone. See Section \ref{non-regular}.

\subsection{Some calculations in type $G_2$}

For $G$ of type $G_2$, the cone $\mathcal{C}_G$ has $48$ facets. Of these, $6$ facets are dominant chamber walls and $3$ are the ``alcove walls'' $\lambda_i(\theta^\vee)\le \ell$. The remaining $39$ facets are regular facets with an associated Levi-movable Gromov-Witten invariant equal to $1$. For the purpose of illustration, we will examine just one 
of these regular facets.

For this root system, the highest root is $\theta = 3\alpha_1 + 2\alpha_2$, and
$$
s_\theta = s_2s_1s_2s_1s_2.
$$

\subsubsection{A regular facet with maximal parabolic index $2$}


The degree $d=1$ deformed Gromov-Witten invariant
$$
\langle \sigma_{s_1s_2s_1s_2}^P,    \sigma_{s_1s_2}^P,     \sigma_{s_2 }^P  \rangle_1^{\circledast_0}
$$
is equal to $1$.  
Here $P = P_2$ is the maximal parabolic whose negative simple root is just $-\alpha_1$. The Levi Dynkin diagram is of course an $A_1$ diagram:

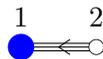
\begin{center}
\begin{tikzpicture}

\draw[-] (0,0) -- (1,0);
\draw[-] (0, 0.05) -- (1, 0.05);
\draw[-] (0, -0.05) -- (1, -0.05);

\draw[-] (0.65, 0.1) -- (0.5, 0) -- (0.65, -0.1);

\fill[color=black] (1,0) circle (3pt);
\node at (1,0) [circle,white,fill,inner sep=1.9pt]{};

\fill[color=blue] (0,0) circle (5pt);

\node at (0,0.4) {$1$};
\node at (1,0.4) {$2$};

\end{tikzpicture}
\end{center}

The five type $I$ rays are listed below, along with the data from which they are derived:

$$
\begin{array}{rrr}
j=1: ~~ & s_2s_1s_2 \xrightarrow{\alpha_1} s_1s_2s_1s_2: & (\omega_1, \omega_2, 2\omega_1, 2)  \\
j=2: ~~ & s_2 \xrightarrow{\alpha_1} s_1s_2: & (0, \omega_1, \omega_1, 1)  \\
     &   s_1s_2 < \overline{s_\theta s_1s_2}: & (\omega_2, \omega_2, 3\omega_1, 3)  \\
j=3: ~~ & e \xrightarrow{\alpha_2} s_2 :  &  ~~~  (0, \omega_2, \omega_2, 2) \\
     & s_2 < \overline{s_\theta s_2}: &  (\omega_2, 2\omega_2, 3\omega_3, 4)
\end{array}
$$~\\

The Dynkin index of $\mathfrak{g}/\mathfrak{p}$ for $SL_{2,\alpha_1}$ is $10$.

The standard choice for $\mu$ is $-\alpha_2^\vee$, for which $k_0 = \omega_2(-\alpha_2^\vee) = -1$ and $w_\mu = s_1$.

The class of $\mathcal{R}_L$, shifted to degree $0$ and restricted to $L'$, corresponds to
$$
(7 \omega_1^L, 4\omega_1^L, 3\omega_1^L) \otimes \mathcal{M}.
$$
The finite part induces to $0$, so we have
$$
\op{Ind}(\mathcal{M}) = 0.
$$

There are $4$ rays of $\op{Pic}^+(\op{Parbun}_{L'})$, and the following table records where they map to under $\op{Ind}$.

$$
\begin{array}{r|r}
\vec r & \op{Ind}(\vec r)  \\ \hline
\mathcal{M} & 0  \\
(\omega_1^L, \omega_1^L, 0) \otimes \mathcal{M} &  0  \\
(0, \omega_1^L, \omega_1^L) \otimes \mathcal{M} &  (\omega_2, \omega_2, 2\omega_1 ,2)  \\
(\omega_1^L, 0,  \omega_1^L) \otimes \mathcal{M} &  0
\end{array}
$$




\subsubsection{Extremal rays not on any regular facet}\label{non-regular}

The extremal rays of the saturated tensor cones (studied in \cite{BRays, BKiers}) have the property that each belongs to a regular face. This is also true of the present cones $\mathcal{C}_G$ when $G$ is of type $A$ \cite{BRigid}. However, beyond these assumptions, the property ceases to hold. In other words, it is possible for an extremal ray to be not on any regular face. This means it is an extremal ray of the dominant chamber itself.

\begin{proposition}
Suppose $\vec r$ is an extremal ray of $\mathcal{C}_G$ with $s=3$ but does not belong to any regular face. Then $\vec r$ has the form
$$(c_i\omega_i, c_j\omega_j, c_k\omega_k, \ell)$$
where
$$
c_i\omega_i(\theta^\vee) = c_j\omega_j(\theta^\vee) = c_k\omega_k(\theta^\vee) = \ell.
$$
\end{proposition}

\begin{proof}
The rays of the dominant chamber are generally of the form 
$$
(\lambda_1, \hdots, \lambda_s, \ell)
$$
where each $\lambda_i$ satisfies either
\begin{enumerate}
\item[($a$)] $\lambda_i$ is a rational multiple of a fundamental weight and $\lambda_i(\theta^\vee) = \ell$, or

\item[($b$)] $\lambda_i = 0$.
\end{enumerate}

The special ray $(0,0,0,1)$ is on every classical face.

There are no rays of $\mathcal{C}_G$ of the form $(\omega_i, 0, 0, \ell)$ (or permutations) since there are no $G$-invariants in $V_{\omega_i}$ even after scaling.

Rays of the form $(c_i \omega_i, c_j \omega_j, 0, \ell)$ require that $c_i\omega_i$ and $c_j \omega_j$ are dual weights to one another.  In particular, $c_i = c_j$. Moreover, dual pairs such as these appear on the classical faces. 

This leaves rays of the form listed above.
\end{proof}

Examples of rays not on any regular face include:
\begin{itemize}
\item $(\omega_2, \omega_2, \omega_2, 2)$ in type $D_4$

\item $(\omega_2, \omega_2, \omega_2, 2)$ in type $G_2$
\end{itemize}

\vspace{0.05 in}

\noindent
P.B.: Department of Mathematics, University of North Carolina, Chapel Hill, NC 27599\\
{email: belkale@email.unc.edu}

\noindent
J.K.: Department of Mathematics, Marian University, Indianapolis, IN 46222\\
{email: jkiers@marian.edu}

\end{document}